\documentclass[11pt]{article}
\usepackage{amsmath,enumerate,amsfonts,amssymb}

\def\RR{{\mathbb R}}

\def\Sphere{{\mathbb S}}
\def\CL{{\mycal{L}}}

\newtheorem{theorem} {\sc  Theorem\rm} [section]
\newtheorem{thm}{\sc  Theorem\rm}

\newtheorem{corollary}[theorem] {\sc  Corollary\rm}
\newtheorem{lemma} [theorem] {\sc  Lemma\rm}
\newtheorem{proposition} [theorem] {\sc  Proposition\rm}

\newtheorem{remark}[theorem]{\sc  Remark\rm}

\def\bproof{\noindent{\bf Proof.\;}}
\def\eproof{\hfill$\square$\medskip}

\newcounter{marnote}

\DeclareFontFamily{OT1}{rsfs}{}
\DeclareFontShape{OT1}{rsfs}{m}{n}{ <-7> rsfs5 <7-10> rsfs7 <10-> rsfs10}{}
\DeclareMathAlphabet{\mycal}{OT1}{rsfs}{m}{n}

\def\dist{{\rm dist}}
\def\tr{{\rm tr}}
\def\IdMa{{\rm Id}}
\def\Sing{{\rm Sing}}
\def\mcS{{\mycal{S}}}
\def\supp{{\rm Supp}\,}

\begin{document}

\title{Refined approximation for a class of Landau-de Gennes energy minimizers}
\author{Luc Nguyen \thanks{OxPDE, Mathematical Institute, University of Oxford, 24--29 St Giles', Oxford OX1 3LB, UK; email: luc.nguyen@maths.ox.ac.uk}~ and Arghir Zarnescu \thanks{OxPDE, Mathematical Institute, University of Oxford, 24--29 St Giles', Oxford OX1 3LB, UK; email: zarnescu@maths.ox.ac.uk}}
\maketitle

\begin{abstract}
We study a class of Landau-de Gennes energy functionals in the asymptotic regime of small elastic constant $L>0$. We revisit and sharpen the results in [18] on the convergence to the limit Oseen-Frank functional. We examine how the Landau-de Gennes global minimizers are approximated by the Oseen-Frank ones by determining the first order term in their asymptotic expansion as $L\to 0$. We identify the appropriate functional setting in which the asymptotic expansion holds, the sharp rate of convergence to the limit and determine the equation for the first order term. We find that the equation has a ``normal component'' given by an algebraic relation and a ``tangential component'' given by a linear system.
\end{abstract}



\section{Introduction}

The complexity of nematic liquid crystals has led to the existence of several major competing theories attempting to describe them, ranging from simple and popular theories such as the Oseen-Frank theory \cite{of} to more involved ones such as the Landau-de Gennes theory \cite{dg}. Because of intertwining features of these two theories, it is of considerable interest to consider the asymptotics of the Landau-de Gennes theory to the Oseen-Frank limit.

The difference between different theories exhibits itself in many ways. For our purpose, we will restrict our discussion within the scope of the Oseen-Frank and the Landau-de Gennes theories. For other theories, we refer the readers to the stimulating review by Lin and Liu \cite{Lin-Liu}.

First and foremost of all, different theories use different mathematical descriptions of the molecular alignment. For example, to describe a certain liquid crystal contained in a region $\Omega$ $\subset$ $\RR^d$, $d$ $=$ $2$, $3$, the Oseen-Frank theory uses vector fields $n:\Omega\to\mathbb{S}^2$, whereas the Landau-de Gennes theory uses matrix-valued functions $Q:\Omega\to \mcS_0$, referred to as $Q$-tensors, where $\mcS_0$ denotes the set of three-by-three symmetric and traceless matrices.

The different number of degrees of freedom used in the different theories reflects in the different predictive capacities of the theories. The simplest and very popular Ossen-Frank description, which uses only two degrees of freedom, has the advantage of being simple but misses certain interesting features of liquid crystal. The first deficiency is that it ignores the important ``head-to-tail'' symmetry of the material. The consequences of this deficiency of the Oseen-Frank theory are analysed in \cite{BallZar}.     Moreover, the Oseen-Frank theory can only explain some of the so-called ``defect'' patterns present in nematic  liquid crystals, namely the ``point'' defects but not the more complicated ``line'' or ``wall'' defects. Part of this limitation of the Oseen-Frank theory could be due to the fact that it only considers uniaxial nematic liquid crystals. The more complex Landau-de Gennes theory allows for both uniaxial and biaxial nematic liquid crystals and has five degrees of freedoms hence it could have better predictive capacities.

The additional complexity of Landau-de Gennes theory can, potentially, change dramatically the interpretation of defects: while in the Oseen-Frank theory the defects are discontinuities of the vector fields, in Landau-de Gennes theory one could interpret defects as discontinuities of eigenvectors.  Although this possibility was considered in the mathematical literature [16], [18] and is consistent with the view of P.G. de Gennes [6], there does not seem to be yet a generally accepted definition of defects, and the relevance of this proposed definition remains to be explored.

Despite many differences, it is of interest from both mathematical and practical points of view to see how and to what extent the mathematically simpler Oseen-Frank theory can be used to ``approximate'' the Landau-de Gennes theory. In \cite{Ma-Za}, Majumdar and Zarnescu showed that in the Oseen-Frank limit, and under suitable boundary conditions, a global ``Landau-de Gennes energy minimizer'' $Q_L$, which is parametrized by an elastic constant $L$, can be approximated in a suitable sense by a global ``Oseen-Frank energy minimizer'' $Q_*$ for $L$ sufficiently small. However, many interesting features of $Q_L$ are not necessarily captured by $Q_*$, which has only two degrees of freedom. For example, $Q_*$ has only point defects and is uniaxial, and therefore does not reflect the appearance of optical ``line defects'' that was observed in Schlieren texture \cite{dg} regardless whether one interpretes defects as discontinuities of eigenvectors or as the uniaxial-biaxial interface. 

From the foregoing discussion, it is therefore necessary to consider the difference $D_L$ $:=$ $Q_L - Q_*$ which potentially encodes the use of the three additional degrees of freedom. Even though $D_L$ is a small quantity, its presence might nevertheless create significant physical effects, e.g. in the eigenvalues and eigenvectors of $Q_L$. Evidence for this behavior has been observed numerically, see \cite[Fig. 1]{SonKilHes}. From the theoretical point of view, as far as the eigenvectors of $Q_L$ are concerned, the following behavior is possible. Consider $A_*, A_L:[0,1]^3\to \mcS_0$ which are defined by
\[
A_*(x,y,z)=\left(\begin{array}{ccc} 1 & y & 0\\ 
y & 1 & 0\\
0 & 0 & -2
\end{array}\right) \text{ and } A_L=\left(\begin{array}{ccc} 1+ Lx & y & 0\\ 
y & 1-Lx & 0\\
0 & 0 & -2
\end{array}\right).
\]
It is straightforward to show that there exist smooth functions $e_i:[0,1]^3\to \mathbb{S}^2,i=1,2,3$ such that at each point in space, $\{e_1, e_2, e_3\}$ is an orthonormal frame of eigenvectors of $A_*$. In constrast, no such smooth eigenframe exists for $A_L$. Moreover, if one adopts defects as discontinuities of eigenvectors, then $A_L$ has a line defect while $A_*$ does not have any defect, even though it is close to $A_L$. We nevertheless note that, in the above example, neither $A_*$ nor $A_L$ is a minimizer of the relevant functionals.

The purpose of this paper is bifold. Firstly, we revisit and sharpen the convergence result in \cite{Ma-Za}. Secondly, we prove in an suitable setting that $D_L$ is of order $L$ and derive the equations that govern the limit of $\frac{D_L}
{L}$.

From a mathematical point of view, our problem bears strong analogies to the Ginzburg-Landau theory for superconductors, which has been investigated intensively in the literature. In fact, the work of B\'{e}thuel, Brezis and H\'{e}lein \cite{BBH} has provided much guidance in our work. The new complexity and challenge in our problem come mainly from having to deal with ``high dimensional'' objects, $Q$-tensors defined on a three-dimensional domain with values into a five-dimensional space. The geometry of the five dimensional linear space $\mcS_0$ and its appropriate decompositions that take into account the behaviour near the two dimensional ``limit manifold'' (see \eqref{LimitSurf::Def} for a precise definition) are crucial in obtaining certain cancellations, that allow to by-pass the singular character of the limit $L\to 0$. In order to obtain estimates independent of the elastic constant $L$ we need to use in various ways the maximum principle for scalar equations and the choice of the right scalar quantities is a highly nontrivial task that is strongly influenced by the understanding of the appropriate geometry of the $Q$-tensor and its relations with the equations. The use of the physically significant ``spectral {\it scalar} quantities'' (i.e., those that depend only on the spectrum of $Q$) was already recognized in \cite{Ma-Za} but the analysis here goes significantly further by also using certain specific {\it tensorial} quantities $X_L$, $Y_L$ and $Z_L$ defined by \eqref{XDef}, \eqref{YDef} and \eqref{ZDef}. Note that all these quantities are directly related to the matrix minimal polynomial associated to (matrices belonging to) the limit manifold. 

The rest of the paper is organized as follows. In Section \ref{Statements}, we set up the mathematical framework of our problem and state our main results. In Section \ref{GeomDes}, we quickly survey some geometry of the limit manifold $\mcS_*$. In the first part of Section \ref{ProjForm}, we derive the equations that governs the limit of $\frac{D_L}{L}$ provided it exists. In the second part we derive the equations for the orthogonal projection of $Q_L$ onto the limit manifold, if such projection exists. These derivations use crucially the geometry of the limit manifold. In Sections \ref{SEC-C1AlphaConv}, \ref{SEC-CjConv} we prove the $C^{1,\alpha}$ and $C^j$ convergence of $Q_L$ to $Q_*$, which extends the previous convergence results in \cite{Ma-Za}. Sections \ref{SEC-HDEst} and \ref{SEC-Rate} are devoted to prove the convergence of $\frac{D_L}{L}$ in a suitable setting.

\subsection*{Acknowledgements}

The authors wish to thank Professors John M. Ball and Fang-Hua Lin for stimulating discussions and valuable encouragements. This research was supported in part by the EPSRC Science and Innovation award to the Oxford Centre for Nonlinear PDE (EP/E035027/1). 

\section{Notations and main results}\label{Statements}

Let $M_{3 \times 3}^{\rm sym}$ denote the set of real $3\times 3$ symmetric matrices, $\mcS_0$ the set of traceless symmetric $3 \times 3$ matrices, and $\Omega$ an open bounded subset of $\RR^3$ with smooth boundary. Consider, a Landau-de Gennes functional of the form
\begin{equation}
I_L[Q]= \int_\Omega \Big[\frac{L}{2}|\nabla Q|^2 + f_B(Q)\Big]\,dx, \qquad Q \in H^1(\Omega, \mcS_0).
\label{L-DGFunc}
\end{equation}
Here $f_B$ is the bulk energy density that accounts for the bulk effects, $|\nabla Q|^2$ is the elastic energy density that penalizes spatial inhomogeneities and $L$ $>$ $0$ is a material-dependent elastic constant. In this paper, we will contrive ourselves to the case where the bulk energy density is a quartic polynomial of the form
\begin{equation}
f_B(Q) = -\frac{a^2}{2}\tr(Q^2) - \frac{b^2}{3}\tr(Q^3) + \frac{c^2}{4}[\tr(Q^2)]^2,
\label{BulkEDensity}
\end{equation}
where $a^2$, $b^2$ and $c^2$ are material- and temperature-dependent positive constants. It is well-known that this type of bulk energy density is the simplest form that allows multiple local minima and a first order nematic-isotropic phase transition \cite{dg}, \cite{Virga}.

The Euler-Lagrange equation for the function $I_L$ is
\begin{equation}
L\,\Delta Q_L = -a^2\,Q_L - b^2[Q_L^2 - \frac{1}{3}\tr(Q_L^2)\,\IdMa] + c^2\,\tr(Q_L^2)\,Q_L,
\label{L-DG::ELEqn}
\end{equation}
where the term $\frac{1}{3}b^2\,\tr(Q_L^2)\,\IdMa$ is a Lagrange multiplier that accounts for the tracelessness constraint. By standard arguments using elliptic theory, it can be seen that $H^1$ solution of \eqref{L-DG::ELEqn} is real analytic, see e.g. \cite[Proposition 13]{Ma-Za}. 

In \cite{Ma-Za}, it was shown that, subjected to a given suitable boundary value $Q_b$, along a subsequence, the minimizers $Q_L$ of $I_L$ converge uniformly away from a set of discrete points to a minimizer $Q_*$ of the functional
\begin{equation}
I_*[Q] = \int_\Omega |\nabla Q|^2\,dx, \qquad Q \in H^1(\Omega,\mcS_*)
\label{LimitFunc}
\end{equation}
where 
\begin{equation}
\mcS_* = \{Q \in N: f_B(Q) = \min f_B(Q)\}.
\label{LimitSurf::Def}
\end{equation}
Even though \eqref{LimitSurf::Def} gives a mysterious definition for $\mcS_*$, this surface can in fact be represented as (see \cite{Maj})
\begin{equation}
\mcS_{*} = \Big\{Q = s_+(n \otimes n - \frac{1}{3}\IdMa), n \in \Sphere^2\Big\}, s_+ = \frac{b^2 + \sqrt{b^4 + 24 a^2 c^2}}{4c^2}.
\label{LimitSurf::Rep}
\end{equation}
It is readily seen that $\mcS_*$ is isometric (modulo a scaling) to the projective plane $\RR P^2$, which is not orientable. If $\Omega$ is simply connected, it follows from \cite{BallZar} that $Q_*$ $=$ $s_+(n_* \otimes n_* - \frac{1}{3}\IdMa)$ where $n_*$ $\in$ $H^1(\Omega,\Sphere^2)$ minimizes the so-called Oseen-Frank functional,
\begin{equation}
I_{OF}[n] = \int_\Omega |\nabla n|^2\,dx, \qquad n \in H^1(\Omega,\Sphere^2).
\label{O-FFunc}
\end{equation}
Note that such $n_*$ is usually called a minimizing $\Sphere^2$-valued harmonic map. Analogous to $n_*$, $Q_*$ is also an ``$\mcS_*$-valued harmonic map''. It satisfies the Euler-Lagrange equation for $I_*$,
\begin{multline}
\Delta Q_* = -\frac{2}{s_+^2}|\nabla Q_*|^2\,Q_* + \frac{2}{s_+}\Big[\sum_{\alpha = 1}^3 (\nabla_\alpha Q_*)^2 - \frac{1}{3}|\nabla Q_*|^2\,\IdMa\Big]\\
	= -\frac{4}{s_+^2}(Q_* - \frac{1}{6}\,s_+\,\IdMa)\sum_{\alpha=1}^3 (\nabla Q_*)^2.
\label{Limit::ELEqn}
\end{multline}

The set of discrete points away from which $Q_L$ converges to $Q_*$ mentioned above consists precisely of points of discontinuity of $Q_*$. Presumably, these points correspond to point defects of the nematics. In this context, there is existing literature on the location of singularities. See \cite{AlLi} and the references therein.

To state our main results, let $Q_b:$ $\partial\Omega$ $\rightarrow$ $\mcS_*$ be a smooth given boundary data and consider the minimization problems
\begin{align}
&\min\big\{I_L[Q]: Q \in H^1(\Omega,\mcS_0), Q\big|_{\partial\Omega} \equiv Q_b\big\},
\label{MinProb}\\
&\min\big\{I_*[Q]: Q \in H^1(\Omega,\mcS_*), Q\big|_{\partial\Omega} \equiv Q_b\big\}.
\label{MinLimProb}
\end{align}

Our first result sharpens the $C^0$ convergence result proved in \cite{Ma-Za}.

\begin{thm}\label{MainThm1}
Let $\Omega$ be an open bounded subset of $\RR^3$ and $Q_L$ be a minimizer of the minimization problem \eqref{MinProb}. For any sequence $L_k$ $\rightarrow$ $0$, there exists a subsequence $L_{k'}$ such that $Q_{L_{k'}}$ converges strongly in $H^1$-norm to a minimizer $Q_*$ of the minimization problem \eqref{MinLimProb}. Moreover, if $\Sing(Q_*)$ is the singular set of $Q_*$, i.e. the set of points $x$ $\in$ $\bar\Omega$ where $Q_*$ is not smooth in any neighborhood of $x$, then
\begin{align*}
&Q_{L_{k'}} \rightarrow Q_* \text{ in } C^{1,\alpha}_{\rm loc}(\bar\Omega\setminus\Sing(Q_*),\mcS_0), \alpha \in (0,1),\\
&Q_{L_{k'}} \rightarrow Q_* \text{ in } C^j_{\rm loc}(\Omega\setminus \Sing(Q_*),\mcS_0), j \geq 2.
\end{align*}
In addition, for $K$ $\Subset$ $\bar\Omega\setminus\Sing(Q_*)$, there exists $\bar L$ $>$ $0$ such that for $L_{k'}$ $<$ $\bar L$, we can rewrite \eqref{L-DG::ELEqn} as
\[
\Delta Q_{L_{k'}} =  -\frac{4}{s_+^2}(Q_{L_{k'}} - \frac{1}{6}\,\IdMa)\sum_{\alpha = 1}^3(\nabla_\alpha Q_{L_{k'}})^2 + \,R_{L_{k'}} \text{ in } K,
\]
where $R_{L_{k'}}$ satisfies
\[
\big|\nabla^j R_{L_{k'}}(y)\big|
	\leq \frac{C(a^2,b^2,c^2,\Omega,K,Q_b,Q_*,j)}{\dist(y,\partial\Omega)^{j+2}}\,L_{k'} \text{ for }y \in K.
\]
\end{thm}

\begin{remark}
(a) For the sake of clarity, let us highlight that, in Theorem \ref{MainThm1}, the $C^{1,\alpha}$-convergence of $Q_{L_k'}$ to $Q_*$ is up to the boundary.

(b) The bound for $R_{L_{k'}}$ can be slightly improved. See Corollary \ref{ImprRemEst}.
\end{remark}

Having a good convergence of $Q_L$ to $Q_*$, we turn to the question of how to extract information about the other three degrees of freedom that $Q_*$ misses. A reasonable way is to look for some asymptotic expansion of the form $Q_L$ $=$ $Q_* + L\,Q_\bullet + O(L^2)$. Intuitively, such expansion is only possible if $Q_*$ is an isolated solution to the limiting harmonic map problem \eqref{Limit::ELEqn}. It is thus reasonable to assume that the linearized operator $\CL_{Q_*}$ of the harmonic map problem has some kind of bijectivity.

On a different aspect, it should be noted that, in general, one does not expect that $\frac{1}{L}(Q_L - Q_*)$ to converge in the energy space $H^1_0(\Omega, \mcS_0)$. For example, consider the case where $\Omega$ is the unit ball, $a^2$ $=$ $b^2$ $=$ $c^2$ $=$ $1$ and $Q_*$ is the so-called ``hedgehog'', i.e.
\[
Q_* = \frac{3}{2}\Big(\frac{x \otimes x}{|x|^2} - \frac{1}{3}\,\IdMa\Big).
\]
By a direct computation using \eqref{FAE::Id01} in Theorem \ref{MainThm2} below, one finds that, if $Q_\bullet$ exists, then the normal component $A$ of $Q_\bullet$ (with respect to the decomposition $\mcS_0$ $\approx$ $T_{Q_*}\mcS_* \oplus (T_{Q_*} \mcS_*)^\perp_{\mcS_0}$) is given by
\[
A = \frac{81}{20|x|^2}\Big(\frac{x \otimes x}{|x|^2} - \frac{1}{3}\,\IdMa\Big).
\]
This implies in particular that $Q_\bullet$ does not vanish on $\partial\Omega$. As the $\frac{1}{L}(Q_L - Q_*)$'s vanish on $\partial\Omega$, they cannot converge to $Q_\bullet$ in $H^1_0(\Omega)$.

For the purpose of getting a convergence result for $\frac{1}{L}(Q_L - Q_*)$, we contrive ourselves to the case where $Q_*$ is smooth. Even though our assumption on $Q_*$ is very restrictive in the sense that it already bans the appearance of point defects in $Q_*$, it is not too restrictive in the study of higher dimensional defects of nematics. For, in practice, it has been observed that line defects can occur either close to point defects or far from point defects. We prove:

\begin{thm}\label{MainThm2}
Let $\Omega$ be an open bounded subset of $\RR^3$ and $Q_L$ be a minimizer of the minimization problem \eqref{MinProb}. Assume that $Q_{L_k}$ converges strongly in $H^1(\Omega, \mcS_0)$ to a minimizer $Q_*$ $\in$ $H^1(\Omega,\mcS_*)$ of the limit minimization problem \eqref{MinLimProb} for some sequence $L_k$ $\rightarrow$ $0$. Assume in addition that $Q_*$ is smooth and that the linearized operator $\CL_{Q_*}:$ $H^1_0(\Omega,M_{3\times 3})$ $\rightarrow$ $H^{-1}(\Omega,M_{3\times 3})$ of the limit harmonic map problem (see \eqref{LinearizedOp}) is bijective. Then there exists $Q_\bullet$ $\in$ $C^\infty(\Omega,\mcS_0) \cap H^s(\Omega,\mcS_0)$ for any $0$ $<$ $s$ $<$ $1/2$ such that
\begin{align*}
&\frac{1}{L_{k}}(Q_{L_{k}} - Q_*) \rightarrow Q_\bullet \text{ in } H^s(\Omega), 0 < s < 1/2,\\
&\frac{1}{L_{k}}(Q_{L_{k}} - Q_*) \rightarrow Q_\bullet \text{ in } C^j_{\rm loc}(\Omega), j \geq 0.
\end{align*}
Moreover, if we split $Q_\bullet$ $=$ $A + B$ where $A$ belongs to the normal space $(T_{Q_*} \mcS_*)^\perp_{\mcS_0}$ to $\mcS_*$ at $Q_*$ with respect to $\mcS_0$ and $B$ belongs to the tangent space $T_{Q_*} \mcS_*$ to $\mcS_*$ at $Q_*$, then
\begin{enumerate}[(i)]
\item ${A}$ is given by
\begin{equation}
{A}
	= -\frac{2}{b^2s_+^2}\Big[\frac{6}{6a^2 + b^2\,s_+}|\nabla Q_*|^2\big(c^2\,Q_* + \frac{1}{3}b^2\IdMa\big)\big(Q_* - \frac{1}{6}s_+\,\IdMa\big) - \sum_{\alpha = 1}^3(\nabla_\alpha Q_*)^2\Big],
\label{FAE::Id01}
\end{equation}
\item and ${B}$ satisfies in $\Omega$ the equations
\begin{multline}
\Delta {B}
	= \Big[-b^2({B}\,{A} + {A}\,{B}) - \frac{6c^2}{6a^2 + b^2\,s_+}|\nabla Q_*|^2\,{B}\Big]\\
		-\frac{4}{s_+^2}\big[\big(\nabla {B}\big)^\parallel\,\nabla Q_* + \nabla Q_*\,\big(\nabla {B}\big)^\parallel\big]\big(Q_* - \frac{1}{6}s_+\,\IdMa\big) - \big(\Delta {A}\big)^\parallel,
\label{FAE::Id02}
\end{multline}
where $\big(\nabla {B}\big)^\parallel$ and $\big(\Delta {A}\big)^\parallel$ are the tangential components of $\nabla {B}$ and $\Delta {A}$ with respect to the decomposition $\mcS_0$ $\approx$ $T_{Q_*}\mcS_* \oplus (T_{Q_*} \mcS_*)^\perp_{\mcS_0}$, respectively.
\end{enumerate}
\end{thm}

\begin{remark}
As an example, we note that, if $\Omega$ is simply connected and 
\begin{equation}
\|Q_b - p\|_{C^2(\partial\Omega)} \leq \epsilon
\label{intr::SmallData}
\end{equation}
for some fixed point $p$ $\in$ $\mcS_*$ and some sufficiently small $\epsilon$ $>$ $0$, then $\CL_{Q_*}:$ $H^1_0(\Omega,M_{3\times 3})$ $\rightarrow$ $H^{-1}(\Omega,M_{3\times 3})$ is bijective. To see this, observe first that the constant map $p$ is the unique minimizing harmonic map that has constant boundary values $p$. Thus, by the stability of minimizing $\Sphere^2$-harmonic maps \cite{HardtLin} and the lifting of maps in $H^1(\Omega,\mcS_{*})$ to maps in $H^1(\Omega,\Sphere^2)$ \cite{BallZar}, the minimization problem \eqref{MinLimProb} has a unique solution $Q_*$ provided that \eqref{intr::SmallData} holds for some small $\epsilon$. (For expository purposes, we note that by \cite{Strw}, when $\Omega$ is a ball, $Q_*$ is actually the unique weakly harmonic map admitting $Q_b$ as its boundary values for all sufficiently small $\epsilon$.) Furthermore, by classical regularity results for harmonic maps \cite{su}, \cite{SU-Bdry}, $Q_*$ is smooth and satisfies $|Q_* - p| + |\nabla Q_*| + |\nabla^2 Q_*|$ $\leq$ $o(1)$ as $\epsilon$ $\rightarrow$ $0$. On the other hand, the linearized operator corresponding to a constant map is the Laplace operator and hence is bijective. Thus, by \cite[Chapter IV, Theorem 5.17]{Kato}, $\CL_{Q_*}$ is also bijective provided $\epsilon$ is chosen appropriately small. 
\end{remark}

\section{Preliminaries}\label{GeomDes}

We begin with a brief study of the geometry of the limit manifold $\mcS_*$ defined in \eqref{LimitSurf::Def}. Using \eqref{LimitSurf::Rep} and some simple algebraic manipulations, we find the following equivalent definitions for $\mcS_*$.

\begin{lemma}\label{DiffRepLimSurf}
Let $s_+$ $=$ $\frac{b^2 + \sqrt{b^4 + 4a^2\,c^2}}{4c^2}$. A matrix $Q$ belongs to $\mcS_*$ if and only if one of the following holds.e
\begin{enumerate}[(i)]
\item $f_B(Q)$ $=$ $\min\{f_B(R): R \in \mcS_0\}$.
\item $Q$ $=$ $s_+(n \otimes n - \frac{1}{3}\IdMa)$ for some $n$ $\in$ $\Sphere^2$.
\item $Q$ $\in$ $\mcS_0$, $\tr(Q^2)$ $=$ $\frac{2s_+^2}{3}$ and $\tr(Q^3)$ $=$ $\frac{2s_+^3}{9}$.
\item $Q$ $\in$ $\mcS_0$ and its minimal polynomial is $\lambda^2 - \frac{1}{3}s_+\,\lambda - \frac{2}{9}s_+^2$.
\end{enumerate}
\end{lemma}

For a matrix $Q$ in $\mcS_*$, let $T_Q \mcS_*$ denote the tangent space to $\mcS_*$ at $Q$ in $\mcS_0$, and $(T_Q \mcS_*)^\perp_{\mcS_0}$ and $(T_Q \mcS_*)^\perp$ denote the orthogonal complements of $T_Q \mcS_*$ in the tangent spaces at $Q$ to $\mcS_0$, $T_Q \mcS_0$, and to $M_{3 \times 3}^{\rm sym}$, $T_Q M_{3\times 3}^{\rm sym}$, respectively. We will often identify $\mcS_0$ with $T_Q \mcS_0$ and $M_{3 \times 3}^{\rm sym}$ with $T_Q M_{3 \times 3}^{\rm sym}$. We have the following characterization of these spaces.

\begin{lemma}\label{Tangent-NormalSpaces}
For a point $Q$ $\in$ $\mcS_*$, the tangent and normal spaces to $\mcS_*$ at $Q$ are
\begin{align}
T_Q \mcS_*
	&= \{\dot Q \in M_{3 \times 3}^{\rm sym}: \frac{1}{3}s_+\dot Q = \dot Q\,Q + Q\,\dot Q\},\label{TNS::Id01}\\
(T_Q \mcS_*)^\perp_{\mcS_0}
	&=  \{Q^\perp \in \mcS_0: Q^\perp\,Q = Q\,Q^\perp\},\label{TNS::Id02}\\
(T_Q \mcS_*)^\perp
	&=  \{Q^\perp \in M_{3\times 3}^{\rm sym}: Q^\perp\,Q = Q\,Q^\perp\}.\label{TNS::Id03}
\end{align}
Moreover, any matrix $A$ $\in$ $T_Q M_{3\times 3}^{\rm sym}$ $\approx$ $M_{3\times 3}^{\rm sym}$ can be decomposed uniquely as $A$ $=$ $\dot A + A^\perp$ $\in$ $T_Q \mcS_* \oplus (T_Q \mcS_*)^\perp$ by
\begin{align}
A^\perp
	&= -\frac{2}{s_+^2}\big(\frac{1}{3}s_+\,A - Q\,A - A\,Q\big)\big(Q - \frac{1}{6}s_+\,\IdMa\big)\nonumber\\
	&= -\frac{2}{s_+^2}\big(Q - \frac{1}{6}s_+\,\IdMa\big)\big(\frac{1}{3}s_+\,A - Q\,A - A\,Q\big).\label{TNS::Id05}
\end{align}
\end{lemma}

\bproof Set
\[
\mycal{A} = \{\dot Q \in M_{3 \times 3}^{\rm sym}: \frac{1}{3}s_+\dot Q = \dot Q\,Q + Q\,\dot Q\}.
\]
By Lemma \ref{DiffRepLimSurf}(ii), we have $Q$ $=$ $s_+\big(n \otimes n - \frac{1}{3}\IdMa\big)$ and
\[
T_Q \mcS_* = \{n \otimes \dot n + \dot n \otimes n: \dot n \in T_{n}\Sphere^2\}.
\]
It is thus immediate that $T_Q \mcS_*$ $\subset$ $\mycal{A}$. To see the converse, pick $\dot Q$ $\in$ $\mycal{A}$. Then by Lemma \ref{DiffRepLimSurf}(iv),  
\begin{equation}
9Q\,\dot Q\,Q = -2\,s_+^2 \dot Q.
\label{TNS::Eqn01}
\end{equation}
Choose an orthonormal frame $\{n, \hat n, \check n\}$ of $\RR^3$ at $n$ $\in$ $\Sphere^2$ $\subset$ $\RR^3$. Using \eqref{TNS::Eqn01}, we get
\[
Q(\dot Q n) = - \frac{s_+}{3}\dot Q n, Q(\dot Q\hat n) = \frac{2s_+}{3} \dot Q\hat n, \text{ and } Q(\dot Q\check n) = \frac{2s_+}{3}\dot Q\check n.
\]
By the explicit form of $Q$, we infer that $\dot Q n$ $=$ $\alpha \hat n + \beta\check n$, $\dot Q\hat n$ $=$ $\gamma_1 n$, and $\dot Q\check n$ $=$ $\gamma_2 n$. As $\dot Q$ is symmetric, it is necessary that $\alpha$ $=$ $\gamma_1$ and $\beta$ $=$ $\gamma_2$. We then arrive at $\dot Q$ $=$ $n \otimes (\alpha \hat n + \beta\check n) + (\alpha \hat n + \beta\check n) \otimes n$, which shows $\dot Q$ $\in$ $T_Q \mcS_*$. \eqref{TNS::Id01} follows. 

Next, let
\[
\mycal{B} := \{Q^\perp \in \mcS_0: Q^\perp\,Q = Q\,Q^\perp\}.
\]
We have
\[
(T_Q \mcS_*)^\perp_{\mcS_0} = \{\alpha n \otimes n + \beta \hat n \otimes \hat n + \gamma \check n \otimes \check n + \lambda(\hat n \otimes \check n + \check n \otimes \hat n), \alpha + \beta + \gamma = 0\} \subset \mycal{B}.
\]
Now pick $Q^\perp$ $\in$ $\mycal{B}$. Then
\[
Q(Q^\perp n) = \frac{2s_+}{3}Q^\perp n, Q(Q^\perp\hat n) = -\frac{s_+}{3} Q^\perp \hat n, \text{ and } Q(Q^\perp \check n) = -\frac{s_+}{3}Q^\perp \check n.
\]
Putting $Q^\perp n$ $=$ $\alpha n$, $Q^\perp\hat n$ $=$ $\beta \hat n + \lambda \check n$, and $Q^\perp \check n$ $=$ $\lambda \hat n + \gamma\check n$, we reach \eqref{TNS::Id02}. The proof of \eqref{TNS::Id03} is similar to the above and is thus omitted.

Now let $A$ $\in$ $T_Q M_{3 \times 3}^{\rm sym}$ and decompose $A$ $=$ $\dot A + A^\perp$. By \eqref{TNS::Id03} and \eqref{TNS::Eqn01},
\[
9Q\,A\,Q = -2\,s_+^2\dot A + 9A^\perp\,Q^2,
\]
and so
\[
2s_+^2\,A + 9Q\,A\,Q = 2s_+^2A^\perp + 9 A^\perp\,Q^2 = A^\perp(2s_+^2\IdMa + 9Q^2).
\]
Hence, by Lemma \ref{DiffRepLimSurf}(iv),
\begin{align*}
A^\perp
	&= (2s_+^2\,A + 9Q\,A\,Q)(2s_+^2\IdMa + 9Q^2)^{-1},\\
	&= \frac{1}{s_+}(2s_+^2\,A + 9Q\,A\,Q)(3Q + 4s_+\IdMa)^{-1}\\
	&= -\frac{1}{18s_+^3}(2s_+^2\,A + 9Q\,A\,Q)(3Q - 5s_+\IdMa)\\
	&= \frac{1}{9s_+^2}(-3s_+(AQ + QA) + 18 QAQ + 5 s_+^2 A)\\
	&= -\frac{2}{s_+^2}\big(\frac{1}{3}s_+\,A - Q\,A - A\,Q\big)\big(Q - \frac{1}{6}s_+\,\IdMa\big)\\
	&= -\frac{2}{s_+^2}\big(Q - \frac{1}{6}s_+\,\IdMa\big)\big(\frac{1}{3}s_+\,A - Q\,A - A\,Q\big).
\end{align*}
The last assertion follows.
\eproof

The next lemma states some interrelations between the tangent and normal spaces of $\mcS_*$.
\begin{lemma}\label{InterRel}
Let $Q$ be a point in $\mcS_*$. For $X$, $Y$ in $T_Q \mcS_*$ and $Z$, $W$ in $(T_Q \mcS_*)^\perp$, we have $XY + YX$, $ZW + WZ$ $\in$ $(T_Q \mcS_*)^\perp$ and $XZ + ZX$ $\in$ $T_Q \mcS_*$.
\end{lemma}

\bproof That $ZW + WZ$ $\in$ $(T_Q \mcS_*)^\perp$ follows immediately from Lemma \ref{Tangent-NormalSpaces}. To see that $XY + YX$ $\in$ $(T_Q \mcS_*)^\perp$, we note that Lemma \ref{Tangent-NormalSpaces} gives
\[
\frac{1}{3}s_+\,X = XQ + QX \text{ and }\frac{1}{3}s_+\,Y = YQ + QY.
\]
Hence
\[
XYQ = X\Big[\frac{1}{3}s_+\,Y - QY\Big] = \Big[\frac{1}{3}s_+\,X - XQ]Y = QXY.
\]
Similarly, we have $YXQ$ $=$ $QYX$. It then follows from Lemma \ref{Tangent-NormalSpaces} that $XY + YX$ $\in$ $(T_Q \mcS_*)^\perp$.

For the last assertion, we calculate
\begin{multline*}
(XZ + ZX)Q = XQZ + Z\Big[\frac{1}{3}s_+\,X - QX\Big]\\
	= \Big[\frac{1}{3}s_+\,X - QX\Big]Z + Z\Big[\frac{1}{3}s_+\,X - QX\Big] = \frac{1}{3}s_+(XZ + ZX) - Q(XZ + ZX).
\end{multline*}
By Lemma \ref{Tangent-NormalSpaces}, this implies that $XZ + ZX$ $\in$ $T_Q \mcS_*$.
\eproof

\begin{lemma}\label{MysIds}
For any $X$, $Y$ $\in$ $T_Q \mcS_*$ and $Z$ $\in$ $(T_Q \mcS_*)^\perp$, there holds
\begin{align}
\tr((XY + YX)Q)
	&= \frac{s_+}{3}\tr(XY),\label{MysIds::Id1}\\
(XY + YX)Q
	&= -\frac{s_+}{3}(XY + YX) + \tr(XY)Q + \frac{s_+}{3}\tr(XY)\,\IdMa,\label{MysIds::Id2}\\
\Big[\frac{1}{s_+}Q + \frac{1}{3}\IdMa\Big]Z
	&= \Big(\frac{1}{s_+}\tr(QZ) + \frac{1}{3}\tr(Z)\Big)\,\Big[\frac{1}{s_+}Q + \frac{1}{3}\IdMa\Big].\label{MysIds::Id3}
\end{align}
\end{lemma}

\bproof Since $X$ $\in$ $T_Q \mcS_*$, $XQ + QX$ $=$ $\frac{1}{3}s_+\,X$. It follows that
\[
XQY + QXY = \frac{1}{3}s_+\,XY.
\]
Taking trace, we get \eqref{MysIds::Id1}.

We next prove \eqref{MysIds::Id3}. Set
\[
P_1 = \frac{1}{s_+}Q + \frac{1}{3}\IdMa \text{ and } P_2 = -\frac{1}{s_+}Q + \frac{2}{3}\IdMa.
\]
Then $P_1^2$ $=$ $P_1$, $P_2^2$ $=$ $P_2$, $P_1\,P_2$ $=$ $0$ and $P_1 + P_2$ $=$ $\IdMa$. Set
\[
V_i = \{W \in (T_Q \mcS_*)^\perp: P_i\,W = W\}, \qquad i = 1,2.
\]
Then $(T_Q \mcS_*)^\perp$ $=$ $V_1 \oplus V_2$, $V_1$ $\perp$ $V_2$, and $P_i$ $\in$ $V_i$. Write $Q$ $=$ $s_+(n \otimes n - \frac{1}{3}\IdMa)$ for some $n$ $\in$ $\Sphere^2$ and pick an orthonormal basis $\{n, \hat n, \check n\}$ of $\RR^3$ at $n$. It is easy to check that $\hat n \otimes \hat n$, $\check n \otimes \check n$, and $\hat n \otimes \check n + \check n \otimes \hat n$ belong to $V_2$. It follows that $\dim V_2$ $\geq$ $3$. Since $\dim V_1 + \dim V_2$ $=$ $\dim (T_Q \mcS_*)^\perp$ $=$ $4$, and $\dim V_1$ $\geq$ $1$ (as $P_1$ $\in$ $V_1$), we infer that $\dim V_1$ $=$ $1$. Now, for any $Z$ $\in$ $(T_Q \mcS_*)^\perp$, $P_1\,Z$ $\in$ $V_1$ and so
\[
P_1\,Z = k\,P_1 \text{ for some } k \in \RR.
\]
Taking trace we get
\[
\frac{1}{s_+}\tr(QZ) + \frac{1}{3}\tr(Z) = k.
\]
\eqref{MysIds::Id3} is established.

Finally, to get \eqref{MysIds::Id2}, we note that $W$ $=$ $XY + YX - \frac{2}{3}\tr(XY)\,\IdMa$ $\in$ $(T_Q \mcS_*)^\perp$. Applying \eqref{MysIds::Id3} to $Z$ $=$ $W$, we get
\[
(XY + YX)P_1 = \frac{2}{3}\tr(XY)P_1 + \frac{1}{s_+}\tr[(XY + YX)Q]\,P_1.
\]
\eqref{MysIds::Id2} follows immediately from the above identity and \eqref{MysIds::Id1}.
\eproof

\begin{lemma}\label{2ndForm}
The second fundamental form of $\mcS_*$ in $\mcS_0$ is given by
\begin{align*}
\mathrm{II}(X,Y)(Q)
	&= -\frac{2}{s_+^2}\,\tr(X\,Y)\,Q + \frac{1}{s_+}\Big[X\,Y + Y\,X - \frac{2}{3}\tr(X\,Y)\,\IdMa\Big]\\
	&= -\frac{1}{s_+^2}\,(X\,Y + Y\,X)(2Q - \frac{1}{3}s_+\,\IdMa) = -\frac{1}{s_+^2}\,(2Q - \frac{1}{3}s_+\,\IdMa)(X\,Y + Y\,X)
\end{align*}
where $X$ and $Y$ are tangent vectors to $\mcS_*$ at $Q$ $\in$ $\mcS_*$.
\end{lemma}

\bproof Since $SO(3)$ acts transitively on $\mcS_*$ by conjugation, it suffices to verify the conclusion at
\[
Q_0 = \frac{1}{3}\mathrm{diag}(2,-1,-1) \in \mcS_*.
\]
Also, by linearity, it suffices to consider $X$ and $Y$ in a set of basis vectors. For the ease in calculation, we isometrically embed $\mcS_0$ into $M_{3\times 3}$ by the obvious embedding. We parametrize $M_{3\times 3}$ by
\[
A = \left[\begin{array}{ccc}
y_1 & y_2 & y_3\\
y_4 & y_5 & y_6\\
y_7 & y_8 & y_9
\end{array}\right].
\]

To pick tangential vectors at a point $Q_0$, we note that $B^t\,Q_0 + Q_0\,B$ is tangent to $\mcS_*$ at $Q_0$ for any skew-symmetric matrix $B$. Choosing $B$ $=$ $e_1 \otimes e_2 - e_2 \otimes e_1$, we get
\begin{align*}
v_1
	&= \left[\begin{array}{ccc}
-y_2 - y_4 & y_1-y_5 & -y_6\\
y_1-y_5 & y_2 + y_4 & y_3\\
-y_8 & y_7 & 0
\end{array}\right]\\
	&= -(y_2 + y_4)\partial_1 + (y_1 - y_5)\partial_2 - y_6\partial_3 + (y_1 - y_5)\partial_4 + (y_2+ y_4)\partial_5 + y_3\partial_6 - y_8\partial_7 + y_7\partial_8.
\end{align*}
Choosing $B$ $=$ $e_1 \otimes e_3 - e_3 \otimes e_1$, we get
\begin{align*}
v_2
	&= \left[\begin{array}{ccc}
-y_3 - y_7 & -y_8 & y_1 - y_9\\
-y_6 & 0 & y_4\\
y_1 - y_9 & y_2 & y_3 + y_7
\end{array}\right] \\
	&= -(y_3 + y_7)\partial_1 - y_8 \partial_2 + (y_1 - y_9)\partial_3 - y_6\partial_4 + y_4\partial_6 + (y_1 - y_9)\partial_7 + y_2\partial_8 + (y_3 + y_7)\partial_9.
\end{align*}
It is readily seen that $v_1$ and $v_2$ form a local frame for $\mcS_*$ in a neighborhood of $Q_0$.

Let $\bar\nabla$ denote the connection of $M^{3\times 3}$ $\cong$ $\RR^9$. We calculate
\begin{align*}
\bar\nabla_{v_1} v_1
	&= 2(y_1 - y_5)[-\partial_1 + \partial_5] + ...,\\
\bar\nabla_{v_1} v_2
	&= (y_1 - y_5)[\partial_6 + \partial_8] + ...,\\
\bar\nabla_{v_2} v_2
	&= 2(y_1 - y_9)[-\partial_1 + \partial_9] + ...
\end{align*}
where the dots comprise of terms that vanish at $Q_0$. Since $\mathrm{II}(X,Y)$ is the normal component of $\bar\nabla_X Y$, we thus have
\begin{align*}
\mathrm{II}(v_1, v_1)(Q_0)
	&= 2s_+(-\partial_1 + \partial_5) = 2s_+\mathrm{diag}(-1,1,0),\\
\mathrm{II}(v_1, v_2)(Q_0)
	&= s_+(\partial_6 + \partial_8) = s_+\left[\begin{array}{ccc}
0 & 0 & 0\\
0 & 0 & 1\\
0 & 1 & 0
\end{array}\right],\\
\mathrm{II}(v_2, v_2)(Q_0)
	&= 2s_+(-\partial_1 + \partial_9) = 2s_+\mathrm{diag}(-1,0,1).
\end{align*}
It is elementary to check that the assertion holds for $X$, $Y$ $\in$ $\{v_1, v_2\}$ and $Q$ $=$ $Q_0$. The proof is complete.
\eproof

\begin{corollary}\label{HarMapEqn}
A map $Q_0:$ $\RR^n$ $\rightarrow$ $\mcS_*$ is harmonic if and only if one of the following occurs.
\begin{align*}
&\text{(i) $\Delta Q_0$ commutes with $Q_0$.}\\
&\text{(ii) } \Delta Q_0 = -\frac{2}{s_+^2}\,|\nabla Q_0|^2\,Q_0 + \frac{2}{s_+}\Big[\sum_{\alpha=1}^n(\nabla_\alpha Q_0)^2 - \frac{1}{3}|\nabla Q_0|^2\,\IdMa\Big].\\
&\text{(iii) } \Delta Q_0 = -\frac{4}{s_+^2}\,\sum_{\alpha=1}^n(\nabla_\alpha Q_0)^2\,(Q_0 - \frac{1}{6}s_+\,\IdMa),\\
&\text{(iv) } \Delta Q_0 = -\frac{4}{s_+^2}\,(Q_0 - \frac{1}{6}s_+\,\IdMa)\sum_{\alpha=1}^n(\nabla_\alpha Q_0)^2.
\end{align*}
\end{corollary}

\bproof The conclusion follows immediately from well-known forms of harmonic map equations (see e.g. \cite{Helein}) and Lemmas \ref{Tangent-NormalSpaces}, \ref{2ndForm}.
\eproof

\section{Equations for the first order term in the asymptotic expansion and for the projection onto the limit manifold}\label{ProjForm}

The goal of this section is bifold. Let $Q_L$ be a critical point of $I_L$. In the first part, we will assume that $Q_L$ has a formal asymptotic form $Q_L$ $=$ $Q_* + L\,Q_\bullet + o(L)$ for some harmonic map $Q_*$ and derive the equations that govern $Q_\bullet$. In the second part, we derive the equations for the orthogonal projection $Q_L^\sharp$ of $Q_L$ onto the limit manifold $\mcS_*$ provideded that such projection is defined. Roughly speaking, these equations are of the form ``harmonic map plus corrector terms''. This latter part will be useful when we prove uniform $C^{1,\alpha}$ bounds for $Q_L$ up to the boundary.

\begin{proposition}\label{FirstApproxEqn}
Assume that $Q_{L_k}$ $\in$ $C^2(\Omega, \mcS_0)$ is a critical point of $I_{L_k}$, and that as $L_k$ $\rightarrow$ $0$, $Q_{L_k}$ converges on compact subsets of $\Omega$ in $C^2$-norm to $Q_*$ $\in$ $C^2(\Omega,\mcS_*)$ which is a critical point of $I_*$ and $\frac{1}{L_k}(Q_{L_k} - Q_*)$ converges in $C^2$-norm to some $Q_\bullet$ $\in$ $C^2(\Omega, \mcS_0)$. If we write $Q_\bullet$ $=$ ${A} + {B}$ with ${A}$ $\in$ $(T_{Q_*} \mcS_*)^\perp$ and ${B}$ $\in$ $T_{Q_*} \mcS_*$, then \eqref{FAE::Id01} and \eqref{FAE::Id02} hold.
\end{proposition}

\bproof For simplicity, we will drop the subscript $L_k$. We begin with the proof of \eqref{FAE::Id01}. Split
\[
Q = Q_* + L_k\,\hat Q.
\]
Noting that $Q_*$ $\in$ $\mcS_*$, we calculate using \eqref{L-DG::ELEqn} and $2c^2\,s_+^2 = 3a^2 + b^2\,s_+$,
\begin{align}
L_k\,\Delta\hat Q
	&= \Big[b^2\big(\frac{1}{3}s_+\hat Q - Q_*\,\hat Q - \hat Q\,Q_*\big) + \frac{2}{3}b^2\tr(Q_*\,\hat Q)\IdMa + 2c^2\tr(Q_*\,\hat Q)Q_* - \Delta Q_*\Big]\nonumber\\
		&\qquad\qquad + L_k\Big[-b^2\big(\hat Q^2 - \frac{1}{3}\,\tr(\hat Q^2)\,\IdMa\big) + c^2\,\tr(\hat Q^2)Q_* + 2c^2\,\tr(Q_*\,\hat Q)\,\hat Q\Big]\nonumber\\
		&\qquad\qquad + L_k^2\,c^2\,\tr(\hat Q^2)\,\hat Q.
\label{FAE::Eqn01}
\end{align}
Sending $L_k$ $\rightarrow$ $0$, we thereby obtain
\begin{equation}
\frac{1}{3}s_+\,Q_\bullet - Q_*\,Q_\bullet - Q_\bullet\,Q_*
	= -\frac{2}{3}\tr(Q_*\,Q_\bullet)\IdMa - \frac{2c^2}{b^2}\tr(Q_*\,Q_\bullet)Q_* + \frac{1}{b^2}\Delta Q_*.
\label{FAE::Eqn02}
\end{equation}
Multiplying \eqref{FAE::Eqn02} by $Q_*$ to the left, taking trace and using Lemma \ref{DiffRepLimSurf}(iv), we get
\begin{equation}
\tr(Q_*\,Q_\bullet)
	= \frac{3}{6a^2 + b^2\,s_+}\tr(Q_*\,\Delta Q_*).
\label{FAE::Eqn03}
\end{equation}
Substituting this into \eqref{FAE::Eqn02} and using Lemma \ref{Tangent-NormalSpaces}, we hence get
\begin{align}
{A}
	&= -\frac{2}{s_+^2}\big(\frac{1}{3}s_+\,Q_\bullet - Q_*\,Q_\bullet - Q_\bullet\,Q_*\big)\big(Q_* - \frac{1}{6}s_+\,\IdMa\big)\nonumber\\
	&= -\frac{2}{s_+^2}\Big[-\frac{6}{6a^2 + b^2\,s_+}\tr(Q_*\,\Delta Q_*)\big(\frac{c^2}{b^2}\,Q_* + \frac{1}{3}\IdMa\big) + \frac{1}{b^2}\Delta Q_*\Big]\big(Q_* - \frac{1}{6}s_+\,\IdMa\big).
\label{FAE::Eqn04}
\end{align}
On the other hand, by Corollary \ref{HarMapEqn}(iii) and Lemma \ref{DiffRepLimSurf}(iv),
\begin{equation}
\Delta Q_*\big(Q_* - \frac{1}{6}s_+\,\IdMa\big) = - \sum_{\alpha=1}^3 (\nabla_\alpha Q_*)^2 =: -(\nabla Q_*)^2.
\label{FAE::Eqn05}
\end{equation}
Taking trace yields
\begin{equation}
\tr(Q_*\,\Delta Q_*) = -|\nabla Q_*|^2.
\label{FAE::Eqn06}
\end{equation}
Substituing \eqref{FAE::Eqn05} and \eqref{FAE::Eqn06} into \eqref{FAE::Eqn04} we get \eqref{FAE::Id01}.

We now turn to the proof of \eqref{FAE::Id02}. We note that
\[
\Delta {B} = \big(\Delta {B}\big)^\perp + \big(\Delta Q_\bullet\big)^\parallel - \big(\Delta {A}\big)^\parallel.
\]
It is therefore enough to establish
\begin{align}
\big(\Delta {B}\big)^\perp
	&= -\frac{4}{s_+^2}\big[\big(\nabla {B}\big)^\parallel\,\nabla Q_* + \nabla Q_*\,\big(\nabla {B}\big)^\parallel\big]\big(Q_* - \frac{1}{6}s_+\,\IdMa\big),
\label{FAE::(ii)-1}\\
\big(\Delta Q_\bullet)^\parallel
	&= -b^2({B}\,{A} + {A}\,{B}) - \frac{6c^2}{6a^2 + b^2\,s_+}|\nabla Q_*|^2\,{B}.
\label{FAE::(ii)-2}
\end{align}

To prove \eqref{FAE::(ii)-1}, we note that, by Lemma \ref{Tangent-NormalSpaces},
\begin{equation}
\frac{1}{3}s_+\,{B} - {B}\,Q_* - Q_*\,{B} = 0.
\label{FAE::Eqn08}
\end{equation}
This implies
\[
\frac{1}{3}s_+\,\Delta {B} - \Delta {B}\,Q_* - Q_*\,\Delta {B} = {B}\,\Delta\,Q_* + \Delta Q_*\,{B} + 2\nabla {B}\,\nabla Q_* + 2\nabla Q_*\,\nabla {B}.
\]
Hence, by \eqref{TNS::Id05} in Lemma \ref{Tangent-NormalSpaces},
\begin{align}
\big(\Delta {B}\big)^\perp
	&= -\frac{1}{s_+^2}\big[{B}\,\Delta\,Q_* + \Delta Q_*\,{B} + 2\nabla {B}\,\nabla Q_* + 2\nabla Q_*\,\nabla {B}\big]\big(Q_* - \frac{1}{6}s_+\,\IdMa\big)\nonumber\\
		&\qquad\qquad -\frac{1}{s_+^2}\big(Q_* - \frac{1}{6}s_+\,\IdMa\big)\big[{B}\,\Delta\,Q_* + \Delta Q_*\,{B} + 2\nabla {B}\,\nabla Q_* + 2\nabla Q_*\,\nabla {B}\big].
\label{FAE::Eqn09}
\end{align}
To proceed, note that $\Delta Q_*$ $\in$ $(T_{Q_*}\mcS_*)^\perp$ (by Lemma \ref{Tangent-NormalSpaces} and Corollary \ref{HarMapEqn}(i)) and $\nabla Q_*$ $\in$ $T_{Q_*}\mcS_*$. Applying Lemma \ref{InterRel} we see that
\begin{multline*}
\big[{B}\,\Delta\,Q_* + \Delta Q_*\,{B} + 2(\nabla {B})^\perp\,\nabla Q_* + 2\nabla Q_*\,(\nabla {B})^\perp\big]\big(Q_* - \frac{1}{6}s_+\,\IdMa\big)\\
		+ \big(Q_* - \frac{1}{6}s_+\,\IdMa\big)\big[{B}\,\Delta\,Q_* + \Delta Q_*\,{B} + 2(\nabla {B})^\perp\,\nabla Q_* + 2\nabla Q_*\,(\nabla {B})^\perp\big] \in T_{Q_*} \mcS_*,
\end{multline*}
and
\begin{multline*}
\big[2(\nabla {B})^\parallel\,\nabla Q_* + 2\nabla Q_*\,(\nabla {B})^\parallel\big]\big(Q_* - \frac{1}{6}s_+\,\IdMa\big)\\
		+ \big(Q_* - \frac{1}{6}s_+\,\IdMa\big)\big[2(\nabla {B})^\parallel\,\nabla Q_* + 2\nabla Q_*\,(\nabla {B})^\parallel\big] \in (T_{Q_*} \mcS_*)^\perp.
\end{multline*}
Hence, by projecting both sides of \eqref{FAE::Eqn09} onto the normal space $(T_{Q_*}\mcS_*)^\perp$, we arrive at
\begin{multline}
\big(\Delta {B}\big)^\perp = -\frac{2}{s_+^2}\big[(\nabla {B})^\parallel\,\nabla Q_* + \nabla Q_*\,(\nabla {B})^\parallel\big]\big(Q_* - \frac{1}{6}s_+\,\IdMa\big)\\
		+ \frac{2}{s_+^2}\big(Q_* - \frac{1}{6}s_+\,\IdMa\big)\big[(\nabla {B})^\parallel\,\nabla Q_* + \nabla Q_*\,(\nabla {B})^\parallel\big] \in (T_{Q_*} \mcS_*)^\perp.
\label{FAE::Eqn09x}
\end{multline}
Observe that, by Lemma \ref{InterRel} again, $(\nabla {B})^\parallel\,\nabla Q_* + \nabla Q_*\,(\nabla {B})^\parallel$ belongs to $(T_{Q_*}\mcS_*)^\perp$ and so commutes with $Q_*$. Thus, the right hand side of \eqref{FAE::Eqn09x} is equal to
\[
-\frac{4}{s_+^2}\big[\big(\nabla {B}\big)^\parallel\,\nabla Q_* + \nabla Q_*\,\big(\nabla {B}\big)^\parallel\big]\big(Q_* - \frac{1}{6}s_+\,\IdMa\big).
\]
\eqref{FAE::(ii)-1} follows.

Finally, we prove \eqref{FAE::(ii)-2}. Recall that by Lemma \ref{Tangent-NormalSpaces} and Corollary \ref{HarMapEqn}, the first bracket term on the right hand side of \eqref{FAE::Eqn01} belongs to $(T_{Q_*} \mcS_*)^\perp$. Thus, by projecting \eqref{FAE::Eqn01} onto $T_{Q_*}\mcS_*$, dividing by $L_k$ and then letting $L_k$ $\rightarrow$ $0$, we get
\[
\big(\Delta Q_\bullet)^\parallel
	= \Big[-b^2Q_\bullet^2 + 2c^2\,\tr(Q_*\,Q_\bullet)\,Q_\bullet\Big]^\parallel.
\]
Using Lemma \ref{InterRel}, \eqref{FAE::Eqn03} and \eqref{FAE::Eqn06}, we then obtain \eqref{FAE::(ii)-2}. The proof is complete.
\eproof

\begin{proposition}\label{ProjEqn}
Let $Q_L$ $\in$ $H^1(\Omega, \mcS_0)$ be a critical point for $I_L$. Assume furthermore that $Q_L(\Omega)$ is contained in a tubular neighborhood of $\mcS_*$ that projects smoothly onto $\mcS_*$. Then the orthogonal projection $Q_L^\sharp$ of $Q_L$ onto $\mcS_*$ satisfies in $\Omega$ the equation
\begin{multline}
\Delta Q_L^\sharp
	= -\frac{2}{s_+^2}|\nabla Q_L^\sharp|^2\,Q_L^\sharp + \frac{2}{s_+}\Big[\sum_{\alpha = 1}^3(\nabla_\alpha Q_L^\sharp)^2 - \frac{1}{3}|\nabla Q_L^\sharp|^2\,\IdMa\Big]\\
		 - \Big[T_L^{-1}\big(\frac{1}{s_+}\,Q_L^\sharp - \frac{2}{3}\,\IdMa\big)W_L - W_L\big(\frac{1}{s_+}\,Q_L^\sharp - \frac{2}{3}\,\IdMa\big)T_L^{-1}\Big],
\label{ProjEqn::Est1}
\end{multline}
where 
\begin{align}
W_L
	&= 2\nabla Q_L^\sharp\,\nabla[(Q_L^\sharp)^{-1}\,Q_L]\,Q_L^\sharp - 2Q_L^\sharp\,\nabla[(Q_L^\sharp)^{-1}\,Q_L]\,\nabla Q_L^\sharp\nonumber\\
		&\qquad\qquad - \frac{1}{s_+}Q_L\sum_{\alpha=1}^3 (\nabla_\alpha Q_L^\sharp)^2 + \frac{1}{s_+}\sum_{\alpha=1}^3 (\nabla_\alpha Q_L^\sharp)^2\,Q_L,\label{ProjEqn::Est2}\\
T_L
	&= Q_L - \frac{2}{9}s_+\,\tr[(Q_L^\sharp)^{-1}\,Q_L]\,\IdMa + \beta\Big[\frac{1}{s_+}\,Q_L^\sharp + \frac{1}{3}\,\IdMa\Big],\label{ProjEqn::Est3}
\end{align}
and $\beta$ is an arbitrary nonzero real number.
\end{proposition}

\bproof We will drop the subscript $L$ for convenience. Being a critical point of $I_L$, $Q$ satisfies \eqref{L-DG::ELEqn}, i.e.
\[
L\Delta Q = -a^2Q - b^2\big[Q^2 - \frac{1}{3}\tr(Q^2)\IdMa\big] + c^2\,\tr(Q^2)\,Q.
\]
Let $K$ $=$ $(Q^\sharp)^{-1}Q$. By definition, $Q - Q^\sharp$ is normal to the tangent plane to $\mcS_*$ at $Q^\sharp$, which implies in view of Lemma \ref{Tangent-NormalSpaces} that $Q$, $Q^\sharp$ and $K$ commutes with one another. In particular,
\[
Q = K\,Q^\sharp = Q^\sharp\,K,
\]
and so,
\[
\Delta Q = K\,\Delta Q^\sharp + 2\nabla K\,\nabla Q^\sharp + \Delta K\,Q^\sharp = \Delta Q^\sharp\, K + 2\nabla Q^\sharp\,\nabla K + Q^\sharp\,\Delta K.
\]
On the other hand, by \eqref{L-DG::ELEqn}, $\Delta Q$ commutes with with $Q^\sharp$. It follows that
\[
Q^\sharp\big[K\,\Delta Q^\sharp + 2\nabla K\,\nabla Q^\sharp + \Delta K\,Q^\sharp\big] = \big[\Delta Q^\sharp\, K + 2\nabla Q^\sharp\,\nabla K + Q^\sharp\,\Delta K\big]Q^\sharp,
\]
which implies
\begin{equation}
Q\,\Delta Q^\sharp - \Delta Q^\sharp\,Q = 2\nabla Q^\sharp\,\nabla K\,Q^\sharp - 2Q^\sharp\,\nabla K\,\nabla Q^\sharp =: \hat W.
\label{ProjEqn01}
\end{equation}

Now write
\begin{equation}
\Delta Q^\sharp = X + Y, \text{ where } X \in T_{Q^\sharp}\mcS_* \text{ and } Y \in (T_{Q^\sharp}\mcS_*)^\perp.
\label{ProjEqn01bis}
\end{equation}
It is well-known that $Y$ $=$ $\textrm{II}(\nabla Q^\sharp,\nabla Q^\sharp)(Q^\sharp)$ (see e.g. \cite{Helein}). Hence, by Lemma \ref{2ndForm},
\begin{equation}
Y = -\frac{2}{s_+^2}|\nabla Q^\sharp|^2\,Q^\sharp + \frac{2}{s_+}\Big[\sum_{\alpha = 1}^3(\nabla_\alpha Q^\sharp)^2 - \frac{1}{3}|\nabla Q^\sharp|^2\,\IdMa\Big].
\label{ProjEqn02}
\end{equation}
Therefore, by \eqref{ProjEqn01} and \eqref{ProjEqn::Est2}, we have
\begin{equation}
Q\,X - X\,Q = \hat W - Q\,Y + Y\,Q = W.
\label{ProjEqn03}
\end{equation}
In addition, as $X$ $\in$ $T_{Q^\sharp}\mcS_*$, Lemma \ref{Tangent-NormalSpaces} gives
\begin{equation}
Q^\sharp\,X + X\,Q^\sharp = \frac{1}{3}s_+\,X.
\label{ProjEqn04}
\end{equation}
The conclusion then can be derived directly from \eqref{ProjEqn02}, \eqref{ProjEqn03} and \eqref{ProjEqn04} using the identities established in Lemma \ref{MysIds}. Since the calculation is lengthy and is somewhat of minor importance to our goal, we deter it to the appendix.
\eproof


\section{$C^{1,\alpha}$-convergence}\label{SEC-C1AlphaConv}

Let $Q_L$ be a minimizer of the minimization problem \eqref{MinProb}. It is customary to show that for any $L_{k}$ $\rightarrow$ $0$, there exists a subsequence $L_{k'}$ such that $Q_{L_{k'}}$ converges strongly in $H^1$ to some $Q_*$ $\in$ $H^1(\Omega,\mcS_*)$ which is $I_*$-minimizing. In this section, we are interested in establishing $C^{1,\alpha}$ convergence for such sequence of minimizers.

We will use the following notation. For a function $u$ defined on $\bar\Omega$, we denote by $\Sing(u)$ the singular set of $u$, i.e. the set of points $x$ $\in$ $\bar\Omega$ such that there is no neighborhood $U$ of $x$ for which $u\big|_{U}$ is smooth.

\begin{proposition}\label{C1AlphaConv}
Assume that $Q_{L_k}$ converges strongly in $H^1(\Omega, \mcS_0)$ to a minimizer $Q_*$ $\in$ $H^1(\Omega,\mcS_*)$ of $I_*$ for some sequence $L_k$ $\rightarrow$ $0$. For any compact subset $K$ of $\bar\Omega\setminus \Sing(Q_*)$, there exists $\bar L$ $=$ $\bar L(a^2,b^2,c^2,\Omega,K,Q_b,Q_*)$ $>$ $0$ such that for any $\alpha$ $\in$ $(0,1)$ there holds
\[
\|Q_{L_k}\|_{C^{1,\alpha}(K)} \leq C(a^2,b^2,c^2,\Omega,K,Q_b,Q_*,\alpha) \text{ for any } L_k \leq \bar L.
\]
In particular, away from $\Sing(Q_*)$, $Q_{L_k}$ converges to $Q_*$ in $C^{1,\alpha}$-norm.
\end{proposition}

It is convenient to introduce
\begin{equation}
\tilde f_B(Q) = f_B(Q) - \min_{\mcS_0} f_B,
\label{VarBulkEDensity}
\end{equation}
and
\begin{equation}
\tilde I_L[Q]= \int_\Omega \Big[\frac{L}{2}|\nabla Q|^2 + \tilde f_B(Q)\Big]\,dx, \qquad Q \in H^1(\Omega, \mcS_0).
\label{L-DGFuncVar}
\end{equation}
Clearly, $Q$ is a minimizer for $I_L$ if and only if it is a minimizer for $\tilde I_L$.

The proof of Proposition \ref{C1AlphaConv} consists of two main steps. First, we prove a uniform bound for $\Delta Q_{L_k}$ under an additionally assumed uniform $C^1$ bound. Second, we prove the required uniform $C^1$ bound. It is readily seen that Proposition \ref{C1AlphaConv} follows from such estimates. We will frequently use the following two results, which are variants of \cite[Lemma 2]{BBH}.

\begin{lemma}\label{BBHMaxPrin}
Let $B_R$ $\subset$ $\RR^3$ be a ball centered at the origin and of radius $R$. Assume that $u$ $\in$ ${\rm Lip}(B_R) \cap C^0(\bar B_R)$ satisfies in the weak sense
\[
-L\,\Delta u + a\,u \leq C \text{ in } B_R,
\]
where $L$, $a$ and $C$ are positive constants. Then, for $\alpha$ $\geq$ $0$, there holds in $B_R$,
\[
u(x) \leq \frac{C}{a} + 2\exp\Big(-\sqrt{\frac{a}{L}}\,\frac{R - |x|}{2}\Big)\,\sup_{\partial B_r}u^+ \leq \frac{C}{a} + \frac{C(\alpha)\,L^\alpha}{a^\alpha\,(R - |x|)^{2\alpha}}\,\sup_{\partial B_R}u^+.
\]
\end{lemma}

\bproof Define a function $\phi$ in $B_R$ by
\[
\phi(x) = \frac{\sinh\big(\sqrt{\frac{a}{L}}|x|\big)\,R}{\sinh\big(\sqrt{\frac{a}{L}}R\big)\,|x|}
\]
It is routine to check that $\phi$ $\in$ $C^2(B_R)$ and $-L\,\Delta\phi + a\,\phi$ $=$ $0$. It thus follows from the maximum principle that
\[
u(x) \leq \frac{C}{a} + \Big[\sup_{\partial B_r}u^+\Big]\,\phi(x) \text{ in } B_R.
\]
Next, as $-L\,\Delta\phi + a\,\phi$ $=$ $0$, the maximum principle implies that, for $|x|$ $\leq$ $R/2$,
\[
\phi(x) \leq \phi(R/2) = \frac{2\sinh\big(\sqrt{\frac{a}{L}}\frac{R}{2}\big)}{\sinh\big(\sqrt{\frac{a}{L}}R\big)} \leq 2\exp\Big(-\sqrt{\frac{a}{L}}\frac{R}{2}\Big) \leq \frac{C(\alpha)\,L^\alpha}{a^\alpha\,R^{2\alpha}}.
\]
In addition, for $R/2$ $<$ $|x|$ $<$ $R$, we have
\[
\phi(x) \leq 2\exp\Big(-\sqrt{\frac{a}{L}}(R - |x|)\Big) \leq \frac{C(\alpha)\,L^\alpha}{a^\alpha\,(R - |x|)^{2\alpha}}.
\]
The assertion follows from the last three estimates.
\eproof

\begin{lemma}\label{BBHMaxPrin::Bdry}
Let $\Omega$ be a domain in $\RR^3$ and $x_0$ $\in$ $\partial\Omega$. Assume that $u$ $\in$ ${\rm Lip}(B_R(x_0) \cap \Omega) \cap C^0(\bar B_R(x_0) \cap \bar \Omega)$ satisfies in the weak sense
\begin{align*}
&-L\,\Delta u + a\,u \leq C \text{ in } B_R(x_0) \cap \Omega,\\
&u = 0 \text{ on } B_R(x_0) \cap \partial\Omega
\end{align*}
where $L$, $a$ and $C$ are positive constants. Then, for $\alpha$ $\geq$ $0$,
\[
u(x) \leq \frac{C}{a} + \frac{C(\alpha)\,L^\alpha}{a^\alpha\,(R - |x|)^{2\alpha}}\,\sup_{\partial B_R(x_0) \cap \Omega}u^+ \text{ in } B_R(x_0) \cap \Omega.
\]
\end{lemma}

\bproof Extend $u$ to a function $\tilde u$ defined on $B_R(x_0)$ by zero in $B_R \setminus \Omega$. The result then follows from Lemma \ref{BBHMaxPrin::Bdry}.
\eproof

We turn now to the proof of Proposition \ref{C1AlphaConv}. For our purpose, it is useful to note that, by Lemma \ref{DiffRepLimSurf}(iv),
\begin{equation}
\mcS_* = \Big\{Q \in \mcS_0: Q^2 - \frac{1}{3}s_+\,Q - \frac{2}{9}s_+^2\,\IdMa = 0\Big\},
\label{LimitSurf::AlterRep}
\end{equation}
where $s_+$ is given in \eqref{LimitSurf::Rep}. In other words, the minimal polynomial of any matrix in $\mcS_*$ is $\lambda^2 - \frac{1}{3}s_+\,\lambda - \frac{2}{9}s_+^2$. Furthermore, for $Q$ belonging to a small tubular neighborhood of the limit manifold $\mcS_*$, the norm of of $Q^2 - \frac{1}{3}s_+\,Q - \frac{2}{9}s_+^2\,\IdMa$ is comparable to the distance from $Q$ to $\mcS_*$, which is a consequence of the following lemma.

\begin{lemma}\label{DistComp}
Let $\alpha$, $\beta$, $\gamma$, $\mu$, $\nu$ and $\delta$ be real numbers and define for $Q$ $\in$ $\mcS_0$,
\begin{multline*}
h(Q) = \alpha\,(\tr(Q^2))^3 + \beta\,s_+\,\tr(Q^2)\,\tr(Q^3)+ \gamma\,s_+^2\,(\tr(Q^2))^2\\
	+ \mu\,s_+^3\,\tr(Q^3) + \nu\,s_+^4\,\tr(Q^2) + \delta\,s_+^6.
\end{multline*}
If
\begin{align*}
&\frac{8}{27}\alpha + \frac{4}{27}\beta + \frac{4}{9}\gamma + \frac{2}{9}\mu + \frac{2}{3}\nu + \delta = 0,\\
&\frac{8}{3}\alpha + \frac{10}{9}\beta + \frac{8}{3}\gamma + \mu + 2\nu = 0,\\
&(16\alpha + 4\beta + 8\gamma)^2 - (16\alpha + 6\beta + 8\gamma + 3\mu)^2 > 0,
\end{align*}
then there exist $\epsilon$ $>$ $0$ and $C$ $>$ $0$ such that
\[
\frac{1}{C}\,\dist(Q,\mcS_*)^2 \leq h(Q) \leq C\,\dist(Q,\mcS_*)^2
\]
for any $Q$ $\in$ $\mcS_0$ satisfying $\dist(Q,\mcS_*)$ $<$ $\epsilon$.
\end{lemma}

\bproof Let $x$, $y$, $-(x+y)$ be the eigenvalues of $Q$. For $Q$ close to $\mcS_*$, we can further assume that $x$ and $y$ are close to $-\frac{s_+}{3}$ while $-(x+y)$ is close to $\frac{2s_+}{3}$. We have
\begin{multline*}
h(Q) = 8\alpha\,(x^2 + y^2 + xy)^3 - 6\beta\,s_+\,(x^2 + y^2 + xy)xy(x+y)\\
	+ 4\gamma\,s_+^2\,(x^2 + y^2 + xy)^2 - 3\,s_+^3\mu\,xy(x+y) + 2\nu\,s_+^4\,(x^2 + y^2 + xy) + \delta\,s_+^6.
\end{multline*}
By a simple calculation using the given constraints on $\alpha$, $\beta$, $\gamma$, $\mu$ and $\nu$ we have
\begin{align*}
h\big(-\frac{s_+}{3},-\frac{s_+}{3}\big)
	&= 0,\\
\partial_x h\big(-\frac{s_+}{3},-\frac{s_+}{3}\big)
	&= \partial_y h\big(-\frac{s_+}{3},-\frac{s_+}{3}\big) = -\Big[\frac{8}{3}\alpha + \frac{10}{9}\beta + \frac{8}{3}\gamma + \mu + 2\nu\Big]s_+^5 = 0,\\
\partial_{xx} h\big(-\frac{s_+}{3},-\frac{s_+}{3}\big)
	&= \partial_{yy} h\big(-\frac{s_+}{3},-\frac{s_+}{3}\big)\\
	&= \Big[2\big(\frac{8}{3}\alpha + \frac{10}{9}\beta + \frac{8}{3}\gamma + \mu + 2\nu\big) + 16\alpha + 4\beta + 8\gamma\Big]s_+^4\\
	&= \big[16\alpha + 4\beta + 8\gamma\big]s_+^4,\\
\partial_{xy} h\big(-\frac{s_+}{3},-\frac{s_+}{3}\big)
	&= \Big[\big(\frac{8}{3}\alpha + \frac{10}{9}\beta + \frac{8}{3}\gamma + \mu + 2\nu\big) + 16\alpha + 6\beta + 8\gamma + 3\mu\Big]s_+^4\\
	&= \big[16\alpha + 6\beta + 8\gamma + 3\mu\big]s_+^4.
\end{align*}
The assertion follows from a simple application of the Taylor expansion theorem.
\eproof

As argued before, it is of relevance to see how $Q_L^2 - \frac{1}{3}s_+\,Q_L - \frac{2}{9}s_+^2\,\IdMa$ converges to zero. It is convenience to introduce
\begin{equation}
X_L = \frac{1}{L} \Big[Q_L^2 - \frac{1}{3}s_+\,Q_L - \frac{2}{9}s_+^2\,\IdMa\Big].
\label{XDef}
\end{equation}
The following result gives a rate of convergence for $Q_L^2 - \frac{1}{3}s_+\,Q_L - \frac{2}{9}s_+^2\,\IdMa$ provided that $Q_L$ is sufficiently close to $\mcS_*$ and a gradient bound is known.

\begin{proposition}\label{MinPolyConv}
There exist $\delta_0$ $=$ $\delta_0(a^2,b^2,c^2)$ $>$ $0$ such that if $\dist(Q_L,\mcS_*)$ $\leq$ $\delta_0$ and $|\nabla Q_L|$ $\leq$ $C_1$ in some $B_r(x_0) \cap \Omega$ for some $x_0$ $\in$ $\bar\Omega$, $L$ $\in$ $(0,1)$ and $C_1$ $>$ $0$, then
\[
|X_L| \leq C(a^2,b^2,c^2,C_1)[1 + r^{-2}] \text{ in } B_{r/2}(x_0) \cap \Omega.
\]
\end{proposition}

\bproof For convenience, we write $Q$ and $X$ for $Q_L$ and $X_L$, respectively. For $\delta_0$ sufficiently small, we define $Q^\sharp$ as the orthogonal projection of $Q$ onto $\mcS_*$, i.e. 
\[
|Q(x) - Q^\sharp(x)| = \dist(Q(x),\mcS_*).
\]

Set $Y = L^2\,|X|^2$. We have
\begin{align*}
Y
	&= \tr\Big(Q^4 - \frac{2}{3}s_+\,Q^3 - \frac{1}{3}s_+^2\,Q^2 + \frac{4}{27}s_+^3\,Q + \frac{4}{81}s_+^4\,\IdMa\Big)\\
	&= \tr(Q^4) - \frac{2}{3}s_+\,\tr(Q^3) - \frac{1}{3}s_+^2\,\tr(Q^2) + \frac{4}{27}\,s_+^4\\
	&= \frac{1}{2}\tr(Q^2)^2 - \frac{2}{3}s_+\,\tr(Q^3) - \frac{1}{3}s_+^2\,\tr(Q^2) + \frac{4}{27}\,s_+^4\\
	&=: g(Q).
\end{align*}
Applying Lemma \ref{DistComp} we can choose $\delta_0$ sufficiently small such that
\begin{equation}
\frac{1}{C}\,Y \leq |Q - Q^\sharp|^2 \leq C\,Y \text{ in } B_r(x_0) \cap \Omega.
\label{MinPolyConv::WeakBound}
\end{equation}

We have
\begin{equation}
\Delta Y
	= g_{ij}(Q)\,\Delta Q_{ij} + g_{ij,pq}(Q)\,\nabla Q_{ij} \cdot \nabla Q_{pq},
\label{MinPolyConv::Eqn01}
\end{equation}
where $g_{ij}(Q)$ $=$ $\frac{\partial g}{\partial Q_{ij}}(Q)$ and $g_{ij,pq}(Q)$ $=$ $\frac{\partial^2 g}{\partial Q_{ij}\partial Q_{pq}}(Q)$.

Since $g$ vanishes on $\mcS_*$ (in view of \eqref{LimitSurf::AlterRep}), it achieves its minimum everywhere on $\mcS_*$ and in particular at $Q^\sharp$ (as a function on $M_{3\times 3}$). Hence, by \eqref{MinPolyConv::WeakBound} and the given gradient bound, we have
\begin{equation}
g_{ij,pq}(Q)\,\nabla Q_{ij} \cdot \nabla Q_{pq} \geq g_{ij,pq}(Q^\sharp)\,\nabla Q_{ij} \cdot \nabla Q_{pq} - C\,\sqrt{Y} \geq -C\,\sqrt{Y}  \text{ in } B_r(x_0) \cap \Omega.
\label{MinPolyConv::Eqn02}
\end{equation}

On the other hand, by recalling \eqref{L-DG::ELEqn}, we have
\begin{align*}
g_{ij}(Q)\,L\,\Delta Q_{ij}
	&= \Big[2\,\tr(Q^2)\,Q_{ji} - 2s_+(Q^2)_{ji} - \frac{2}{3}s_+^2\,Q_{ji}\Big]\\
		&\qquad\qquad \times \Big[- a^2 Q_{ij} - b^2[(Q^2)_{ij} - \frac{1}{3}\tr(Q^2)\delta_{ij}] + c^2\,\tr(Q^2)\,Q_{ij}\Big]\\
	&= 2c^2\,(\tr(Q^2))^3 + (- 2b^2 - 2c^2\,s_+)\,\tr(Q^2)\,\tr(Q^3)\\
		&\qquad\qquad  + (-2a^2 + \frac{1}{3}b^2\,s_+ - \frac{2}{3}c^2\,s_+^2)\,(\tr(Q^2))^2 + (2a^2\,s_+ + \frac{2}{3}b^2s_+^2)\tr(Q^3)\\
		&\qquad\qquad + \frac{2}{3}a^2\,s_+^2\,\tr(Q^2)\\
	&= (3a^2\frac{1}{s_+^2} + b^2\frac{1}{s_+})\,(\tr(Q^2))^3 + (-3a^2\frac{1}{s_+} - 3b^2)\,\tr(Q^2)\,\tr(Q^3)\\
		&\qquad\qquad  - 3a^2\,(\tr(Q^2))^2 + (2a^2\,s_+ + \frac{2}{3}b^2s_+^2)\tr(Q^3) + \frac{2}{3}a^2\,s_+^2\,\tr(Q^2).
\end{align*}
Applying Lemma \ref{DistComp} and using \eqref{MinPolyConv::WeakBound}, we get, for $\delta_0$ sufficiently small, 
\begin{equation}
g_{ij}(Q)\,L\,\Delta Q_{ij} \geq C\,|Q - Q^\sharp|^2 \geq CY  \text{ in } B_r(x_0) \cap \Omega.
\label{MinPolyConv::Eqn03}
\end{equation}

Combining \eqref{MinPolyConv::Eqn01}, \eqref{MinPolyConv::Eqn02} and \eqref{MinPolyConv::Eqn03} we get
\[
L\,\Delta Y \geq 2C_2\,Y - C_2'\,L\,\sqrt{Y} \geq C_2\,Y - C_3\,L^2  \text{ in } B_r(x_0) \cap \Omega,
\]
where $C_2$, $C_2'$, $C_3$ are positive constants. We thus have
\begin{equation}
-L\,\Delta Y + C_2\,Y \leq C_3\,L^2  \text{ in } B_r(x_0) \cap \Omega.
\label{MinPolyConv::Eqn04}
\end{equation}
It is readily seen that the assertion follows from Lemmas \ref{BBHMaxPrin} and \ref{BBHMaxPrin::Bdry}.
\eproof

As a consequence of Proposition \ref{MinPolyConv}, we have the following estimate for $\Delta Q_L$.

\begin{corollary}\label{UniLapBnd}
There exist $\delta_0$ $=$ $\delta_0(a^2,b^2,c^2)$ $>$ $0$ such that if $\dist(Q_L,\mcS_*)$ $\leq$ $\delta_0$ and $|\nabla Q_L|$ $\leq$ $C_1$ in some $B_r(x_0) \cap \Omega$ for some $x_0$ $\in$ $\bar\Omega$, $L$ $\in$ $(0,1)$ and $C_1$ $>$ $0$, then
\[
|\Delta Q_L| \leq C(a^2,b^2,c^2,C_1)[1 + r^{-2}] \text{ in } B_{r/2}(x_0) \cap \Omega.
\]
\end{corollary}

\bproof By \eqref{L-DG::ELEqn}, 
\[
L^2|\Delta Q_L|^2
	= \Big|-a^2Q - b^2\big(Q^2 - \frac{1}{3}\tr(Q^2)\,\IdMa\big) + c^2\,\tr(Q^2)\,Q\Big|^2 =: h(Q).
\]
A direct computation using Lemma \ref{DistComp} shows that $h(Q)$ and $|Q^2 - \frac{1}{3}s_+Q - \frac{2}{9}s_+^2\,\IdMa|^2$ are comparable near $\mcS_*$. The assertion follows from Proposition \ref{MinPolyConv}.
\eproof

\begin{corollary}\label{TildeFBConv}
There exist $\delta_0$ $=$ $\delta_0(a^2,b^2,c^2)$ $>$ $0$ such that if $\dist(Q_L,\mcS_*)$ $\leq$ $\delta_0$ and $|\nabla Q_L|$ $\leq$ $C_1$ in some $B_r(x_0) \cap \Omega$ for some $x_0$ $\in$ $\bar\Omega$, $L$ $\in$ $(0,1)$ and $C_1$ $>$ $0$, then
\[
\tilde f_B(Q_{L_k}) \leq C(a^2,b^2,c^2,C_1)[1 + r^{-4}]\,L^2 \text{ in } B_{r/2}(x_0) \cap \Omega.
\]
\end{corollary}

\bproof Using Lemma \ref{DistComp}, we see that $\tilde f_B(Q)$ and $|Q^2 - \frac{1}{3}s_+Q - \frac{2}{9}s_+^2\,\IdMa|^2$ are comparable near $\mcS_*$. The assertion follows from Proposition \ref{MinPolyConv}.
\eproof

\begin{corollary}\label{|Q2|Conv}
There exist $\delta_0$ $=$ $\delta_0(a^2,b^2,c^2)$ $>$ $0$ such that if $\dist(Q_L,\mcS_*)$ $\leq$ $\delta_0$ and $|\nabla Q_L|$ $\leq$ $C_1$ in some $B_r(x_0) \cap \Omega$ for some $x_0$ $\in$ $\bar\Omega$, $L$ $\in$ $(0,1)$ and $C_1$ $>$ $0$, then
\[
0 \leq  \frac{2}{3}s_+^2 - \tr(Q_{L}^2) \leq C[1 + r^{-2}]\,L \text{ in } B_{r/2}(x_0) \cap \Omega.
\]
\end{corollary}

\bproof The second inequality follows from Proposition \ref{MinPolyConv} and that $\tr(Q^2) - \frac{2}{3}s_+^2$ $=$ $\tr(Q^2 - \frac{1}{3}s_+\,Q - \frac{2}{9}s_+^2\,\IdMa)$. The first inequality follows from \cite[Proposition 3]{Ma-Za}.
\eproof

We turn to establishing a uniform gradient bound for $Q_L$. We will use Bochner technique. Note that uniform interior gradient estimate was already established in \cite{Ma-Za}. To adapt the argument therein to our situation here, we need some additional information about the gradient on the boundary. To this end, we split $Q_L$ into $Q_L^\sharp$, the projection of $Q_L$ onto $\mcS_*$, and $Q_L - Q_L^\sharp$, and obtain boundary gradient estimates for each of them separately. It turns out that $Q_L - Q_L^\sharp$ can be controlled rather easily by the minimal polynomial. For the other part, we use Proposition \ref{ProjEqn}. We start with:

\begin{lemma}\label{C1BndAuxRslt}
Given $x$ $\in$ $\Omega$ and $r$ $>$ $10\,\dist(x,\partial\Omega)$, assume that $Q_L$ satisfies
\begin{equation}
\sup_{B_r(x) \cap \Omega} |\nabla Q_L| \leq C_1.
\label{CBAR::Hyp1}
\end{equation}
There exist $\delta_0$ $>$ $0$ and $C$ $>$ $0$ such that if $\dist(Q_L,\mcS_*)$ $\leq$ $\delta_0$ in $B_r(x) \cap \Omega$, then
Uni\begin{multline*}
\sup_{\partial\Omega \cap B_{r/2}(x)}|\nabla (Q_L - Q_L^\sharp)| \leq C\Big[r^{-1}\|Q - Q^\sharp\|_{L^\infty(B_r(x) \cap \Omega)}\\
	+ r^{-3/2}\|\nabla Q\|_{L^2(B_r(x) \cap \Omega)} + r^{1/4}\|\nabla Q\|_{L^2(B_r(x) \cap \Omega)}^{1/2}\Big],
\end{multline*}
where $Q_L^\sharp$ is the orthogonal projection of $Q_L$ onto $\mcS_*$, defined by
\[
|Q_L(y) - Q_L^\sharp(y)| = \min\big\{|Q_L(y) - R|: R \in \mcS_*\big\}.
\]
\end{lemma}

We will need a few simple results whose proof will be given after the proof of Lemma \ref{C1BndAuxRslt}.

\begin{lemma}\label{CalcLemma}
Let $g:$ $B_1(0)$ $\subset$ $\RR^m$ $\rightarrow$ $\RR$ be a smooth function which has a local minimum at the origin and $g(0)$ $=$ $0$. Then there exist some positive constants $\epsilon$ $>$ $0$ and $C$ $>$ $0$ depending only on $g$ such that
\[
\Big(2g\,\nabla_{ij} g - \nabla_i g\,\nabla_j g\Big)(x)\,a_i\,a_j \geq -C\,|x|^3\,|a|^2
\]
for any $x$ $\in$ $B_\epsilon(0)$ and $a$ $\in$ $\RR^m$.
\end{lemma}

\begin{lemma}\label{SimpleEllEst}
Assume that $f$ $\in$ $L^p(B_r^+(0))$, $p$ $>$ $3$, $g$ $\in$ $C^2(B_r(0) \cap \{x_3 = 0\})$ and $u$ $\in$ $W^{1,\infty}(B_r^+(0))$ satisfy
\begin{align*}
-a_{ij}\,\nabla_{ij} u + b_i\,\nabla_i u &= f \text{ in } B_r^+(0),\\
u &= g \text{ on } B_r(0) \cap \{x_3 = 0\},
\end{align*}
where the coefficients $a_{ij}$ and $b_i$ are continuous and satisfy $\lambda$ $\leq$ $(a_{ij})$ $\leq$ $\Lambda$ and $|b_i|$ $\leq$ $\Lambda$ for some $0$ $<$ $\lambda$ $\leq$ $\Lambda$ $<$ $\infty$. Then
\begin{multline*}
\sup_{B_{r/2}^+(0)} |\nabla u|
	\leq C\big[\sup_{B_r(0) \cap \{x_3 = 0\}} |\nabla_T g| + r\,\sup_{B_r(0) \cap \{x_3 = 0\}}|\nabla_T^2 g|\\
		+ r^{-3/2}\|\nabla u\|_{L^2(B_r^+(0))} + r^{1-3/p}\,\|f\|_{L^p(B_r^+(0))}\big],
\end{multline*}
where $\nabla_T$, $\nabla_T^2$ denote the horizontal gradient and Hessian, respectively.
\end{lemma}

\noindent{\bf Proof of Lemma \ref{C1BndAuxRslt}.} We will drop the subscript $L$. Also we will abbreviate $B_r$ for $B_r(0)$. In the proof we will use without explicit mentioning the following two facts: (a) $|Q|$ $\leq$ $\sqrt{\frac{2}{3}}s_+$ in $\Omega$ (see \cite[Proposition 3]{Ma-Za}), and (b) $|R|$ $=$ $\sqrt{\frac{2}{3}}s_+$ for any $R$ $\in$ $\mcS_*$.

Set $\hat X = L\,X$. Arguing as in the proof of Proposition \ref{MinPolyConv} (cf. \eqref{MinPolyConv::Eqn01}), for $\delta_0$ sufficiently small, there holds
\begin{equation}
\frac{1}{C}\,|\hat X| \leq |Q - Q^\sharp| \leq C\,|\hat X| \text{ in } B_r \cap \Omega.
\label{CBAR::Eqn02}
\end{equation}

Note that $|\hat X|$ is smooth whenver $|\hat X|$ $>$ $0$ and, by \eqref{CBAR::Hyp1},
\begin{equation}
\nabla|\hat X| = \frac{\nabla |\hat X|^2}{2|\hat X|} = \frac{\tr(\nabla \hat X\,\hat X)}{|\hat X|} \leq |\nabla \hat X| \leq C|\nabla Q| \leq C\text{ in } B_r \cap \Omega \cap \{|\hat X| > 0\}.
\label{CBAR::Eqn01}
\end{equation}
Since $|\hat X|$ is continuous in $B_r \cap \Omega$, it follows that $|\hat X|$ $\in$ $W^{1,\infty}(B_r \cap \Omega)$. In particular (see e.g. \cite[p. 84]{EG}),
\begin{equation}
\nabla |\hat X| = 0 \text{ a.e. in } B_r \cap \Omega \cap \{|\hat X| = 0\}.
\label{CBAR::Eqn01x}
\end{equation}

As in the proof of Proposition \ref{MinPolyConv}, we have
\[
|\hat X|^2
	= \frac{1}{2}\tr(Q^2)^2 - \frac{2}{3}s_+\,\tr(Q^3) - \frac{1}{3}s_+^2\,\tr(Q^2) + \frac{4}{27}\,s_+^4 =: g(Q).
\]
Hence, in the region where $|\hat X|$ $>$ $0$,
\begin{equation}
\Delta |\hat X|
	= \frac{g_{ij}(Q)\,\Delta Q_{ij}}{2|\hat X|} + \frac{2g(Q)\,g_{ij,pq}(Q)\,\nabla Q_{ij} \cdot \nabla Q_{pq} - g_{ij}(Q)\,\nabla Q_{ij}\,g_{pq}(Q)\,\nabla Q_{pq}}{4|\hat X|^3},
\label{CBAR::Eqn03}
\end{equation}
where $g_{ij}(Q)$ $=$ $\frac{\partial g}{\partial Q_{ij}}(Q)$ and $g_{ij,pq}(Q)$ $=$ $\frac{\partial^2 g}{\partial Q_{ij}\partial Q_{pq}}(Q)$.

Since $g$ vanishes on $\mcS_*$ (in view of \eqref{LimitSurf::AlterRep}), it achieves its minimum everywhere on $\mcS_*$ and in particular at $Q^\sharp$ (as a function on $M_{3\times 3}$). Hence, by \eqref{CBAR::Eqn02} and Lemma \ref{CalcLemma}, for $\delta_0$ sufficiently small, we have in $B_r \cap \Omega$,
\begin{equation}
\frac{2g(Q)\,g_{ij,pq}(Q)\,\nabla Q_{ij} \cdot \nabla Q_{pq} - g_{ij}(Q)\,\nabla Q_{ij}\,g_{pq}(Q)\,\nabla Q_{pq}}{|\hat X|^3} \geq -C\,|\nabla Q|^2.
\label{CBAR::Eqn04}
\end{equation}

On the other hand, by recalling \eqref{MinPolyConv::Eqn03}, we have for $\delta_0$ sufficiently small that
\begin{equation}
g_{ij}(Q)\,L_k\,\Delta Q_{ij} \geq C\,|Q - Q^\sharp|^2 \geq 0  \text{ in } B_r \cap \Omega.
\label{CBAR::Eqn05}
\end{equation}

We now fix $\delta_0$. Combining \eqref{CBAR::Eqn01}, \eqref{CBAR::Eqn03}, \eqref{CBAR::Eqn04} and \eqref{CBAR::Eqn05} we get
\begin{equation}
-\Delta|\hat X| \leq C|\nabla Q|^2 \text{ in } B_r \cap \Omega \cap \{|\hat X| > 0\}.
\label{CBAR::Eqn05x1}
\end{equation}
We next show that the above inequality can be extended to
\begin{equation}
-\Delta|\hat X| \leq C|\nabla Q|^2 \text{ in } B_r \cap \Omega \text{ in the sense of distribution}.
\label{CBAR::Eqn05x2}
\end{equation}
Indeed, since $|\hat X|$ is smooth in $\{|\hat X| > 0\}$, by Sard's theorem, we can select a decreasing sequence $\delta_j$ $\rightarrow$ $0$ such that the level sets $\{|\hat X| = \delta_j\}$ are smooth. Fix some $\varphi$ $\in$ $C^{c,\infty}(B_r \cap \Omega)$ where $\varphi$ $\geq$ $0$ in $B_r \cap \Omega$. Then, by \eqref{CBAR::Eqn01x}, we compute
\begin{align*}
\int_{B_r \cap \Omega} \nabla|\hat X|\,\nabla \varphi\,dx
	&= -\int_{B_r \cap \Omega \cap \{|\hat X \geq \delta_j\}}\Delta |\hat X|\,\varphi\,dx + \int_{B_r \cap \Omega \cap \{|\hat X| = \delta_j\}} \nabla|\hat X| \cdot \nu\,\varphi\,d\sigma(x)\\
		&\qquad\qquad + \int_{B_r \cap \Omega \cap \{0 < |\hat X| < \delta_j\}} \nabla|\hat X|\,\nabla \varphi\,dx,
\end{align*}
where $\nu$ is the outer normal to $\{|\hat X| = \delta_j\}$ relative to the set $B_r \cap \Omega \cap \{|\hat X| > \delta_j\}$. Observe that $\nabla |\hat X| \cdot \nu$ $\leq$ $0$ on $\{|\hat X| = \delta_j\}$. Thus, by \eqref{CBAR::Eqn01} and \eqref{CBAR::Eqn05x1},
\[
\int_{B_r \cap \Omega} \nabla|\hat X|\,\nabla \varphi\,dx
	\leq C\int_{B_r \cap \Omega}|\nabla Q|^2\,\varphi\,dx + C\sup_{B_r \cap \Omega}|\nabla \varphi|\,\big|B_r \cap \Omega \cap \{0 < |\hat X| < \delta_j\}\big|.
\]
Sending $j$ $\rightarrow$ $\infty$, we obtain \eqref{CBAR::Eqn05x2}.

By \eqref{CBAR::Eqn05x2}, $|\hat X|$ $\leq$ $w$ where $w$ is the solution to
\begin{align*}
-\Delta w 
	&= C|\nabla Q|^2 \text{ in } B_r \cap \Omega,\\
w
	&= |\hat X| \text{ on } \partial (B_r \cap \Omega).
\end{align*}
Using \eqref{CBAR::Hyp1} and applying Lemma \ref{SimpleEllEst}, we obtain, for $y \in B_{r/2} \cap \Omega$,
\begin{align}
w(y)
	&\leq C\Big[r^{-3/2}\|\nabla w\|_{L^2(B_{3r/4} \cap \Omega)} + r^{1/4}\|\nabla Q\|_{L^8(B_{3r/4} \cap \Omega)}^2\Big]\dist(y,\partial\Omega)\nonumber\\
	&\leq C\Big[r^{-3/2}\|\nabla w\|_{L^2(B_{3r/4} \cap \Omega)} + r^{1/4}\|\nabla Q\|_{L^2(B_{3r/4} \cap \Omega)}^{1/2}\Big]\dist(y,\partial\Omega).
\label{CBAR::Eqn06}
\end{align}
To proceed we split $w$ $=$ $w_1 + w_2$ where $w_1$ is the solution to
\begin{align*}
-\Delta w_1 
	&= C|\nabla Q|^2 \text{ in } B_r \cap \Omega,\\
w_1
	&= 0 \text{ on } \partial (B_r \cap \Omega).
\end{align*}
By elliptic estimates and \eqref{CBAR::Hyp1}, we have
\begin{equation}
\|\nabla w_1\|_{L^2(B_{3r/4} \cap \Omega)} \leq C\,\|\Delta w_1\|_{L^2(B_r \cap \Omega)} \leq  C\,\|\nabla Q\|_{L^2(B_r \cap \Omega)}.
\label{CBAR::Eqn07}
\end{equation}
Also, as $w_2$ satisfies
\begin{align*}
-\Delta w_2
	&= 0 \text{ in } B_r \cap \Omega,\\
w_2
	&= 0 \text{ on } B_r \cap \partial \Omega,\\
w_2
	&= |\hat X| \text{ on } \partial B_r \cap \Omega,
\end{align*}
we infer that
\begin{multline}
\|\nabla w_2\|_{L^2(B_{3r/4} \cap \Omega)} \leq C\,r^{-1}\,\|w_2\|_{L^2(B_r \cap \Omega)}\\
	\leq C\,r^{1/2}\|w_2\|_{L^\infty(B_r \cap \Omega)} \leq C\,r^{1/2}\|\hat X\|_{L^\infty(B_r \cap \Omega)}.
\label{CBAR::Eqn08}
\end{multline}
Taking \eqref{CBAR::Eqn02}, \eqref{CBAR::Eqn06}-\eqref{CBAR::Eqn08} into account altogether we arrive at
\begin{multline*}
|Q - Q^\sharp|(y)
	\leq C\Big[r^{-1}\,\|Q - Q^\sharp\|_{L^\infty(B_r \cap \Omega)} + r^{-3/2}\|\nabla Q\|_{L^2(B_{3r/4} \cap \Omega)}\\
		+ r^{1/4}\|\nabla Q\|_{L^2(B_{3r/4} \cap \Omega)}^{1/2}\Big]\dist(y,\partial\Omega) \text{ in } B_{r/2} \cap \Omega.
\end{multline*}
As $|Q - Q^\sharp|$ vanishes on $\partial\Omega$, the conclusion follows.
\eproof

\noindent{\bf Proof of Lemma \ref{CalcLemma}.} Set
\[
A = 2g\,\nabla^2 g - \nabla g \otimes \nabla g.
\]

Since $g$ has a local minimum at the origin, we have the following Taylor's expansions,
\begin{align*}
g(x)
	&= \frac{1}{2}\nabla_{ij} g(0)\,x_i\,x_j + O(|x|^3)\\
\nabla_{ij} g(x)
	&= \nabla_{ij} g(0) + O(|x|),\\
\nabla_i g(x)\,\nabla_j g(x)
	&= \nabla_{ip}g(0)\,x_p\,\nabla_{jq}g(0)\,x_q + O(|x|^3),
\end{align*}
where the error terms are meant for small $|x|$. Hence
\[
A_{ij}(x)\,a_i\,a_j
	= \big[\nabla_{ij} g(0)\,x_i\,x_j\big]\big[\nabla_{ij} g(0)\,a_i\,a_j\big] - \big[\nabla_{ip}g(0)\,x_p\,a_i\big]\big[\nabla_{jq}g(0)\,x_q\,a_j\big] + O(|x|^3\,|a|^2).
\]
Now observe that $\nabla^2 g(0)$ is non-negative and so the difference of the first two terms on the right hand side is non-negative as well. The assertion follows.
\eproof

\noindent{\bf Proof of Lemma \ref{SimpleEllEst}.} By scaling, it suffices to consider $r$ $=$ $1$. Also, for simplicity, we will only present a proof for the case $(a_{ij})$ $=$ $(\delta_{ij})$ and $b_i$ $=$ $0$. The general case can be done in exactly the same way.

We first observe that we can assume without loss of generality that $g$ $\equiv$ $0$. Indeed, extend $g$ to a function $G$ on $B_1^+$ by setting $G(x_1,x_2,x_3)$ $=$ $g(x_1,x_2)$. Then the function $\tilde u$ $:=$ $u - G$ belongs to $W^{1,\infty}(B_1^+)$ and satisfies $-\Delta \tilde u$ $=$ $\tilde f$ where $\tilde f$ $=$ $f + \Delta G$. Moreover, there hold
\begin{align*}
|\nabla u|
	&\leq |\nabla \tilde u| + \sup_{B_1 \cap \{x_3 = 0\}} |\nabla_T g|,\\
|\tilde f|
	&\leq |f| + \sup_{B_1 \cap \{x_3 = 0\}} |\nabla_T^2 g|.
\end{align*}
Hence, if the assertion holds for $g$ $\equiv$ $0$, we can apply it to $\tilde u$ and then use the above inequalities to recover the general case. We therefore assume henceforth that $g$ $\equiv$ $0$. 

Fix $0$ $<$ $r_1$ $\leq$ $r_2$ $<$ $1$ and select a smooth cut-off function $\eta$ which is identically $1$ in $B_{r_1}$ and vanishes outside $B_{r_2}$. It is readily seen that the function $\hat u$ $:=$ $\eta u$ belongs to $W^{1,\infty}(B_1^+) \cap W^{1,1}_0(B_1^+)$, vanishes near $\partial B_1 \cap \{x_n > 0\}$ and satisfies $-\Delta \hat u$ $=$ $\hat f$ where $\hat f$ $=$ $\eta\,f - 2\nabla\eta\,\nabla u - u\,\Delta\eta$ $\in$ $L^p(B_1^+)$. Hence, by the unique solvability in $W^{2,p}$ and by $W^{2,p}$-estimates (see \cite[Theorems 9.13, 9.15]{GT}), $\hat u$ $\in$ $W^{2,p}(B_{1}^+) \cap W^{1,p}_0(B_{1}^+)$ and
\begin{align*}
\|\nabla^2 \hat u\|_{L^p(B_{r_1}^+)}
	&\leq C\,\|\hat f\|_{L^p(B_1^+)}\\
	&\leq C\Big[\frac{1}{(r_2 - r_1)^2}\|u\|_{L^p(B_{r_2}^+)} + \frac{1}{r_2 - r_1}\|\nabla u\|_{L^p(B_{r}^+)} + \|f\|_{L^p(B_1^+)}\Big].
\end{align*}
Thus, since $u|_{B_1 \cap \{x_3 = 0\}}$ $\equiv$ $0$, the Poincar\'{e}'s inequality gives
\[
\|\nabla^2 \hat u\|_{L^p(B_{1}^+)}
	\leq C\Big[\frac{r_2}{(r_2 - r_1)^2}\|\nabla u\|_{L^p(B_{r_2}^+)} + \|f\|_{L^p(B_1^+)}\Big].
\]
Applying Morrey's inequality (see \cite[Theorem 7.19]{GT}), we hence get
\begin{align*}
\mathop\textrm{osc}_{B_{r_1}^+}|\nabla u|
	&= \mathop\textrm{osc}_{B_{r_1}^+}|\nabla \hat u| \leq C\,r_1^{1-3/p}\|\nabla^2 \hat u\|_{L^p(B_{r_1}^+)}\\
	&\leq C\Big[\frac{r_2^{2 - 3/p}}{(r_2 - r_1)^2}\|\nabla u\|_{L^p(B_{r_2}^+)} + \|f\|_{L^p(B_1^+)}\Big],
\end{align*}
and so
\begin{align*}
\sup_{B_{r_1}^+}|\nabla u|
	&\leq \mathop\textrm{osc}_{B_{r_1}^+}|\nabla u| + C\,r_1^{-3/p}\|\nabla u\|_{L^p(B_{r_1}^+)}\\
	&\leq C\Big[\frac{r_2^{2 - 3/p}}{(r_2 - r_1)^2}\|\nabla u\|_{L^p(B_{r_2}^+)} + \|f\|_{L^p(B_1^+)}\Big].
\end{align*}
Using Young's inequality, we hence arrive at
\begin{equation}
\sup_{B_{r_1}^+}|\nabla u|
	\leq \frac{1}{2}\sup_{B_{r_2}^+}|\nabla u| + C_1\Big[\frac{1}{(r_2 - r_1)^p}\|\nabla u\|_{L^2(B_{1}^+)} + \|f\|_{L^p(B_1^+)}\Big].
\label{SEE::Eqn01}
\end{equation}

To derive the conclusion from \eqref{SEE::Eqn01}, we use standard iteration procedure (see e.g. \cite[Lemma 6.1]{Giu}). Fix some $\lambda$ $>$ $1$ and set $r_i$ $=$ $\frac{3}{4} - \frac{1}{4\lambda^i}$. Applying \eqref{SEE::Eqn01} repeatedly, we get 
\begin{align*}
\sup_{B_{1/2}^+}|\nabla u|
	&= \sup_{B_{r_0}^+}|\nabla u|\\
	&\leq \frac{1}{2}\sup_{B_{r_{1}}^+}|\nabla u| + C_1\Big[\frac{4^p}{(\lambda-1)^p}\|\nabla u\|_{L^2(B_{1}^+)} + \|f\|_{L^p(B_1^+)}\Big]\lambda^p\\
	&\leq \frac{1}{4}\sup_{B_{r_{2}}^+}|\nabla u| + C_1\Big[\frac{4^p}{(\lambda-1)^p}\|\nabla u\|_{L^2(B_{1}^+)} + \|f\|_{L^p(B_1^+)}\Big]\big\{\lambda^p + 2^{-1}\lambda^{2p}\big\}\\
	&\leq \ldots \leq \frac{1}{2^i}\sup_{B_{r_{i}}^+}|\nabla u| + C_1\Big[\frac{4^p}{(\lambda-1)^p}\|\nabla u\|_{L^2(B_{1}^+)} + \|f\|_{L^p(B_1^+)}\Big]\big\{\lambda^p + \ldots + 2^{1-i}\lambda^{ip}\big\}.
\end{align*}
Choosing $\lambda$ $<$ $2^{1/p}$ and letting $i$ $\rightarrow$ $\infty$, we thus get
\begin{align*}
\sup_{B_{1/2}^+}|\nabla u|
	&\leq C_2\Big[\|\nabla u\|_{L^2(B_{1}^+)} + \|f\|_{L^p(B_1^+)}\Big],
\end{align*}
which completes the proof.
\eproof

The next lemma is an extension of \cite[Lemma 7]{Ma-Za} to cover boundary estimates.
\begin{lemma}\label{BochnerEst}
There exist $\epsilon_0$ $>$ $0$, $\delta_0$ $>$ $0$ and $C$ $>$ $0$ such that if $Q_L$ satisfies $\dist(Q_L,\mcS_*)$ $\leq$ $\delta_0$ in $B_{2r_0}(x_0) \cap \Omega$ for some $x_0$ $\in$ $\bar \Omega$ and
\[
\bar E := \sup_{x \in B_{r_0}(x_0) \cap \Omega} \sup_{0 < r < r_0}\frac{1}{r}\int_{B_r(x) \cap \Omega} \Big[\frac{1}{2}|\nabla Q_L|^2 + \frac{1}{L}\tilde f_B(Q_L)\Big]\,dy \leq \epsilon_0,
\]
then
\[
\sup_{B_{r_0/2}(x_0) \cap \Omega} |\nabla Q_L| \leq C\Big[\sup_{\partial\Omega}\big[|\nabla_T Q_b| + |\nabla_T^2 Q_b|\big]  + \frac{1}{r_0}\|\dist(Q_L,\mcS_*)\|_{L^\infty(B_{r_0}(x_0))} + \frac{1}{r_0}\bar E^{1/2}\Big].
\]
Here $\nabla_T$ denotes the tangential derivative along $\partial\Omega$.
\end{lemma}

\bproof For simplicity, we write $Q$ $=$ $Q_L$ and define
\begin{align*}
&v = \frac{1}{2}|\nabla Q|^2 + \frac{1}{L}\tilde f_B(Q), B = \sup_{\partial\Omega} \big[|\nabla_T Q_b| + |\nabla_T^2 Q_b|\big],\\
&\bar d = \sup_{B_{r_0}(x_0) \cap \Omega} \dist(Q,\mcS_*) \text{ and } R = \bar E^{1/2}.
\end{align*}
By \cite[Lemma 6]{Ma-Za}, for $\delta_1$ $>$ $0$ sufficiently small, the following Bochner-type inequality holds
\begin{equation}
-\Delta v \leq C\,v^2 \text{ in } B_{2r_0}(x_0) \cap \Omega \text{ provided } \delta_0 < \delta_1, \text{ which we will assume}.
\label{BE::Eqn01}
\end{equation}

The proof then proceeds by a standard Bochner-type argument. It suffices to establish
\begin{equation}
\sup_{B_{r_0/2}(x_0) \cap \Omega} v \leq C\Big[B^2 + \frac{1}{r_0^2}\bar d^2 +  \frac{1}{r_0^2}\,R^2\Big].
\label{BE::Eqn02}
\end{equation}
To this end, we show for appropriate choice of $\epsilon_0$, $\delta_0$ and $r_0$ that
\begin{equation}
M := \max_{0 \leq r \leq r_0}(r_0 - r)^2\sup_{B_r(x_0) \cap \Omega} \big[v - C_*\,B^2\big] \leq C\,(R + \bar d)^2.
\label{BE::Eqn03}
\end{equation}
where $C_*$ is some positive constant to be determined.
Evidently, \eqref{BE::Eqn03} implies \eqref{BE::Eqn02}.

Pick $r_1$ $\in$ $[0,r_0]$ and $x_1$ $\in$ $\overline{B_{r_1}(x) \cap \Omega}$ such that
\begin{align*}
(r_0 - r_1)^2\sup_{B_{r_1}(x_0) \cap \Omega} \big[v - C_*\,B^2\big]
	&= \max_{0 \leq r \leq r_0}(r_0 - r)^2\sup_{B_r(x_0) \cap \Omega} \big[v - C_*\,B^2\big],\\
V
	&:= v(x_1) = \sup_{B_{r_1}(x_0) \cap \Omega} v.
\end{align*}
Then for $r_2$ $=$ $\frac{1}{2}(r_0 - r_1)$,
\begin{align}
\sup_{B_{r_2}(x_1) \cap \Omega} v
	&\leq \sup_{B_{r_2 + r_1}(x_0)} v \leq 4V - 3C_*\,B^2 \leq 4V,\label{BE::Eqn04}\\
M 
	&= 4r_2^2\,\big[V - C_*\,B^2\big].\label{BE::Eqn05}
\end{align}

We do the following rescaling:
\begin{align*}
\tilde \Omega
	&= \big\{y \in \RR^3: x_1 + V^{-1/2}\,y \in \Omega\big\},\\
\rho_0
	&= V^{1/2}\,r_2 \geq \sqrt{M}/2,\\
\tilde v(y)
	&= \frac{1}{V}\,v(x_1 + V^{-1/2}\,y) \text{ for } y \in B_{\rho_0}(0) \cap \tilde\Omega.
\end{align*}
Then by \eqref{BE::Eqn01}, \eqref{BE::Eqn04} and our hypotheses, we have
\begin{align}
-\Delta \tilde v
	&\leq C_0\,\tilde v \text{ in } B_{\rho_0}(0) \cap \tilde\Omega,\label{BE::Eqn06}\\
\sup_{B_{\rho_0}(0) \cap \tilde\Omega} \tilde v
	&\leq 4\,\tilde v(0) = 4,\label{BE::Eqn07}\\
\frac{1}{\rho}\int_{B_\rho(0) \cap \tilde\Omega}\,\tilde v\,dy
	&\leq R^2 \leq \epsilon_0 \text{ for } 0 < \rho \leq \rho_0. \label{BE::Eqn08}
\end{align}

For simplicity, in the sequel, we will write $B_r$ for $B_r(0)$.

Let $\rho_1$ $=$ $\min\{\dist(0,\partial\tilde\Omega),\frac{1}{10}\rho_0\}$. We first claim that $\rho_1$ $\leq$ $1$ for $\epsilon_0$ sufficiently small. For otherwise, by the Harnack inequality for \eqref{BE::Eqn06} and by \eqref{BE::Eqn08} apply to the ball $B_1(0)$, we have
\[
1 = \tilde v(0) \leq C\int_{B_1(0)} \tilde v(y)\,dy \leq C\,R^2 \leq C\,\epsilon_0,
\]
which is impossible for $\epsilon_0$ small. The claim thus follows. Now let $\hat v(z)$ $=$ $\tilde v(\rho_1\,z)$ and apply the Harnack inequality for $\hat v$ on $B_1$ to get
\[
1 = \tilde v(0) = \hat v(0) \leq C\int_{B_1} \hat v(z)\,dz = \frac{C}{\rho_1^3}\int_{B_{\rho_1}} v(y)\,dy \leq \frac{C\,R^2}{\rho_1^2}.
\]
Hence
\[
\min\big\{\dist(0,\partial\tilde\Omega),\frac{1}{10}\rho_0\big\} = \rho_1 \leq C_1\,R.
\]
If $\rho_1$ $=$ $\frac{1}{10}\rho_0$, this implies that $M$ $\leq$ $4\rho_0^2$ $\leq$ $400C_1^2\,R^2$, which proves \eqref{BE::Eqn03} and we are done. Hence, we assume that
\begin{equation}
\dist(0,\partial\tilde\Omega) = \rho_1 < \frac{1}{10}\rho_0.
\label{BE::Eqn08bis}
\end{equation}
Note that we can assume in addition that
\begin{equation}
V \geq C_*\,B^2,
\label{BE::Eqn08Add}
\end{equation}
for otherwise \eqref{BE::Eqn03} follows from \eqref{BE::Eqn05} and we are also done.

Next, we note that if we define $\tilde Q(y)$ $=$ $Q(x_1 + V^{-1/2}\,y)$ for $y$ $\in$ $B_{\rho_0} \cap \tilde\Omega$ and $\tilde L$ $=$ $V\,L$, then
\begin{equation}
\tilde v = \frac{1}{2}|\nabla \tilde Q|^2 + \frac{1}{\tilde L}\,\tilde f_B(\tilde Q).
\label{BE::Eqn09}
\end{equation}
and so, by \eqref{BE::Eqn07},
\begin{equation}
|\nabla\tilde Q| \leq 4 \text{ in } B_{\rho_0} \cap \tilde \Omega.
\label{BE::Eqn10}
\end{equation}
Hence, by \eqref{BE::Eqn08bis} and Lemma \ref{C1BndAuxRslt}, we have for $\rho$ $\in$ $(0,\rho_0)$ that
\begin{equation}
\sup_{\partial\tilde\Omega \cap B_{\rho/2}} |\nabla(\tilde Q - \tilde Q^\sharp)|
	\leq C\Big[\rho^{-1}\|\tilde Q - \tilde Q^\sharp\|_{L^\infty(B_\rho \cap \tilde\Omega)} + \rho^{-3/2}\|\nabla \tilde Q\|_{L^2(B_{\rho} \cap \tilde\Omega)} + \rho^{1/4}\|\nabla \tilde Q\|_{L^2(B_{\rho} \cap \tilde\Omega)}^{1/2}\Big].
\label{BE::Eqn11}
\end{equation}
Here $\tilde Q^\sharp$ denotes the orthogonal projection of $\tilde Q$ onto $\mcS_*$. Also, note that, as $\tilde Q^\sharp$ $=$ $P \circ \tilde Q$ for some smooht projection map $P$,
\[
|\nabla \tilde Q^\sharp| \leq C\,|\nabla \tilde Q|. 
\]
Thus, by Proposition \ref{ProjEqn}, 
\[
|\Delta \tilde Q^\sharp| \leq C\,|\nabla \tilde Q|^2.
\]
Noting \eqref{BE::Eqn08bis} and $\rho\,V^{-1/2}$ $\leq$ $\rho_0\,V^{-1/2}$ $=$ $r_2$ $\leq$ $r_0$, we can apply Lemma \ref{SimpleEllEst} to get
\begin{align}
\sup_{B_{\rho/2} \cap \tilde\Omega} |\nabla \tilde Q^\sharp|
	&\leq C\Big[\|\nabla_T\tilde Q\|_{L^\infty(\partial\tilde\Omega \cap B_{\rho})} + \rho\,\|\nabla_T^2 \tilde Q\|_{L^\infty(\partial\tilde\Omega \cap B_{\rho})}\nonumber\\
		&\qquad\qquad + \rho^{-3/2}\|\nabla \tilde Q\|_{L^2(B_\rho \cap \tilde\Omega)} + \rho^{1/4}\|\Delta \tilde Q^\sharp\|_{L^4(B_{\rho} \cap \tilde\Omega)}\Big]\nonumber\\
	&\leq C\Big[V^{-1/2}B + \rho^{-3/2}\|\nabla \tilde Q\|_{L^2(B_{\rho} \cap \tilde\Omega)} + \rho^{1/4}\|\nabla \tilde Q\|_{L^2(B_{\rho} \cap \tilde\Omega)}^{1/2}\Big].
\label{BE::Eqn12}
\end{align}
Summing up \eqref{BE::Eqn11} and \eqref{BE::Eqn12} and using \eqref{BE::Eqn08}, \eqref{BE::Eqn08Add} and \eqref{BE::Eqn09}, we infer that
\[
\sup_{\partial\tilde\Omega \cap B_{\rho/2}} \tilde v \leq C_2\Big[\frac{1}{C_*} + \frac{R^2 + \bar d^2}{\rho^2} + R\,\rho\Big].
\]
We therefore conclude that the function
\[
w_\rho := \left\{\begin{array}{l}
\max\{C_2\big[\frac{1}{C_*} + \frac{R^2 + \bar d^2}{\rho^2} + R\,\rho\big], \tilde v\} \text{ in } B_{\rho/2} \cap \tilde\Omega,\\
C_2\big[\frac{1}{C_*} + \frac{R^2  + \bar d^2}{\rho^2} + R\,\rho\big] \text{ in } B_{\rho/2} \setminus \tilde\Omega
\end{array}\right.
\]
is Lipschitz in $B_{1/2}$ and satisfies
\[
-\Delta w_\rho \leq C\,w_\rho \text{ in } B_{\rho/2}.
\]

As before, the Harnack inequality implies that $\rho_0$ $<$ $1$ if we chose $\epsilon_0$ sufficiently small and $C_*$ sufficiently large. Now, define $\hat w_\rho(z)$ $=$ $w_\rho(\rho\,z)$ and apply the Harnack inequality together with \eqref{BE::Eqn08} to get
\begin{multline*}
1 = \tilde v \leq w_\rho(0) = \hat w_\rho(0) \leq C\int_{B_{1/2}} \hat w_\rho(z)\,dz\\
	= \frac{C}{\rho^3}\int_{B_{\rho/2}} w_\rho(y)\,dy \leq C\Big[\frac{1}{C_*} + \frac{R^2  + \bar d^2}{\rho^2} + R\,\rho\Big].
\end{multline*}
To sum up we have shown that
\begin{equation}
1 \leq C_3\Big[\frac{1}{C_*} + \frac{(R + \bar d)}{\rho^2} + (R + \bar d)\,\rho\Big] \text{ for any } \rho \in (0,\rho_0).
\label{BE::Eqn13}
\end{equation}
Now select $C_*$ $>$ $4C_3$, $\epsilon_0$ $<$ $32^{-1}\,C_3^{-3/2}$ and $\delta_0$ $<$ $\epsilon_0^{1/2}$ so that
\[
C_3\Big[\frac{1}{C_*} + \frac{(R + \bar d)^2}{\rho_*^2} + (R + \bar d)\,\rho_*\Big] < \frac{1}{2} \text{ for } \rho_* = (R + \bar d)^{1/3}.
\]
\eqref{BE::Eqn13} then implies that
\[
\rho_0 < \rho_* \leq 2\,\epsilon_0^{1/6},
\]
and so
\[
1 \leq C_3\Big[\frac{1}{C_*} + \frac{(R + \bar d)^2}{\rho_0^2} + \epsilon_0^{2/3}\Big] \leq \frac{1}{2} + C_3\,\frac{(R + \bar d)^2}{\rho_0^2}.
\]
Therefore, by \eqref{BE::Eqn05},
\[
M \leq \rho_0^2 \leq 2C_3\,(R + \bar d)^2.
\]
The proof is complete.
\eproof

\begin{remark}\label{BochnerEst::Interior}
Regarding interior estimate, the proof given above yields the following statement. There exist $\epsilon_0$ $>$ $0$, $\delta_0$ $>$ $0$ and $C$ $>$ $0$ such that if $Q_L$ satisfies $\dist(Q_L,\mcS_*)$ $\leq$ $\delta_0$ in $B_{2r_0}(x_0)$ $\Subset$ $\Omega$ and
\[
\bar E := \sup_{x \in B_{r_0}(x_0)} \sup_{0 < r < r_0}\frac{1}{r}\int_{B_r(x)} \Big[\frac{1}{2}|\nabla Q_L|^2 + \frac{1}{L}\tilde f_B(Q_L)\Big]\,dy \leq \epsilon_0,
\]
then
\[
\sup_{B_{r_0/2}(x_0)} |\nabla Q_L| \leq \frac{C}{r_0}\bar E^{1/2}.
\]
\end{remark}

\begin{proposition}\label{UniC1Bnd}
Assume that $Q_{L_k}$ converges strongly in $H^1(\Omega, \mcS_0)$ to a minimizer $Q_*$ $\in$ $H^1(\Omega,\mcS_*)$ of $I_*$ for some sequence $L_k$ $\rightarrow$ $0$. For any compact subset $K$ of $\bar\Omega\setminus \Sing(Q_*)$, there exists $\bar L$ $=$ $\bar L(a^2,b^2,c^2,\Omega,K,Q_b,Q_*)$ $>$ $0$ such that
\[
\sup_K |\nabla Q_{L_k}| \leq C(a^2,b^2,c^2,\Omega,K,Q_b,Q_*) \text{ for any } L_k \leq \bar L.
\]
\end{proposition}

\bproof It suffices to consider $K$ $=$ $K_{2\eta}$ $:=$ $\{x \in \bar\Omega: \dist(x,\Sing(Q_*)) \geq 2\eta\}$ where $\eta$ is some arbitrary small positive number. Let $\epsilon_0$, $\delta_0$ and $r_0$ be as in Lemma \ref{BochnerEst}.

First, by \cite[Propositions 4 and 6]{Ma-Za} and Lemma \ref{DistComp}, we can select $\bar L$ such that
\[
\dist(Q_{L_k},\mcS_*) \leq C\sqrt{\tilde f_B(Q_{L_k})} \leq \delta_0 \text{ for any } x \in K \text{ and } L_k \leq \bar L.
\]
where here and below $C$ denotes a constant that may depend on $\Omega$, $a^2$, $b^2$, $c^2$, and $Q_b$ but is otherwise independent of $L_k$ and $Q_{L_k}$.

By standard regularity results for harmonic maps (see e.g. \cite{su} and \cite{SU-Bdry}), for some $\epsilon_1$ $>$ $0$ to be chosen, there exists $0$ $<$ $r_0$ $<$ $\frac{1}{4}\eta$ such that
\begin{equation}
\frac{1}{r_0}\int_{B_{2r_0}(x) \cap \Omega}|\nabla Q_*|^2\,dx \leq \epsilon_1 \text{ for any }x \in K_\eta.
\label{UniC1Bnd::Eqn01}
\end{equation}
Also, note that $Q_{L_k} - Q_*$ $\in$ $H^1_0(\Omega, \mcS_0)$. Hence since $Q_{L_k}$ is $\tilde I_{L_k}$-minimizing,
\[
\int_{\Omega} \Big[|\nabla Q_{L_k}|^2 + \frac{1}{L_k}\tilde f_B(Q_{L_k})\Big]\,dx \leq \int_{\Omega} |\nabla Q_*|^2\,dx.
\]
On the other hand, since $Q_{L_k}$ $\rightarrow$ $Q_*$ in $H^1$, we have
\begin{equation}
\int_\Omega |\nabla Q_*|^2\,dx = \lim_{k \rightarrow \infty}\int_{\Omega} |\nabla Q_{L_k}|^2\,dx.
\label{UniC1Bnd::Eqn02}
\end{equation}
It follows that
\begin{equation}
\lim_{k \rightarrow 0} \frac{1}{L_k}\int_\Omega \tilde f_B(Q_{L_k})\,dx = 0.
\label{UniC1Bnd::Eqn03}
\end{equation}
From \eqref{UniC1Bnd::Eqn01}-\eqref{UniC1Bnd::Eqn03}, we infer that there exists some $\bar L$ $>$ $0$ such that
\[
\frac{1}{r_0}\int_{B_{2r_0}(x) \cap \Omega} e_{L_k}\,dx \leq \epsilon_1 \text{ for any } x \in K \text{ and } L_k \leq \bar L,
\]
where
\[
e_{L_k} = \frac{1}{2}|\nabla Q_{L_k}|^2 + \frac{1}{L_k}\tilde f_B(Q_{L_k}).
\]
Applying the monotonicity formulas in \cite[Lemmas 2,9]{Ma-Za}, we arrive at
\[
\frac{1}{r}\int_{B_{r}(x) \cap \Omega} e_k\,dx \leq C\,\epsilon_1 \text{ for any } x \in K, r \in (0,r_0] \text{ and } L_k \leq \bar L.
\]
We now fix $\epsilon_1$ so that $C\epsilon_1$ $\leq$ $\epsilon_0$. The assertion follows from Lemma \ref{BochnerEst} with the appropriate choice of $\bar L$ so that the above argument goes through.
\eproof

\noindent{\bf Proof of Proposition \ref{C1AlphaConv}.} Fix $K$ $\Subset$ $K'$ $\Subset$ $\bar\Omega\setminus\Sing(Q_*)$. By Proposition \ref{UniC1Bnd} and Corollary \eqref{UniLapBnd}, there exists $\bar L$ $>$ $0$ and $C$ $>$ $0$ depending only on $a^2$, $b^2$, $c^2$, $\Omega$, $K$, $K'$, $Q_b$ and $Q_*$ such that
\[
\sup_{K'}\big[|\nabla Q_{L_k}| + |\Delta Q_{L_k}|\big]  \leq C  \text{ for } L_k \leq \bar L.
\]
Also, by \cite[Proposition 3]{Ma-Za},
\[
\sup_\Omega |Q_{L_k}| \leq \sqrt{\frac{2}{3}}s_+.
\]
The conclusion follows from standard $W^{2,p}$-estimates for Poisson equations and Morrey's inequality (see e.g. \cite{GT}).
\eproof

\section{$C^{j}$-convergence}\label{SEC-CjConv}

In the previous section, we showed that the $H^1$-convergence of a sequence of minimizers $Q_{L_k}$ $\in$ $H^1(\Omega, \mcS_0)$ of $I_{L_k}$ to a minimizer $Q_*$ $\in$ $H^1(\Omega,\mcS_*)$ of $I_*$ improves itself to a $C^{1,\alpha}$-convergence on compact subsets of $\bar\Omega\setminus \Sing(Q_*)$. In this section, we study $C^j$-convergence. Like in Ginzburg-Landau theory, we do not expect to have $C^2$ convergence up to the boundary; for $\Delta Q_L$ vanishes on $\partial \Omega$ while $\Delta Q_*$ needs not do.

\begin{proposition}\label{CjConv}
Assume that $Q_{L_k}$ converges strongly in $H^1(\Omega, \mcS_0)$ to a minimizer $Q_*$ $\in$ $H^1(\Omega,\mcS_*)$ of $I_*$ for some sequence $L_k$ $\rightarrow$ $0$. Then $Q_{L_k}$ converges uniformly in $C^j$-norm on compact subsets of $\Omega\setminus\Sing(Q_*)$ to $Q_*$.
\end{proposition}

The approach we take closely follows \cite{BBH}. A key step is to show the convergence of the ``normal component'', i.e. the ``minimal polynomial'' $Q_L^2 - \frac{1}{3}s_+\,Q_L - \frac{2}{9}s_+^2\,\IdMa$. It is clear that Proposition \ref{CjConv} follows from Propositions \ref{MinPolyConv}, \ref{UniC1Bnd} and the following lemma.

\begin{lemma}\label{PrepUniHighDerEst}
Assume for some $B_r(x)$ $\Subset$ $\Omega$, $L$ $\in$ $(0,1)$ and $C_1$ $>$ $0$ that
\[
\sup_{B_r(x)} \Big[|\nabla Q_L| + |X_L|\Big]
	\leq C_1.
\]
Then, for $0$ $<$ $s$ $<$ $r$ and $j$ $\geq$ $0$,
\begin{align}
\sup_{B_s(x)}\big|\nabla^j X_L\big|
	&\leq C(a^2,b^2,c^2,C_1,j)[1 + (r-s)^{-j}],\label{PUHDE::Est01}\\
\sup_{B_s(x)}|\nabla^{j+1} Q_{L}|
	&\leq C(a^2,b^2,c^2,C_1,j)[1 + (r-s)^{-j}].\label{PUHDE::Est02}
\end{align}
\end{lemma}

In the proof of the above lemma, obtaining \eqref{PUHDE::Est01} is central. As we have said before, $X_L$ ``measures'' how fast $Q_L$ approaches the limit manifold. Hence, intuitively speaking, $X_L$ should ``behave'' similarly to a vector in the normal space $(T_{Q_*} \mcS_*)^\perp$. Now, consider the linear map $\Lambda$ on $(T_{Q_*} \mcS_*)^\perp$ created by a left multiplcation by $Q_*$. This map has two eigenvalues $-\frac{1}{3}s_+$ and $\frac{2}{3}s_+$, where the former is double and the latter is simple. (See the proof of Lemma \ref{MysIds}.) The eigenspace decomposition of $(T_{Q_*} \mcS_*)^\perp$ with respect to $\Lambda$ is given by
\[
(T_{Q_*} \mcS_*)^\perp = \Big[\big(Q_* + \frac{1}{3}s_+\,\IdMa\big)(T_{Q_*} \mcS_*)^\perp\Big] \oplus \Big[\big(Q_* - \frac{2}{3}s_+\,\IdMa\big)(T_{Q_*} \mcS_*)^\perp\Big].
\]
We will momentarily see that this decomposition is very useful in the analysis of $X_L$. (See Steps 4 and 5 of the following proof.)

\bigskip
\bproof Like usual, we drop the subscript $L$ in $Q_{L}$ and $X_L$. Also, we will abbreviate $B_r(x)$ to $B_r$. Throughout the proof we will frequently use the estimate $|Q|$ $\leq$ $\sqrt{\frac{2}{3}}s_+$ without explicitly mentioning it (see \cite[Proposition 3]{Ma-Za}).

For $j \geq 0$ and $s > 0$, define
\begin{align*}
m(s,j)
	&= \sup_{B_{s}}\sum_{l=0}^j\big|\nabla^l X| + \sum_{l=0}^{j+1}|\nabla^l Q|,\\
M(s,j)
	&= m(s,j) + \sum_{l_1 + \ldots + l_p = j+1,\, l_q > 0,\,q > 1}m(s,l_1)\ldots m(s,l_p).
\end{align*}

Fix $\epsilon$ $\in$ $(0,1)$ for the moment. In the sequel, $C$ denotes a constant that depends only on $a^2$, $b^2$, $c^2$, $C_1$, $j$ and $\epsilon$. We will consecutively show the following six estimates for $0$ $<$ $r_1$ $<$ $r_2$ $<$ $r$.
\begin{align}
\|\nabla^{j+2} Q\|_{L^\infty(B_{r_1})}
	&\leq \epsilon\,\|\nabla^{j+1} X\|_{L^\infty(B_{r_2})} + C\Big[1 + \frac{1}{r_2 - r_1}\Big]M(r_2,j),
\label{PUHDE::Step1-2}\\
\|\nabla^{j+1} X\|_{L^\infty(B_{r_1})}
	&\leq C\Big[\frac{1}{\sqrt{L}} + \frac{1}{r_2 - r_1}\Big]M(r_2,j),
\label{PUHDE::Step2}\\
\|\nabla^{j+1}\tr(X)\|_{L^\infty(B_{r_1})}
	&\leq \epsilon\,\|\nabla^{j+1} X\|_{L^\infty(B_{r_2})} + C\Big[1 + \frac{1}{r_2 - r_1}\Big]M(r_2,j),
\label{PUHDE::Step3}\\
\big\|\nabla^{j+1} Y^{(-2)}\big\|_{L^\infty(B_{r_1})}
	&\leq \epsilon\,\|\nabla^{j+1} X\|_{L^\infty(B_{r_2})} + C\Big[1 + \frac{1}{r_2 - r_1}\Big]M(r_2,j),
\label{PUHDE::Step4}\\
\|\nabla^{j+1}X\|_{L^\infty(B_{r_1})}
	&\leq \epsilon\,\|\nabla^{j+1} X\|_{L^\infty(B_{r_2})} + C\Big[1 + \frac{1}{r_2 - r_1}\Big]M(r_2,j),
\label{PUHDE::Step5}\\
m(r_1,j+1)
	&\leq C\Big[1 + \frac{1}{r_2 - r_1}\Big]M(r_2,j),
\label{PUHDE::Step6}
\end{align}
where in \eqref{PUHDE::Step4}, $Y^{(-2)}$ $=$ $\big(Q - \frac{2}{3}s_+\,\IdMa\big)X$. It is readily seen that the conclusion follows from \eqref{PUHDE::Step6}. Among the above, estimate \eqref{PUHDE::Step4} is crucial. 

\medskip
\noindent\underline{Step 1:} Proof of \eqref{PUHDE::Step1-2}. By \eqref{L-DG::ELEqn}, we have
\begin{align}
\Delta Q
	&= \frac{1}{L}\Big[-a^2\,Q - b^2[Q^2 - \frac{1}{3}\tr(Q^2)\,\IdMa] + c^2\,\tr(Q^2)\,Q\Big]\nonumber\\
	&= c^2\,\tr(X)\,Q + \frac{1}{3}b^2\,\tr(X)\,\IdMa - b^2\,X.
\label{PUHDE::Eqn01}
\end{align}
Differentiating \eqref{PUHDE::Eqn01} up to $j+1$ order, we get
\[
\|\Delta (\nabla^{j+1} Q)\|_{L^\infty(B_{r_2})} \leq C\,\|\nabla^{j+1} X\|_{L^\infty(B_{r_2})} + M(r_2,j).
\]
Therefore, by a standard interpolation inequality (see e.g. \cite[Lemma A.1]{BBH}),
\begin{align*}
\|\nabla^{j+2} Q\|_{L^\infty(B_{r_1})}^2
	&\leq C\|\nabla^{j+1}Q\|_{L^\infty(B_{r_2})}\big[\|\Delta (\nabla^{j+1} Q)\|_{L^\infty(B_{r_2})} + \frac{1}{(r_2 - r_1)^2}\|\nabla^{j+1}Q\|_{L^\infty(B_{r_2})}\big]\\
	&\leq C\,M(r_2,j)\big[\|\nabla^{j+1} X\|_{L^\infty(B_{r_2})} + [1 + (r_2 - r_1)^{-2}]M(r_2,j)\big].
\end{align*}
This implies \eqref{PUHDE::Step1-2}.

\medskip
\noindent\underline{Step 2:} Proof of \eqref{PUHDE::Step2}. Using \eqref{PUHDE::Eqn01}, we calculate
\begin{multline}
L\,\Delta X = 2Q\Big[c^2\,\tr(X)\,Q + \frac{1}{3}b^2\,\tr(X)\,\IdMa - b^2\,X\Big]\\
		- \frac{1}{3}s_+\Big[c^2\,\tr(X)\,Q + \frac{1}{3}b^2\,\tr(X)\,\IdMa - b^2\,X\Big] + 2\sum_{\alpha = 1}^3(\nabla_\alpha Q)^2.
\label{PUHDE::Eqn02}
\end{multline}
It thus follows from a standard interpolation inequality (see e.g. \cite[Lemma A.1]{BBH}) that
\begin{align*}
\|\nabla^{j+1} X\|_{L^\infty(B_{r_1})}^2
	&\leq C\|\nabla^j X\|_{L^\infty(B_{r_2})}\big[\|\Delta \nabla^j X\|_{L^\infty(B_{r_2})} + \frac{1}{(r_2 - r_1)^2}\|\nabla^j X\|_{L^\infty(B_{r_2})} \big]\\
	&\leq C\big[L^{-1} + (r_2 - r_1)^{-2}\big]M(r_2,j)^2,
\end{align*}
which proves \eqref{PUHDE::Step2}.

\medskip
\noindent\underline{Step 3:} Proof of \eqref{PUHDE::Step3}. Taking trace of \eqref{PUHDE::Eqn02}, we get
\begin{equation}
L\,\Delta \tr(X) = 2c^2\,\tr(X)\,\tr(Q^2) - 2b^2\,\tr(QX) + 2|\nabla Q|^2.
\label{PUHDE::Eqn03}
\end{equation}
Using
\[
\tr(Q^3) = \frac{3}{4s_+}L^2(\tr(X))^2 + \frac{1}{2}s_+\,\tr(Q^2) - \frac{1}{9}s_+^3 - \frac{3}{2s_+}L^2\,\tr(X^2),
\]
we write
\begin{align}
\tr(QX)
	&= \frac{1}{L}\Big[\tr(Q^3) - \frac{1}{3}s_+\,\tr(Q^2)\big]\nonumber\\
	&= \frac{1}{L}\Big[\frac{3}{4s_+}L^2(\tr(X))^2 + \frac{1}{6}s_+\,\tr(Q^2) - \frac{1}{9}s_+^3 - \frac{3}{2s_+}L^2\,\tr(X^2)\Big]\nonumber\\
	&= \frac{3}{4s_+}L(\tr(X))^2 + \frac{1}{6}s_+\,\tr(X) - \frac{3}{2s_+}L\,\tr(X^2).
\label{PUHDE::Eqn04}
\end{align}
We thus rewrite \eqref{PUHDE::Eqn03} as
\begin{align}
L\,\Delta \tr(X)
	&= (\frac{4}{3}c^2\,s_+^2 - \frac{1}{3}b^2s_+)\tr(X) + 2|\nabla Q|^2\nonumber\\
		&\qquad\qquad + L\big[2c^2\,(\tr(X))^2 - \frac{3}{2s_+}\,b^2\,(\tr(X))^2 + \frac{3}{s_+}\,b^2\,\tr(X^2)\big]\nonumber\\
	&= (2a^2 + \frac{1}{3}b^2s_+)\tr(X) + 2|\nabla Q|^2\nonumber\\
		&\qquad\qquad + L\big[2c^2\,(\tr(X))^2 - \frac{3}{2s_+}\,b^2\,(\tr(X))^2 + \frac{3}{s_+}\,b^2\,\tr(X^2)\big].
\label{PUHDE::Eqn04Bis}
\end{align}
We arrive at
\begin{equation}
-L\,\Delta \tr(X) + (2a^2 + \frac{1}{3}b^2s_+)\tr(X)
	= A + L B,
\label{PUHDE::Eqn05}
\end{equation}
where, by \eqref{PUHDE::Step1-2} and \eqref{PUHDE::Step2}
\begin{multline*}
\|\nabla^{j+1}A\|_{L^\infty(B_{(r_1 + r_2)/2})} + \sqrt{L}\|\nabla^{j+1}B\|_{L^\infty(B_{(r_1 + r_2)/2})}\\
	\leq \epsilon\,\|\nabla^{j+1} X\|_{L^\infty(B_{r_2})} + C[1 + (r_2 - r_1)^{-1}]M(r_2,j).
\end{multline*}
Hence, by differentiating \eqref{PUHDE::Eqn05} to $(j+1)$ order and applying Lemma \ref{BBHMaxPrin}, we infer that
\begin{align*}
\|\nabla^{j+1}\tr(X)\|_{L^\infty(B_{r_1})}
	&\leq C\,\|\nabla^{j+1}(A + LB)\|_{L^\infty(B_{(r_1 + r_2)/2})}\\
		&\qquad\qquad + C\min\Big\{1,\frac{L^{1/2}}{r_2 - r_1}\Big\}\,\|\nabla^{j+1}\tr(X)\|_{L^\infty(B_{(r_1 + r_2)/2})}\\
	&\leq \epsilon\,\|\nabla^{j+1} X\|_{L^\infty(B_{r_2})} + C[1 + (r_2 - r_1)^{-1}]M(r_2,j)\\
		&\qquad\qquad + C\,\min\Big\{1,\frac{L^{1/2}}{r_2 - r_1}\Big\}\,[L^{-1/2} + (r_2 - r_1)^{-1}]\,M(r_2,j)\\
	&\leq \epsilon\,\|\nabla^{j+1} X\|_{L^\infty(B_{r_2})} + C[1 + (r_2 - r_1)^{-1}]M(r_2,j).
\end{align*}
Estimate \eqref{PUHDE::Step3} is proved.

\medskip
\noindent\underline{Step 4:} Proof of \eqref{PUHDE::Step4}. Set
\[
Y^{(\alpha)} = \Big(Q\, + \frac{\alpha}{3}\,s_+\,\IdMa\Big)\,X
\]
where $\alpha$ is to be determined. We derive from \eqref{PUHDE::Eqn01} that
\begin{align*}
L\,\Delta Y^{(\alpha)}
	&= \Delta\big[Q^3 - \frac{1}{3}s_+(1 - \alpha)Q^2 - \frac{1}{9}s_+^2(2 + \alpha)Q\big]\nonumber\\
	&= -3b^2\,Q^2\,X + \frac{2}{3}b^2\,s_+(1 - \alpha)\,Q\,X\nonumber\\
		&\qquad\qquad+ \frac{1}{9}b^2\,s_+^2(2 + \alpha)\,X + R(Q, \nabla Q, \tr(X)),
\end{align*}
where $R$ is an explicit polynomial. Noting that
\begin{equation}
Q^2\,X = L\,X^2 + \frac{1}{3}s_+\,Q\,X + \frac{2}{9}\,s_+^2\,X,
\label{PUHDE::Eqn06Bis}
\end{equation}
we get
\begin{align*}
L\,\Delta Y^{(\alpha)}
	&= -\frac{1}{3}b^2\,s_+(1 + 2\alpha)\,Q\,X\nonumber\\
		&\qquad\qquad - \frac{1}{9}b^2\,s_+^2(4 - \alpha)\,X - 3b^2\,L\,X^2 + R(Q, \nabla Q, \tr(X)),
\end{align*}
Pick $\alpha$ $=$ $-2$, we arrive at
\begin{equation}
-L\,\Delta Y^{(-2)}_{ij} + b^2\,s_+\,Y^{(-2)}_{ij}
	= \big[3b^2\,L\,X^2 - R(Q, \nabla Q, \tr(X))\big]_{ij}.
\label{PUHDE::Eqn08}
\end{equation}
As before, by differentiating \eqref{PUHDE::Eqn08} to $j+1$ order and applying Lemma \ref{BBHMaxPrin} together with \eqref{PUHDE::Step1-2}-\eqref{PUHDE::Step3}, we get
\begin{align*}
\|\nabla^{j+1} Y^{(-2)}\|_{L^\infty(B_{r_1})}
	&\leq \epsilon\,\|\nabla^{j+1} X\|_{L^\infty(B_{r_2})} + C[1 + (r_2 - r_1)^{-1}]M(r_2,j).
\end{align*}
Recalling the definition of $Y^{(-2)}$, we arrive at \eqref{PUHDE::Step4}.

\medskip
\noindent\underline{Step 5:} Proof of \eqref{PUHDE::Step5}. We have
\[
\tr(Y^{(-2)}) = \frac{1}{L}\big[\tr(Q^3) - s_+\tr(Q^2) + \frac{4}{9}s_+^3\big] = \frac{1}{L}\big[\tr(Q^3) - \frac{2}{9}s_+^3\big] - s_+\,\tr(X).
\]
Hence, by \eqref{PUHDE::Step3} and \eqref{PUHDE::Step4},
\begin{equation}
\Big\|\nabla^{j+1}\big(\frac{1}{L}\big[\tr(Q^3) - \frac{2}{9}s_+^3\big]\big)\Big\|_{L^\infty(B_{r_1})} \leq \epsilon\,\|\nabla^{j+1} X\|_{L^\infty(B_{r_2})} + C[1 + (r_2 - r_1)^{-1}]M(r_2,j).
\label{PUHDE::Eqn09}
\end{equation}

On the other hand, as $Q$ is a traceless $3 \times 3$ matrix, Cayley's theorem gives
\[
Q^3 - \frac{1}{2}\,\tr(Q^2)\,Q - \frac{1}{3}\tr(Q^3)\,\IdMa = 0.
\]
Therefore
\begin{align*}
Y^{(1)}
	&= (Q + \frac{1}{3}s_+\,\IdMa)X = \frac{1}{L}(Q^3 - \frac{1}{3}s_+^2\,Q - \frac{2}{27}s_+^3\,\IdMa)\\
	&= \frac{1}{2}\tr(X)Q + \frac{1}{3L}\big[\tr(Q^3) - \frac{2}{9}s_+^3\big]\,\IdMa.
\end{align*}
It thus follows from \eqref{PUHDE::Step3} and \eqref{PUHDE::Eqn09} that
\[
\|\nabla^{j+1}Y^{(1)}\|_{L^\infty(B_{r_1})} \leq \epsilon\,\|\nabla^{j+1} X\|_{L^\infty(B_{r_2})} + C[1 + (r_2 - r_1)^{-1}]M(r_2,j).
\]
As $s_+X$ $=$ $Y^{(1)} - Y^{(-2)}$, we get \eqref{PUHDE::Step5}.

\medskip
\noindent\underline{Step 6:} Proof of \eqref{PUHDE::Step6}. Let
\[
\Phi(s) = \|\nabla^{j+1} X\|_{L^\infty(B_s)}, \qquad 0 \leq s \leq r.
\]
Then $\Phi$ is monotonically non-decreasing and, by \eqref{PUHDE::Step5} with $\epsilon$ $=$ $\frac{1}{2}$,
\begin{equation}
\Phi(r_1) \leq \frac{1}{2}\,\Phi(r_2) + C[1 + (r_2 - r_1)^{-1}]M(r_2,j), \qquad 0 \leq r_1 < r_2 \leq r.
\label{PUHDE::Eqn10}
\end{equation}
The proof proceeds by a standard iteration argument (cf. \cite[Lemma 6.1]{Giu}). Define a sequence $\rho_j$ by
\begin{align*}
\rho_0
	&= r_1,\\
\rho_{i+1} - \rho_i
	&= \frac{2^i}{3^{i+1}}(r_2 - r_1).
\end{align*}
By a simple induction we get from \eqref{PUHDE::Eqn10} that
\[
\Phi(r_1) \leq \frac{1}{2^i}\Phi(\rho_i) + C[1 + (r_2 - r_1)^{-1}]M(r_2,j)\sum_{l=0}^{i-1}\frac{3^l}{4^l}.
\]
Sending $i$ $\rightarrow$ $\infty$, we get
\[
\Phi(r_1) \leq C[1 + (r_2 - r_1)^{-1}]M(r_2,j).
\]
Estimate \eqref{PUHDE::Step6} follows immediately in view of \eqref{PUHDE::Step1-2}. The proof is complete.
\eproof

\begin{remark}
The proof given above can be adapted to give a different (though much longer) proof for Proposition \ref{MinPolyConv}.
\end{remark}

The following is an easy consequence of Lemma \ref{PrepUniHighDerEst}.

\begin{corollary}\label{PUHDE::Cor}
Assume for some $x$ $\in$ $\partial \Omega$, $r$ $>$ $0$, $L$ $\in$ $(0,1)$ and $C_1$ $>$ $0$ that
\[
\sup_{B_r(x) \cap \Omega} \Big[|\nabla Q_L| + |X_L|\Big]
	\leq C_1.
\]
Then, for $y$ $\in$ $B_{r/2}(x) \cap \Omega$, there holds
\begin{align}
|\nabla^j X_L|(y)
	&\leq \frac{C(a^2,b^2,c^2,C_1,j)}{\dist(y,\partial(B_r(x) \cap\Omega))^{j}},\label{PUHDE::Cor::Est01}\\
|\nabla^{j+1} Q_{L}|(y)
	&\leq \frac{C(a^2,b^2,c^2,C_1,j)}{\dist(y,\partial(B_r(x) \cap\Omega))^{j}}.\label{PUHDE::Cor::Est02}
\end{align}
\end{corollary}

So far, we have shown the first part of Theorem \ref{MainThm1}. It remains to write $\Delta Q_L$ in the form stated in the theorem. By \eqref{PUHDE::Eqn01}, this can be done if we know what the limit of $X_L$ is. By a careful inspection in the proof of Lemma \ref{PrepUniHighDerEst} (namely \eqref{PUHDE::Eqn02} and \eqref{PUHDE::Eqn04Bis}), one expects that
\begin{align*}
\tr(X_L) \text{ ``converges to'' } & - \frac{6}{6a^2 + b^2s_+}|\nabla Q_*|^2
\end{align*}
and
\begin{align*}
X_L \text{ ``converges to'' } & \frac{1}{b^2\,s_+^2}\Big\{-\frac{6s_+^2}{6a^2 + b^2s_+}|\nabla Q_*|^2\big[c^2\,Q_* + \frac{1}{3}\,b^2\,\IdMa\big]\nonumber\\
	&\qquad\qquad + 4\big(Q_* - \frac{1}{6}s_+\,\IdMa\big)\sum_{\alpha = 1}^3(\nabla_\alpha Q_*)^2\Big\}.
\end{align*}
A more natural and systematic way of getting this information is to assume that there is a formal asymptotic expansion $Q_L$ $=$ $Q_* + L\,Q_\bullet + O(L^2)$ and expand \eqref{L-DG::ELEqn} accordingly. Compare \eqref{FAE::Id01}.

Motivated on the above discussion, we consider
\begin{align}
Y_L
	&= \tr(X_L) + \frac{6}{6a^2 + b^2s_+}|\nabla Q_L|^2,\label{YDef}\\
Z_L
	&= X_L - \frac{1}{b^2\,s_+^2}\Big\{-\frac{6s_+^2}{6a^2 + b^2s_+}|\nabla Q_L|^2\big[c^2\,Q_L + \frac{1}{3}\,b^2\,\IdMa\big]\nonumber\\
	&\qquad\qquad + 4\big(Q_L - \frac{1}{6}s_+\,\IdMa\big)\sum_{\alpha = 1}^3(\nabla_\alpha Q_L)^2\Big\}.\label{ZDef}
\end{align}
The following lemma establishes the rate of convergence of $Y_L$ and $Z_L$.

\begin{lemma}\label{EqnConv-Ez}
Assume for some $x$ $\in$ $\bar\Omega$, $r$ $>$ $0$, $L$ $\in$ $(0,1)$ and $C_1$ $>$ $0$ that 
\[
\sup_{B_r(x) \cap \Omega} \big[|\nabla Q_L| + |X_L|\big] \leq C_1.
\]
Then, for $y$ $\in$ $B_{r/2}(x) \cap \Omega$, there holds
\begin{align}
|\nabla^j Y_L| + |\nabla^j Z_L| \leq \frac{C(a^2,b^2,c^2,C_1,j)}{\dist(y,\partial(B_r(x) \cap\Omega))^{j+2}}L,
\label{EC-Ez::Est01}
\end{align}
\end{lemma} 

\bproof We will drop the subscript $L$ in $Q_{L}$, $X_L$, $Y_L$ and $Z_L$.

By Corollary \ref{PUHDE::Cor} and \eqref{PUHDE::Eqn04Bis}, we have for $y$ $\in$ $B_{r/2} \cap \Omega$,
\begin{equation}
|\nabla^j Y|(y) \leq \frac{C}{\dist(y,\partial(B_r(x) \cap\Omega))^{j+2}}\,L.
\label{EC-Ez::Eqn01}
\end{equation}

To continue we introduce the following convention: $R(A_1, \ldots, A_t; A_{t+1}, \ldots, A_s)$ will be used to denote various explicitly computable polynomials in the $A_l$'s which are (jointly) linear in $(A_{t+1}, \ldots, A_s)$.

By \eqref{PUHDE::Eqn02}, we have
\begin{align}
L\,\Delta X
	&= -b^2\big(2Q - \frac{1}{3}s_+\,\IdMa\big)X + 2(\nabla Q)^2\nonumber\\
		&\qquad\qquad + \tr(X)\big[\frac{1}{3}(2b^2 + c^2\,s_+)Q  + \frac{1}{9}(-b^2\,s_+ + 4c^2\,s_+^2)\,\IdMa\big] + L\,R(X),
\label{EC-Ez::Eqn02}
\end{align}
where
\[
(\nabla Q)^2 = \sum_{\alpha = 1}^3 (\nabla_\alpha Q)^2.
\]
Multiplying \eqref{EC-Ez::Eqn02} to the left by $Q - \frac{1}{6}s_+\,\IdMa$ we get
\begin{align}
L\big(Q - \frac{1}{6}s_+\,\IdMa\big)\Delta X
	&= -\frac{1}{2}b^2\,s_+^2\,X + 2\big(Q - \frac{1}{6}s_+\,\IdMa\big)(\nabla Q)^2\nonumber\\
		&\qquad\qquad+ \tr(X)\big[\frac{1}{2}c^2\,s_+^2\,Q + \frac{1}{6}\,b^2\,s_+^2\,\IdMa\big] + L\,R(Q,X)\nonumber\\
	&= -\frac{1}{2}b^2\,s_+^2\,Z + L\,R(Q,X;\frac{1}{L}Y)).
\label{EC-Ez::Eqn03}
\end{align}
It is readily seen that the assertion follows from Corollary \ref{PUHDE::Cor}, \eqref{EC-Ez::Eqn01} and \eqref{EC-Ez::Eqn03}.
\eproof

\medskip
\noindent{\bf Proof of Theorem \ref{MainThm1}.} The strong convergence in $H^1(\Omega,\mcS_0)$ follows from \cite[Lemma 3]{Ma-Za}. (The statement therein requires an additional hypothesis that $\Omega$ be simply connected. However, such condition was used there only to obtain a lifting of a map in $H^1(\Omega,\mcS_*)$ to a map in $H^1(\Omega,\Sphere^2)$, which we do not need here.) The convergence in $C^{1,\alpha}$ and $C^j$ follows directly from Propositions \ref{C1AlphaConv} and \ref{CjConv}. 

It remains to show the last assertion. Using \eqref{PUHDE::Eqn01}, \eqref{YDef} and \eqref{ZDef}, we compute
\begin{align}
\Delta Q
	&= c^2\,\tr(X)\,Q + \frac{1}{3}b^2\,\tr(X)\,\IdMa - b^2\,X\nonumber\\
	&= - \frac{4}{s_+^2}\big(Q - \frac{1}{6}s_+\,\IdMa\big)\sum_{\alpha = 1}^3(\nabla_\alpha Q)^2 
		+ c^2\,Y\,Q + \frac{1}{3}b^2\,Y\,\IdMa - b^2\,Z.
\label{ExplicitApproxEqn}
\end{align}
The conclusion then follows from Lemma \ref{EqnConv-Ez}.
\eproof

\section{Higher derivative estimates near the boundary}\label{SEC-HDEst}

By Proposition \ref{UniC1Bnd} and Corollary \ref{UniLapBnd}, $\Delta Q_L$ is uniformly bounded up to the boundary. This suggests that estimate \eqref{PUHDE::Cor::Est02} for the derivatives of $Q_L$ might not be optimal near $\partial\Omega$. In this section, we will first study higher derivatives estimate near $\partial\Omega$ and then use it to slightly improve estimate \eqref{EC-Ez::Est01} in Lemma \ref{EqnConv-Ez}. This will be useful when we study the existence of the ``first order term'' in the asymptotic expansion of $Q_L$.

\begin{proposition}\label{BdryCjEst}
Assume for some $x$ $\in$ $\partial\Omega$, $r$ $>$ $0$, $L$ $\in$ $(0,1)$ and $C_1$ $>$ $0$ that $\Omega \cap B_r(x)$ is connected and
\[
\sup_{B_r(x) \cap \Omega} \Big[|\nabla Q_L| + |X_L|\Big] \leq C_1.
\]
Then, for any $\mu$ $\in$ $(0,1]$ and $j$ $\geq$ $0$, there holds
\[
|\nabla^{j+2} Q_L|(y) \leq \frac{C(a^2,b^2,c^2,C_1,\Omega,Q_b,r,\mu,j)}{\dist(y,\partial\Omega)^{j + \mu}} \text{ in } B_{r/2}(x) \cap \Omega.
\]
\end{proposition}

We start with a simple lemma, which is a local version of \cite[Theorem 4.9]{GT}.

\begin{lemma}\label{VanNearBdry}
Assume that $u$ $\in$ $C^2(B_r^+) \cap C^0(\overline{B_r^+})$ satisfies $u\big|_{B_r \cap \{x_n = 0\}}$ $\equiv$ $0$ and
\[
|\Delta u|(x) \leq x_n^{-\mu-1} \text{ in } B_r^+,
\]
for some $\mu$ $\in$ $(0,1)$. Then there exists a universal constant $C$ such that
\[
u(x) \leq C\,r^{-1}\,\sup_{B_r^+} |u|\,x_n + \frac{1}{\mu - \mu^2}\,x_n^{1-\mu} \text{ in } B_{r/2}^+.
\]
\end{lemma}

\bproof It is enough to consider $r$ $=$ $1$. Let $w$ be the solution to
\begin{align*}
\Delta w
	&= 0 \text{ in } B_1^+,\\
w
	&= 0 \text{ on } B_1 \cap \{x_n = 0\},\\
w
	&= 1 \text{ on } \partial B_1 \cap \{x_n > 0\}.
\end{align*}
By standard elliptic estimates, we have
\[
w(x) \leq C\,x_n \text{ in } B_{1/2}^+.
\]
Since the function $v$ $=$ $\sup_{B_r^+}|u|\,w + \frac{1}{\mu - \mu^2}x_n^{1-\mu}$ satisfies
\begin{align*}
\Delta v
	&= - x_n^{-\mu - 1} \text{ in } B_1^+,
v \geq |u| \text{ on } \partial B_1^+,
\end{align*}
the assertion follows from the maximum principle.
\eproof

\begin{lemma}\label{BdryC2Est}
Under the assumption of Proposition \ref{BdryCjEst}, we have for any $\mu$ $\in$ $(0,1]$ that
\begin{equation}
|\nabla^2 Q_{L}|(y)
	\leq \frac{C(a^2, b^2, c^2, C_1,\Omega,Q_b,\mu)}{r^{1-\mu}\,\dist(y,\partial\Omega)^\mu} \text{ for } y \in B_{r/4}(x) \cap \Omega.\label{BC2E::Est01}
\end{equation}
\end{lemma}

\bproof Like usual, we will write $Q$ for $Q_L$ and $B_r$ for $B_r(x)$. We can assume that $r$ $<$ $1$. Also, the conclusion for $\mu$ $=$ $1$ is true in view of Corollary \ref{PUHDE::Cor}. We will now fix some $\mu$ $\in$ $(0,1)$.

For simplicity in presenting, we will only consider the case where $x$ $=$ $0$ and $B_r \cap \Omega$ $=$ $B_r^+$. The general case requires minor changes due to the procedure of flattening the boundary.

By \eqref{PUHDE::Eqn01} and our hypothesis, $Q$ satisfies an equation of the form
\[
\Delta Q
	= f \text{ in } B_r \cap \Omega,
\]
where $f$ satisfies $|f|(y)$ $\leq$ $C$ in $B_r \cap \Omega$. Moreover, by Corollary \ref{PUHDE::Cor},
\begin{equation}
|\nabla f|(y) \leq \frac{C}{\dist(y,\partial\Omega)} \text{ in } B_{3r/4} \cap \Omega.
\label{BC2E::Eqn00}
\end{equation}

Now, pick $\alpha$ $\in$ $\{1,2\}$ and let $U$ $=$ $\nabla_\alpha Q$. Let $W$ be the solution to
\begin{align*}
\Delta W
	&= 0 \text{ in } B_{3r/4}^+,\\
W
	&= U \text{ on } \partial B_{3r/4}^+.
\end{align*}
Define $W_0$ by $W(x_1, x_2, x_3)$ $=$ $U(x_1, x_2, 0)$. Applying Lemma \ref{SimpleEllEst} to $W - W_0$ we get
\[
\sup_{B_{r/2}^+}|\nabla W|
	\leq C\Big[\sup_{B_{3r/4} \cap \{x_3 = 0\}} |\nabla_T U| + r\,\sup_{B_{3r/4} \cap \{x_3 = 0\}} |\nabla_T^2 U| + \frac{1}{r^{3/2}}\|\nabla(W - W_0)\|_{L^2(B_{2r/3}^+)}\Big].
\]
On the other hand, as $W - W_0$ vanishes on $B_{3r/4}^+ \cap \{x_3 = 0\}$, a simple integration by parts using a cut-off function supported in $B_{3r/4}^+$ which is identically $1$ in $B_{2r/3}^+$ gives
\[
\frac{1}{r^{3/2}}\|\nabla(W - W_0)\|_{L^2(B_{2r/3}^+)} \leq C\Big[r\,\sup_{B_{3r/4}^+} |\Delta(W - W_0)| + \frac{1}{r}\sup_{B_{3r/4}^+} |W - W_0|\Big].
\]
Taking the last two inequalities into account and using the maximum principle we get
\begin{equation}
\sup_{B_{r/2}^+}|\nabla W|
	\leq C\Big[\sup_{B_{3r/4} \cap \{x_3 = 0\}} |\nabla_T U| + r\,\sup_{B_{3r/4} \cap \{x_3 = 0\}} |\nabla_T^2 U| + \frac{1}{r}\sup_{\partial B_{3r/4}^+} |U|\Big] \leq \frac{C}{r}.
\label{BC2E::Eqn01}
\end{equation}
Let $V$ $=$ $U - W$. Then $|V(y)|$ $\leq$ $C\sup_{B_{3r/4}^+}|U|$ $\leq$ $C$ and, by \eqref{BC2E::Eqn00},
\begin{align*}
|\Delta V(y)|
	&\leq C\,r^{\mu}\,y_3^{-\mu - 1} \text{ in } B_{3r/4}^+,\\
V
	&= 0 \text{ on } \partial B_{3r/4}^+.
\end{align*}
Hence, by Lemma \ref{VanNearBdry},
\[
|V|(y) \leq \frac{C}{r}\sup_{B_r^+}|V|\,y_3 + C\,r^{\mu}\,y_3^{1-\mu} \leq C\,r^{\mu-1}\,y_3^{1-\mu} \text{ in } B_{r/2}^+.
\]
It thus follows by interpolating (see e.g. \cite[Lemma A.1]{BBH}) that
\begin{equation}
|\nabla V(y)|^2 \leq C\sup_{B_{y_3/2}(y)}|V|\Big[\sup_{B_{y_3/2}(y)}|\Delta V| + \frac{1}{y_3^2}\sup_{B_{y_3/2}(y)}|V|\Big] \leq C\,r^{2\mu-2}\,y_3^{-2\mu} \text{ in } B_{r/4}^+.
\label{BC2E::Eqn02}
\end{equation}
Summing up \eqref{BC2E::Eqn01} and \eqref{BC2E::Eqn02} we get
\[
|\nabla U|(y) \leq C\,r^{\mu-1}\,y_3^{-\mu} \text{ in } B_{r/4}^+.
\]
As $U$ $=$ $\nabla_\alpha Q$ with $\alpha$ $\in$ $\{1,2\}$, we have obtained the required bound for $\nabla_{ij} Q$ with $(i,j)$ $\neq$ $(3,3)$. The remaining estimate for $\nabla_{33} Q$ follows from $|\Delta Q|$ $=$ $|f|$ $\leq$ $C$ in $B_r \cap \Omega$. We conclude the proof.
\eproof

\medskip
\noindent{\bf Proof of Proposition \ref{BdryCjEst}.} We will do an induction on $j$. The assertion for $j$ $=$ $0$ follows from Lemma \ref{BdryC2Est}. By differentiating \eqref{PUHDE::Eqn01} to $(j+2)$-nd order, we see that
\[
\Delta(\nabla^{j+2} Q) = f_j
\]
where, by Corollary \ref{PUHDE::Cor},
\[
|f_j|(y) \leq \frac{C}{\dist(y,\partial\Omega)^{j+2}} \text{ in } B_{r/2}(x) \cap \Omega.
\]
Therefore, by interpolation (see e.g. \cite[Lemma A.1]{BBH}), we get the required estimate for $\nabla^{j+3}%
Q$.
\eproof

\begin{proposition}\label{EqnConv-Part2}
Assume for some $B_r(x)$ $\Subset$ $\Omega$, $L$ $\in$ $(0,1)$ and $C_1$ $>$ $0$ that 
\[
\sup_{B_r(x)} \big[|\nabla Q_L| + |X_L|\big] \leq C_1.
\]
Then, there exist constants $C$ and $C'$ depending on $(a^2,b^2,c^2,C_1,j)$ such that for $0$ $<$ $s$ $<$ $r$, there holds
\begin{align}
\sup_{B_s(x)}\Big(|\nabla^j Y_L| + |\nabla^j Z_L|\Big) \leq \frac{C}{(r-s)^j}\Big[(r-s)^{-1}\,L + \exp\Big(-\frac{C'}{\sqrt{L}}(r-s)\Big)\Big],
\label{ECP2P::Est02}
\end{align}
where $Y_L$ and $Z_L$ are defined by \eqref{YDef} and \eqref{ZDef}, respectively.
\end{proposition}

\bproof We use some idea from the proof of Lemma \ref{PrepUniHighDerEst}. Like usual, we write $Q$ $=$ $Q_L$, $X$ $=$ $X_L$, $Y$ $=$ $Y_L$ and $Z$ $=$ $Z_L$.

Fix some $\mu$ $\in$ $(0,1)$ for the moment. By Corollary \ref{PUHDE::Cor} and Proposition \ref{BdryCjEst}, we have the following estimates,
\begin{align}
\sup_{B_s(x)}|\nabla^{j+2}Q| \leq \frac{C}{(r-s)^{j+\mu}} \text{ and } 
		\sup_{B_s(x)}\big[|\nabla^{j}X|(y) + |\nabla^j Y|(y) + |\nabla^j Z|\big] \leq \frac{C}{(r-s)^{j}}.
\label{ECP2P::Eqn00}
\end{align}

We split the proof into three steps, which bear some analogy to Steps 3-5 in the proof of Lemma \ref{PrepUniHighDerEst}. In the first step, we prove an estimate for $\tr(Z)$. The remaining steps are to prove \eqref{ECP2P::Est02}. In the proof, $R(A_1, A_2, ..., A_t; A_{t+1}, \ldots A_s)$ will be used to denote various explicitly computable polynomials in the $A_l$'s which are (jointly) linear in $(A_{t+1}, \ldots, A_s)$.

\medskip
\noindent\underline{Step 1:} We show that
\begin{equation}
\sup_{B_s(x)}|\nabla^j Y| \leq \frac{C}{(r-s)^j}\Big[\frac{L}{r-s} + \exp\Big(-\frac{C'}{\sqrt{L}}(r-s)\Big)\Big].
\label{ECP2P1::Eqn01}
\end{equation}
By \eqref{PUHDE::Eqn04Bis}, we have
\[
L\Delta\tr(X) = (2a^2 + \frac{1}{3}b^2\,s_+)Y + L\,R(X).
\]
Hence, by noting that $Y$ $=$ $\tr(X) + \frac{6}{6a^2 + b^2s_+}|\nabla Q|^2$ and using \eqref{PUHDE::Eqn01}, 
\begin{align*}
L\Delta Y
	&= (2a^2 + \frac{1}{3}b^2\,s_+)Y +  L\,R(Q,\nabla Q,X;|\nabla^2 Q|^2,\nabla X).
\end{align*}
Note that by choosing $\mu$ sufficiently small in \eqref{ECP2P::Eqn00}, we have
\[
\sup_{B_{(r+s)/2}(y)}|\nabla^{j} R(Q,\nabla Q,X;|\nabla^2 Q|^2,\nabla X)| \leq \frac{C}{(r-s)^{j+1}}.
\]
Hence, by Lemma \ref{BBHMaxPrin} and \eqref{ECP2P::Eqn00},
\[
\sup_{B_{s}(x)}|\nabla^{j}Y|
	\leq \frac{C\,L}{(r-s)^{j+1}} + \frac{C}{(r-s)^j}\,\exp\Big(-\frac{C_1}{\sqrt{L}}(r-s)\Big),
\]
which implies \eqref{ECP2P1::Eqn01} and so complete Step 1.

\medskip
\noindent\underline{Step 2:} We will show that
\begin{equation}
\sup_{B_s(x)}\Big|\nabla^j\Big(\big(Q - \frac{2}{3}s_+\,\IdMa\big)Z\Big)\Big| \leq \frac{C}{(r-s)^j}\Big[\frac{L}{r-s} + \exp\Big(-\frac{C'}{\sqrt{L}}(r-s)\Big)\Big].
\label{ECP2P2::Est02-1}
\end{equation}

Multiplying \eqref{EC-Ez::Eqn03} to the left by $Q - \frac{1}{6}s_+\,\IdMa$, we get
\begin{equation}
L\big(\frac{1}{4}s_+^2\,\IdMa + \frac{1}{2}\,L\,X\big)\Delta X
	= -\frac{1}{2}b^2\,s_+^2\,\big(Q - \frac{1}{6}s_+\,\IdMa\big)Z + L\,\big(Q - \frac{1}{6}s_+\,\IdMa\big)\,R(Q,X;\frac{1}{L}Y).
\label{ECP2P2::Est02-2}
\end{equation}
Recall the definition of $X$ and \eqref{PUHDE::Eqn01}, the above equation implies
\begin{align}
L\,\Delta X
	&=\frac{4}{s_+^2}\Big[\text{Right hand side of \eqref{ECP2P2::Est02-2}} - \frac{1}{2}L\,X\,\Delta(L\,X)\Big]\nonumber\\
	&= - 2b^2\,\big(Q - \frac{1}{6}s_+\,\IdMa\big)Z + L\,R(Q, \nabla Q, X; \frac{1}{L}Y).
\label{ECP2P2::Eqn03}
\end{align}
It then follows from \eqref{EC-Ez::Eqn03}, \eqref{ECP2P2::Eqn03} and \eqref{PUHDE::Eqn01} that
\[
L\big(Q - \frac{2}{3}s_+\,\IdMa\big)\Delta X
	= b^2\,s_+\big(Q - \frac{2}{3}s_+\,\IdMa\big)Z + L\,R(Q, \nabla Q, X; \frac{1}{L}Y),
\]
and
\[
L\Delta\Big[\big(Q - \frac{2}{3}s_+\,\IdMa\big)X\Big]
	= b^2\,s_+\big(Q - \frac{2}{3}s_+\,\IdMa\big)Z + L\,R(Q, \nabla Q, X, \nabla X; \frac{1}{L}Y).
\]
Recalling the definition of $Z$ and \eqref{PUHDE::Eqn01}, we then arrive at
\begin{multline}
L\Delta\Big[\big(Q - \frac{2}{3}s_+\,\IdMa\big)Z\Big]
	= b^2\,s_+\big(Q - \frac{2}{3}s_+\,\IdMa\big)Z\\
		+ L\,R(Q, \nabla Q, X; \nabla^2 Q, (\nabla^2 Q)^2, \nabla X, \frac{1}{L}Y).
\label{ECP2P2::Eqn04}
\end{multline}
Applying Lemma \ref{BBHMaxPrin} and using \eqref{ECP2P::Eqn00} and\eqref{ECP2P1::Eqn01}, we arrive at \eqref{ECP2P2::Est02-1} (by choosing $\mu$ appropriately small). This completes Step 2.

\medskip
\noindent\underline{Step 3:} We will show that
\begin{equation}
\sup_{B_s(x)}\Big|\nabla^j\Big(\big(Q + \frac{1}{3}s_+\,\IdMa\big)Z\Big)\Big| \leq \frac{C}{(r-s)^j}\Big[\frac{L}{r-s} + \exp\Big(-\frac{C'}{\sqrt{L}}(r-s)\Big)\Big].
\label{ECP2P0::Est02-2}
\end{equation}
It is readily seen that \eqref{ECP2P::Est02} is a consequence of \eqref{ECP2P2::Est02-1} and \eqref{ECP2P0::Est02-2}.

Since $Q$ is a traceless $3\times 3$ matrix, Cayley's theorem gives
\begin{equation}
Q^3 - \frac{1}{2}\tr(Q^2)\,Q - \frac{1}{3}\tr(Q^3)\,\IdMa = 0.
\label{ECP2P0::Eqn01}
\end{equation}

Noting that $\Delta Q$ commutes with $Q$ by \eqref{L-DG::ELEqn}, we compute
\begin{equation}
\Delta Q^3 = 3\,Q^2\,\Delta Q + 2\sum_{\alpha=1}^3\big[\nabla_\alpha Q\,Q\,\nabla_\alpha Q + (\nabla_\alpha Q)^2\,Q + Q\,(\nabla_\alpha Q)^2\big] =: 3\,Q^2\,\Delta Q + 2V.
\label{ECP2P0::Eqn02}
\end{equation}
To simplify $V$, we note that the definition of $X$ implies
\begin{equation}
Q\,\nabla_\alpha Q + \nabla_\alpha Q\,Q = \frac{1}{3}s_+\,\nabla_\alpha Q + L\,\nabla_\alpha X,
\label{ECP2P0::Eqn03}
\end{equation}
and so
\begin{align*}
V
	&= \sum_{\alpha=1}^3\big[\nabla_\alpha Q\,(Q\,\nabla_\alpha Q + \nabla_\alpha Q\,Q) + Q\,(\nabla_\alpha Q)^2\big]\\
	&= \big(Q +\frac{1}{3}s_+\,\IdMa\big)(\nabla Q)^2 + L\,R(\nabla Q;\nabla X) =: U + L\,R(\nabla Q;\nabla X).
\end{align*}
Substituting this into \eqref{ECP2P0::Eqn02} and using \eqref{PUHDE::Eqn01}, we get
\begin{align}
\Delta Q^3
	&= 3\,Q^2\,\Delta Q + 2U + L\,R(\nabla Q;\nabla X)\nonumber\\
	&= s_+\big(Q + \frac{2}{3}s_+\,\IdMa\big)\,\Delta Q + 2U + L\,R(Q,\nabla Q, X;\nabla X)\nonumber\\
	&= s_+\big(\frac{1}{3}b^2 + c^2\,s_+\big)\tr(X)\,Q + \frac{2}{9}s_+^2(b^2 + c^2\,s_+)\tr(X)\,\IdMa - b^2\,s_+\,Q\,X - \frac{2}{3}b^2\,s_+^2\,X\nonumber\\
		&\qquad\qquad+ 2U + L\,R(Q,\nabla Q, X;\nabla X).
\label{ECP2P0::Eqn04}
\end{align}
Taking trace of \eqref{ECP2P0::Eqn04}, we get
\begin{multline}
\Delta \tr(Q^3) = \frac{2}{3}s_+^2(b^2 + c^2\,s_+)\tr(X) - b^2\,s_+\,\tr(Q\,X) - \frac{2}{3}b^2\,s_+^2\,\tr(X)\\
		\qquad\qquad+ 2\,\tr(U) + L\,R(Q,\nabla Q, X;\nabla X).
\label{ECP2P0::Eqn05}
\end{multline}

Similarly, we have
\[
\Delta Q^2 = 2\,Q\,\Delta Q + 2\sum_{\alpha = 1}^3 (\nabla_\alpha Q)^2,
\]
and so, in view of \eqref{ECP2P0::Eqn03} and \eqref{PUHDE::Eqn01},
\begin{align}
\Delta\big[\tr(Q^2)Q\big]
	&= 2\,\tr(Q\,\Delta Q)\,Q + 2|\nabla Q|^2\,Q + \tr(Q^2)\,\Delta Q + 2\sum_{\alpha = 1}^3 \nabla_\alpha \tr(Q^2)\,\nabla_\alpha Q\nonumber\\
	&= 2\,\tr(Q\,\Delta Q)\,Q + 2|\nabla Q|^2\,Q + \frac{2}{3}s_+^2\,\Delta Q + L\,R(Q,\nabla Q, X;\nabla X)\nonumber\\
	&= 2\,\big[c^2\,\tr(X)\,\tr(Q^2) - b^2\tr(Q\,X)\big]\,Q + 2|\nabla Q|^2\,Q\nonumber\\
		&\qquad\qquad + \frac{2}{3}s_+^2\big[c^2\,\tr(X)\,Q + \frac{1}{3}b^2\,\tr(X)\,\IdMa - b^2\,X\big] + L\,R(Q,\nabla Q, X;\nabla X)\nonumber\\
	&= 2c^2\,s_+^2\,\tr(X)\,Q - 2b^2\tr(Q\,X)\,Q + \frac{2}{9}b^2\,s_+^2\,\tr(X)\,\IdMa - \frac{2}{3}b^2\,s_+^2\,X\nonumber\\
		&\qquad\qquad  + 2|\nabla Q|^2\,Q + L\,R(Q,\nabla Q, X;\nabla X).
\label{ECP2P0::Eqn06}
\end{align}

Combining \eqref{ECP2P0::Eqn01}, \eqref{ECP2P0::Eqn04}, \eqref{ECP2P0::Eqn05} and \eqref{ECP2P0::Eqn06} together, we get
\begin{align*}
0
	&= \Big[s_+\big(\frac{1}{3}b^2 + c^2\,s_+\big)\tr(X)\,Q + \frac{2}{9}s_+^2(b^2 + c^2\,s_+)\tr(X)\,\IdMa - b^2\,s_+\,Q\,X - \frac{2}{3}b^2\,s_+^2\,X + 2U\Big]\\
		&\qquad\qquad - \frac{1}{2}\Big[2c^2\,s_+^2\,\tr(X)\,Q - 2b^2\tr(Q\,X)\,Q + \frac{2}{9}b^2\,s_+^2\,\tr(X)\,\IdMa - \frac{2}{3}b^2\,s_+^2\,X + 2|\nabla Q|^2\,Q\Big]\\
		&\qquad\qquad - \frac{1}{3}\Big[\frac{2}{3}s_+^2(b^2 + c^2\,s_+)\tr(X) - b^2\,s_+\,\tr(Q\,X) - \frac{2}{3}b^2\,s_+^2\,\tr(X) + 2\,\tr(U)\Big]\IdMa\\
		&\qquad\qquad + L\,R(Q,\nabla Q, X;\nabla X),
\end{align*}
which simplifies to
\begin{align}
0
	&= \big(Q + \frac{1}{3}s_+\,\IdMa\big)\Big[-b^2\,s_+\,X + \frac{1}{3}b^2s_+\,\tr(X)\,\IdMa + b^2\,\tr(Q\,X)\,\IdMa\Big]+ 2U\nonumber\\
		&\qquad\qquad - |\nabla Q|^2\,Q - \frac{2}{3}\,\tr(U)\,\IdMa + L\,R(Q,\nabla Q, X;\nabla X).
\label{ECP2P0::Eqn07}
\end{align}
Multiplying \eqref{ECP2P0::Eqn07} to the left by $Q + \frac{1}{3}s_+\,\IdMa$ and noting that $(Q + \frac{1}{3}s_+\,\IdMa)^2$ $=$ $s_+\big(Q + \frac{1}{3}s_+\,\IdMa\big) + L\,X$, we get
\begin{align}
0
	&= \big(Q + \frac{1}{3}s_+\,\IdMa\big)\Big[-b^2\,s_+\,X + \frac{1}{3}b^2s_+\,\tr(X)\,\IdMa + b^2\,\tr(Q\,X)\,\IdMa + 2\sum_{\alpha=1}^3 (\nabla_\alpha Q)^2\nonumber\\
		&\qquad\qquad - \frac{2}{3}|\nabla Q|^2\,\IdMa - \frac{2}{3s_+}\,\tr(U)\,\IdMa\Big] + L\,R(Q,\nabla Q, X;\nabla X).
\label{ECP2P0::Eqn07bis}
\end{align}
Comparing \eqref{ECP2P0::Eqn07} and \eqref{ECP2P0::Eqn07bis} we get
\[
|\nabla Q|^2\,Q + \frac{2}{3}\,\tr(U)\,\IdMa = \Big[\frac{2}{3}|\nabla Q|^2 + \frac{2}{3s_+}\,\tr(U)\Big]\big(Q + \frac{1}{3}s_+\,\IdMa\big) + L\,R(Q,\nabla Q, X;\nabla X).
\]
Taking trace of the above yields
\begin{equation}
\tr(U) = \frac{1}{2}s_+|\nabla Q|^2 + L\,R(Q,\nabla Q, X;\nabla X).
\label{ECP2P0::Eqn08}
\end{equation}
Subtituting \eqref{ECP2P0::Eqn08} into \eqref{ECP2P0::Eqn07bis} and using \eqref{PUHDE::Eqn04} and $\tr(X)$ $=$ $Y  - \frac{6}{6a^2 + b^2\,s_+}\,|\nabla Q|^2$, we then obtain
\begin{align}
0
	&= \big(Q + \frac{1}{3}s_+\,\IdMa\big)\Big[-b^2\,s_+\,X - \frac{2b^2\,s_+ + 4c^2\,s_+^2}{6a^2 + b^2\,s_+}|\nabla Q|^2\,\IdMa + 2\sum_{\alpha=1}^3 (\nabla_\alpha Q)^2\Big]\nonumber\\
		&\qquad\qquad + L\,R(Q,\nabla Q, X;\nabla X,\frac{1}{L}Y)\nonumber\\
	&= -b^2\,s_+\big(Q + \frac{1}{3}s_+\,\IdMa\big)Z + L\,R(Q,\nabla Q, X;\nabla X,\frac{1}{L}Y).
\end{align}
Estimate \eqref{ECP2P0::Est02-2} follows, whence \eqref{ECP2P::Est02}.
\eproof

As a consequence of Proposition \ref{EqnConv-Part2}, we have

\begin{corollary}\label{ECP2P::Cor}
Assume 
\[
\sup_{B_r(x) \cap \Omega} \big[|\nabla Q_L| + |X_L|\big] \leq C_1.
\]
for some $x$ $\in$ $\bar\Omega$, $r$ $>$ $0$, $L$ $\in$ $(0,1)$ and $C_1$ $>$ $0$. Then, for $j$ $\geq$ $0$, there exist constants $C$ and $C'$ depending on $(a^2,b^2,c^2,C_1,j)$ such that, in $B_{r/2}(x) \cap \Omega$, there hold
\[
|\nabla^j Y_L|(y) + |\nabla^j Z_L|(y)
	\leq \frac{C}{\tilde d_{r/2}(y)^{j}}\Big[\frac{L}{\tilde d_{r/2}(y)} + \exp\Big(-\frac{C'}{\sqrt{L}}\,\tilde d_{r/2}(y)\Big)\Big].
\]
where $\tilde d_s(y)$ $=$ $\dist(y,\partial(B_s(x) \cap \Omega))$.
\end{corollary}

\begin{corollary}\label{ImprRemEst}
Under the setting of Theorem \ref{MainThm1}, the $R_{L_{k'}}$ satisfies for $y$ $\in$ $K$,
\[
|\nabla^j R_{L_{k'}}|(y) \leq \frac{C}{d(y)^j}\Big[\frac{L}{d(y)} + \exp\Big(-\frac{C'}{\sqrt{L}}\,d(y)\Big)\Big],
\]
where $d(y)$ $=$ $\dist(y,\partial\Omega)$.
\end{corollary}

\bproof The assertion follows from \eqref{ExplicitApproxEqn} and Corollary \ref{ECP2P::Cor}.
\eproof

\section{Rate of convergence for a class of minimizers}\label{SEC-Rate}

In this section, we will show that for a class of limit map $Q_*$, the difference $Q_L - Q_*$ is of order $L$. It is conceivable that such statement may not hold when $Q_*$ is not an isolated minimizer of the limit functional $I_*$. Recall that the Euler-Lagrange equation for $I_*$ is
\[
\Delta Q_* = -\frac{4}{s_+^2}(Q_* - \frac{1}{6}\,s_+\,\IdMa)\sum_{\alpha=1}^3 (\nabla_\alpha Q_*)^2.
\]
It is thus reasonable to assume some kind of bijectivity for the ``linearized operator''
\begin{equation}
\CL_{Q_*}\Psi := \Delta\Psi + \frac{4}{s_+^2}(Q_* - \frac{1}{6}\,s_+\,\IdMa)\sum_{\alpha=1}^3 (\nabla_\alpha Q_*\,\nabla_\alpha \Psi + \nabla_\alpha \Psi\,\nabla_\alpha Q_*) + \frac{4}{s_+^2} \Psi\sum_{\alpha=1}^3 (\nabla_\alpha Q_*)^2.
\label{LinearizedOp}
\end{equation}

In this paper, we will restrict our attention to the case where $Q_*$ is smooth in the entire domain $\Omega$. The general case will be considered in a forthcoming paper. We establish:

\begin{proposition}\label{Rate}
Let $\Omega$ be an open bounded subset of $\RR^3$. Assume that $Q_{L_k}$ converges strongly in $H^1(\Omega, \mcS_0)$ to a minimizer $Q_*$ $\in$ $H^1(\Omega,\mcS_*)$ of $I_*$ for some sequence $L_k$ $\rightarrow$ $0$. Furthermore, assume that $Q_*$ is smooth and the linearized harmonic map operator $\CL_{Q_*}:$ $H^1_0(\Omega,M_{3\times 3})$ $\rightarrow$ $H^{-1}(\Omega,M_{3\times 3})$ defined by \eqref{LinearizedOp} is bijective. Then there exists $\bar L > 0$ such that, for $L_k < \bar L$,
\begin{equation}
\|Q_{L_{k}}- Q_*\|_{H^s(\Omega)}\le C(a^2, b^2, c^2, \Omega, Q_b, \mu, s)\,L_{k} \text{ for } 0 < s < \frac{1}{2},
\label{est:Hsrate}
\end{equation}
and
\begin{equation}
\|\nabla^j(Q_{L_{k}}- Q_*)\|_{L^\infty(K)}\le C(a^2, b^2, c^2, \Omega, K, Q_b, \mu, j)L_{k} \text{ for } K \Subset \Omega, j \geq 0.
\label{est:Higherrate}
\end{equation}
\end{proposition}

We start with a few lemmas.
\begin{lemma}\label{Trick}
Let $\Omega$ $\subset$ $\RR^3$ be a bounded domain with smooth boundary . Consider the elliptic operator
\begin{equation}
(\CL u)_\alpha := \Delta u_\alpha + b_\alpha^\beta \cdot \nabla u_\beta + c_\alpha^\beta\,u_\beta,
\label{LES}
\end{equation}
where $b_\alpha^\beta$ $\in$ $W^{1,p}(\bar\Omega)$ for some $p \geq 2$, $c_\alpha^\beta$ $\in$ $L^\infty(\Omega)$. If $\CL:$ $H^1_0(\Omega)$ $\rightarrow$ $H^{-1}(\Omega)$ is bijective, then
\begin{equation}
\|u\|_{L^2(\Omega)} \leq C\,\|\CL u\|_{\tilde H^{-2}(\Omega)} \text{ for any } u \in H^1_0(\Omega),
\label{HsEst::Est01}
\end{equation}
where $\tilde H^{-2}(\Omega)$ denotes the dual space of $H^2(\Omega) \cap H^1_0(\Omega)$, and $C$ is a constant that depends only on $\Omega$ and the spectral of $\CL$.
\end{lemma}

\bproof Let $\CL^*:$ $H^1_0(\Omega)$ $\rightarrow$ $H^{-1}(\Omega)$ denote the adjoint of $\CL$, i.e.
\[
(\CL^* v)^\alpha = \Delta v^\alpha + \div(b_\gamma^\alpha\,v^\gamma) + c^\alpha_\gamma\,v^\gamma.
\]
Since $\CL$ is bijective, so is $\CL^*$. Let $v$ $\in$ $H^1_0(\Omega)$ be the unique solution to
\[
\CL^* v = u.
\]
By regularity theory for elliptic operators (see e.g. \cite[Chapter 3, Theorem 10.1]{L-U}), we have $v$ $\in$ $H^2(\Omega)$ and
\[
\|v\|_{H^2(\Omega)} \leq C(\|u\|_{L^2(\Omega)} + \|v\|_{L^2(\Omega)}).
\]
Moreover, because the kernel of $\CL^*:$ $H^1_0(\Omega)$ $\rightarrow$ $H^{-1}(\Omega)$ is trivial, the $\|v\|_{L^2(\Omega)}$-term on the right hand side can be dropped. Hence
\[
\|v\|_{H^2(\Omega)} \leq C\,\|u\|_{L^2(\Omega)}.
\]
It thus follows from integrating by parts that
\begin{align*}
\|u\|_{L^2(\Omega)}^2
	&= \int_{\Omega} u\,\CL^* v\,dx = \left<\CL u, v\right>_{(H^{-1}(\Omega),H^1_0(\Omega))}\\
	&\leq \|\CL u\|_{\tilde H^{-2}(\Omega)}\,\|v\|_{H^2(\Omega) \cap H^1_0(\Omega)}\\
	&\leq C\|\CL u\|_{\tilde H^{-2}(\Omega)}\,\|u\|_{L^2(\Omega)}.
\end{align*}
The assertion follows.
\eproof

\begin{lemma}\label{NormEst}
Let $\Omega\subset\mathbb{R}^3$ be a smooth bounded domain. We denote by $\tilde H^{-2}(\Omega)$ the dual space of $H^2(\Omega)\cap H^1_0(\Omega)$. If $|R^\sharp(x)| \le Ce^{-\frac{Cd(x)}{\sqrt{L}}}$ where $d(x)$ denotes the distance to the boundary $\partial \Omega$, then
\begin{equation}
 \|R^\sharp \|_{\tilde H^{-2}(\Omega)}\le CL
 \label{Rintilde}
 \end{equation}
with the constant $C$ independent of $L$.
\end{lemma}

\bproof \underline{Step 1:} We consider first the case when $\Omega=[0,1]^3$ and ``$d(x)=x_3$''. We show that  if $|R^\sharp (x)|\le Ce^{-\frac{Cx_3}{\sqrt{L}}}$ then (\ref{Rintilde}) holds.

We take a test function $\varphi\in H^2(\Omega) \cap H^1_0(\Omega)$. We compute
\begin{align}
 |\langle R^\sharp, \varphi\rangle|
	&\le \underbrace{\Big|\int_{[0,1]^3} \Big(\int_{x_3}^1\big(\int_z^1 R^\sharp(x',w)\,dw\big)\,dz\Big)\partial^2_{x_3}\varphi(x',x_3)\,dx_3\,dx'\Big|}_{\mathcal{I}}\nonumber\\
		&\qquad\qquad +\underbrace{\Big|\int_{[0,1]^2}\Big(\int_0^1\big(\int_z^1 R^\sharp(x',w)\,dw\big)\,dz\Big)\partial_{x_3}\varphi(x',0)\,dx'\Big|}_{\mathcal{II}}.
\label{Rest}
\end{align}
$\mathcal{I}$ can be estimated as follows.
\begin{align}
|\mathcal{I}|
	&\le \bigg\{\int_{[0,1]^3}\Big(\int_{x_3}^1\big(\int_z^1 |R^\sharp(x',w)|\,dw\big)\,dz\Big)^2\,dx\bigg\}^{\frac{1}{2}}\bigg\{\int_{[0,1]^3}|\Delta\varphi|^2\,dx\bigg\}^{\frac{1}{2}}\nonumber\\
	&\le C\bigg\{\int_{[0,1]^3}\Big(\int_{x_3}^1\big(\int_z^1 e^{-\frac{Cw}{\sqrt{L}}}\,dw\big)\,dz\Big)^2\,dx\bigg\}^{\frac{1}{2}}\bigg\{\int_{[0,1]^3}|\Delta\varphi|^2\,dx\bigg\}^{\frac{1}{2}}\nonumber\\
	&\le C\sqrt{L}\bigg\{\int_{[0,1]^3}\big(\int_{x_3}^1 e^{-\frac{Cz}{\sqrt{L}}}\,dz\big)^2\,dx\bigg\}^{\frac{1}{2}}\bigg\{\int_{[0,1]^3}|\Delta\varphi|^2\,dx\bigg\}^{\frac{1}{2}}\nonumber\\
	&\le CL \|\varphi\|_{H^2([0,1]^3)},
  \label{-221}
\end{align}
Similarly,
\begin{align}
 |\mathcal{II}|
	&\le \bigg\{\int_{[0,1]^2}\Big(\int_0^1\big(\int_z^1 R^\sharp(x',w)dw\big)dz\Big)^2\,dx'\bigg\}^{\frac{1}{2}}\bigg\{\int_{[0,1]^2} |\partial_{x_3} \varphi(x',0)|^2\,dx'\bigg\}^{\frac{1}{2}}\nonumber\\
	&\le C\sqrt{L}\bigg\{\int_{[0,1]^2}\big(\int_0^1 e^{-\frac{z}{\sqrt{L}}}\,dz\big)^2\,dx'\bigg\}^{\frac{1}{2}}\|\textrm{Tr}(\nabla\varphi)\|_{L^2(\partial\Omega)}\nonumber\\
	&\le CL\|\varphi\|_{H^2(\Omega)}.
 \label{-222}
\end{align}
Using (\ref{-221}) and (\ref{-222}) in (\ref{Rest})  we obtain the claimed estimate (\ref{Rintilde}) 
 with the constant $C$ independent of $L$.

\medskip
\noindent\underline{Step 2:} We now consider the general case. Select $\delta$ $>$ $0$ sufficiently small so that for any $x$ $\in$ $\Omega$ with $d(x)$ $<$ $\delta$, there is a unique $\pi(x)$ $\in$ $\partial\Omega$ such that $d(x)$ $=$ $|x - \pi(x)|$. Cover $\partial\Omega$ by finitely many open set $O_i$ such that each $O_i$ is diffeomorphic to $[0,1]^2$ by some map $\psi_i:$ $[0,1]^2$ $\rightarrow$ $O_i$. Define $\phi_i:$ $[0,1]^2 \times [0,\delta]$ $\rightarrow$ $\Omega$ by mapping $(\xi,\eta)$ $\in$ $[0,1]^2 \times [0,\delta]$ to the unique point $x$ in $\Omega$ such that $\pi(x)$ $=$ $\psi_i(\xi)$ and $d(x)$ $=$ $\eta$. Let $U_i$ $=$ $\phi_i([0,1]^2 \times [0,\delta])$ and $U_*$ $=$ $\{x\in\Omega: d(x) > \delta/2\}$. Then $\{U_*\} \cup \{U_i\}$ forms a covering of $\Omega$. Let $\mathbf{1}$ $=$ $\theta_* + \sum \theta_i$ be a partition of unity of $\Omega$ associated to this covering. 

Fix $\varphi$ $\in$ $H^2(\Omega) \cap H^1_0(\Omega)$.

By our hypothesis, in $U_*$, $|R^\sharp|$ $\leq$ $C\,L$ where $C$ depends on $\delta$. Thus,
\begin{equation}
\Big|\int_{\Omega} R^\sharp\,\varphi\,\theta_*\,dx\Big| \leq C\,L\,\|\varphi\|_{L^2(\Omega)}.
\label{NE::Est1}
\end{equation}
Next, by performing a change of variables, we have
\[
\int_{\Omega} R^\sharp\,\varphi\,\theta_i\,dx = \int_{[0,1]^2 \times [0,\delta]} \tilde R_i\,\varphi \circ \phi_i\,dx,
\]
where $\tilde R_i$ satisfies $|\tilde R_i(x)|$ $\leq$ $C\,e^{-\frac{Cx_3}{\sqrt{L}}}$. Thus, by Step 1,
\begin{equation}
\Big|\int_{\Omega} R^\sharp\,\varphi\,\theta_i\,dx\Big| \leq C\,L\,\|\varphi\|_{H^2(\Omega)}.
\label{NE::Est2}
\end{equation}
Summing up \eqref{NE::Est1} and \eqref{NE::Est2}, we get \eqref{Rintilde}.
\eproof

It turns out that we will also need an estimate similar to that in Lemma \ref{NormEst} but with the norm measure in fractional Sobolev spaces $H^s(\Omega)$. For the clarity in exposition, we state the estimate here while deter its proof until after the proof of Proposition \ref{Rate}.

\begin{lemma}\label{Remainder::HsNorm}
Let $\Omega\subset\mathbb{R}^3$ be a smooth bounded domain. If $|R^\sharp(x)|\le Ce^{-\frac{Cd(x)}{\sqrt{L}}}$ where $d(x)$ denotes the distance to the boundary then
\begin{equation}
 \|R^\sharp\|_{H^s(\Omega)}\le CL
 \label{RinX}
 \end{equation}
for all $s\in [-2,-\frac{3}{2})$.
\end{lemma}

\begin{lemma}\label{NormEstEz}
Let $\Omega\subset\mathbb{R}^3$ be a smooth bounded domain. If $|R^\flat(x)| \le \frac{CL}{d(x)}$ where $d(x)$ denotes the distance to the boundary $\partial \Omega$, then
\begin{equation}
 \|R^\flat\|_{H^{-1}(\Omega)}\le CL
 \label{NEE::Est1}
 \end{equation}
with the constant $C$ independent of $L$.
\end{lemma}

\bproof The assertion is a consequence of the well-known Hardy inequality. We give a proof here for completeness. Similar to the proof of Lemma \ref{NormEst}, it suffices to show that if $\Omega$ $=$ $[0,1]^3$ and $|R(x)|$ $\leq$ $\frac{CL}{x_3}$, then \eqref{NEE::Est1} holds. Moreover, we can assume that $R$ is non-negative. For $\varphi$ $\in$ ${\rm Lip}_0(\Omega)$, we estimate
\begin{align*}
\Big|\int_{[0,1]^3} R\,\varphi\,dx\Big|
	&\leq \Big|\int_{[0,1]^3} \int_{x_3}^1 R(x',z)\,dz\,\partial_{x_3} \varphi(x',x_3)\,dx_3\,dx'\Big|\\
	&\leq C\Big|\int_{[0,1]^3} \ln x_3\,\partial_{x_3} \varphi(x',x_3)\,dx_3\,dx'\Big|
		\leq C\|\varphi\|_{H^1_0(\Omega)}.
\end{align*}
Using this inequality, one can argue using density, Fatou's lemma and the splitting $f$ $=$ $f^+ - f^-$ to show that $R\,\varphi$ is integrable for any $\varphi$ $\in$ $H^1_0(\Omega)$ and the above inequality extends to $\varphi$ $\in$ $H^1_0(\Omega)$. The assertion follows.
\eproof

\medskip
\noindent{\bf Proof of Proposition \ref{Rate}.} Like usual, we drop the subscript $L_k$. By Propositions \ref{C1AlphaConv} and \ref{CjConv} and the smoothness of $Q_*$, we can assume that $Q$ converges in $C^{1,\alpha}(\bar\Omega)$ and $C^j_{\rm loc}(\Omega)$, $j$ $\geq$ $2$, to $Q_*$. Moreover, by Theorem \ref{MainThm1} and Corollary \ref{ImprRemEst}, we have
\begin{equation}
\Delta Q = -\frac{4}{s_+^2}\big(Q - \frac{1}{6}s_+\,\IdMa\big)\sum_{\alpha = 1}^3 (\nabla_\alpha Q)^2 + R,
\label{Rate::ApproxEqn}
\end{equation}
where $R$ $\in$ $C^\infty(\Omega)$ and satisfies the estimate
\begin{equation}
|R(y)| \leq C\,\Big[\frac{L}{d(y)} + \exp\Big(-\frac{d(y)}{\sqrt{L}}\Big)\Big] \text{ and } |\nabla^l R(y)| \leq \frac{C\,L_k}{d(y)^{l+2}} \text{ for } y \in \Omega, l\geq 0.
\label{Rate::ErrEst}
\end{equation}
Here $d(y)$ $=$ $\dist(y,\partial\Omega)$. 

Let $S$ $=$ $Q - Q_*$. By \eqref{Rate::ApproxEqn} and Corollary \ref{HarMapEqn}(iv), 
\begin{multline}
\Delta S + \frac{4}{s_+^2}\left(Q-\frac{1}{6}\,s_+\,\IdMa\right)\sum_{\alpha=1}^3\big[\nabla_\alpha Q\,\nabla_{\alpha} S +\nabla_\alpha S\,\nabla_\alpha Q_*\big]\\
		+\frac{4}{s_+^2} S\,\sum_{\alpha=1}^3(\nabla_\alpha Q_*)^2 = R,
\label{Rate::DeltaS}
\end{multline}
For our purpose, we use the convergence of $Q$ to $Q_*$ to rewrite \eqref{Rate::DeltaS} in the following form:
\begin{equation}
\tilde \CL\,S_{ij} := \CL_{Q_*}\,S_{ij} + \tilde b_{ij}^{\alpha\beta}\,\nabla S_{\alpha\beta} + \tilde c_{ij}^{\alpha\beta}\,S_{\alpha\beta} = R_{ij},
\label{Rate::DSLinearForm}
\end{equation}
where, by Propositions \ref{C1AlphaConv}, \ref{CjConv} and \ref{BdryCjEst}, $\tilde b_{ij}^{\alpha\beta}$ $\in$ $C^{\alpha}(\bar\Omega) \cap C^\infty(\Omega) \cap W^{1,p}(\Omega)$, $\tilde c_{ij}^{\alpha\beta}$ $\in$ $C^{\alpha}(\bar\Omega) \cap C^\infty(\Omega)$ for any $\alpha$ $\in$ $(0,1)$ and $p$ $>$ $1$, and satisfy
\begin{align}
\|\tilde b_{ij}^{\alpha\beta}\|_{W^{1,p}(\Omega)} \leq C \text{ and } \|\tilde b_{ij}^{\alpha\beta}\|_{C^{\alpha}(\bar\Omega)} + \|\tilde c_{ij}^{\alpha\beta}\|_{C^{\alpha}(\bar\Omega)} \leq o(1).\label{Rate::CoefEst1}
\end{align}
Here $o(1)$ is such that $\lim_{L_k \rightarrow 0} o(1)$ $=$ $0$. In particular, by our hypotheses, for all $L_k$ sufficiently small and $t$ $\in$ $[0,1]$, the operators $\tilde\CL_t:$ $H^1_0(\Omega)$ $\rightarrow$ $H^{-1}(\Omega)$ defined by
\[
\tilde \CL_t\,S_{ij} := \CL_{Q_*}\,S_{ij} + t\,\tilde b_{ij}^{\alpha\beta}\,\nabla S_{\alpha\beta} + t\,\tilde c_{ij}^{\alpha\beta}\,S_{\alpha\beta}
\]
are injective. The stability of the Fredholm index (see e.g. \cite[Chapter IV, Theorem 5.17]{Kato}) then implies that the $\tilde \CL_t$'s are bijective. Thus, by Lemma \ref{Trick},
\[
\|S\|_{L^2(\Omega)} \leq C\,\|R\|_{(H^2(\Omega) \cap H^1_0(\Omega))^*}.
\]
Now, by \eqref{Rate::ErrEst}, Lemmas \ref{NormEst} and \ref{NormEstEz}, we have
\[
\|R\|_{(H^2(\Omega) \cap H^1_0(\Omega))^*} \leq CL_k.
\]
We thus have
\begin{equation}
\|S\|_{L^2(\Omega)} \leq C\,L_k.
\label{Rate::L2Est}
\end{equation}
On the other hand, note that, by \eqref{Rate::DeltaS},
\[
\|S\|_{H^s(\Omega)} \leq C\,\big[\|\Delta S\|_{H^{s-2}(\Omega)} + \|S\|_{H^{s-2}(\Omega)}\big] \leq \frac{1}{2}\|S\|_{H^s(\Omega)} + C\big[\|R\|_{H^{s-2}(\Omega)} + \|S\|_{H^{s-2}(\Omega)}\big]
\]
for any $s$ $>$ $0$. Thus, by \eqref{Rate::ErrEst} and Lemmas \ref{Remainder::HsNorm} and \ref{NormEstEz},
\[
\|S\|_{H^s(\Omega)} \leq C\,L \text{ for any } 0 < s < \frac{1}{2}.
\]
Estimate \eqref{est:Hsrate} follows.

To obtain the $C^j$-convergence, we use \eqref{Rate::L2Est}, \eqref{Rate::ErrEst} and apply standard elliptic estimates to \eqref{Rate::DSLinearForm}to get 
\[
\|S\|_{H^m(K)} \leq C\,L_k,
\]
for any $K$ $\Subset$ $\Omega$ and $m$ $\geq$ $1$. This completes the proof.
\eproof

\medskip
\noindent{\bf Proof of Theorem \ref{MainThm2}.} The conclusion follows immediately from Propositions \ref{Rate} and \ref{FirstApproxEqn}.
\eproof

To finish this section, we furnish the proof of Lemma \ref{Remainder::HsNorm}. For the convenience of the reader, we recall that, for $s\in [-2,-\frac{3}{2})$, $H^s(\Omega)$ can be viewed as an interpolation space between $H^{-2}(\Omega)$ and $H^{-1}(\Omega)$ (see e.g. \cite[Theorem 6.2.4]{interpolation}), 
\begin{align*}
H^s(\Omega)
	&= [H^{-1}(\Omega),H^{-2}(\Omega)]_{s+2,2}\\
	&= \Big\{ u \in H^{-2}(\Omega)+H^{-1}(\Omega) \Big| u = \sum_{i\in\mathbb{Z}}u_i \text{ with } u_i\in H^{-1}(\Omega)\\
		&\qquad\qquad \text{ and } \sum_{i\in\mathbb{Z}} 2^{-2(s+2)i}\|u_i\|_{H^{-2}(\Omega)}^2+\sum_{i\in\mathbb{Z}}2^{-2(s+1)i}\|u_i\|_{H^{-1}(\Omega)}^2 < \infty \Big\}.
\end{align*}
The $H^s$-norm is defined to be
\[
\|u\|_{H^s} = \inf_{u=\sum_{i\in\mathbb{Z}}u_i}\Big(\sum_{i\in\mathbb{Z}} 2^{-2(s+2)i}\|u_i\|_{H^{-2}(\Omega)}^2+\sum_{i\in\mathbb{Z}}2^{-2(s+1)i}\|u_i\|_{H^{-1}(\Omega)}^2\Big)^{\frac{1}{2}}.
\]

\medskip
\noindent{\bf Proof of Lemma \ref{Remainder::HsNorm}.} As in the proof of Lemma \ref{NormEst}, it suffices to consider the case where $\Omega=[0,1]^3$ and ``$d(x)=x_n$''. We show that  if $|R^\sharp(x)|\le Ce^{-\frac{Cx_n}{\sqrt{L}}}$ then (\ref{RinX}) holds.  
 
To this end we consider first a sequence $(y_i)_{i\in\mathbb{N}}$ with $y_i\stackrel{def}{=}2^{-i}$.  We take $\chi\in C_c^\infty\big([0,1]^2\times [-1,1]\big)$ with $\chi\equiv 1$ on $[0,1]^2\times [-\frac{1}{2},\frac{1}{2}]$ and $\chi\equiv 0$ on $[0,1]^2\times [-1,-\frac{3}{4}]\cup [0,1]^2\times [\frac{3}{4},1]$ such that $\frac{\partial\chi}{\partial x_n}(x',x_n)<0$ for $x_n>0$ . Define
\[
R_i(x)\stackrel{def}{=}\left\{\begin{array}{ll}
 R^\sharp(x)\Big(\chi(x',x_n2^{-i})-\chi(x',x_n2^{-i+1})\Big),& \,i\in\mathbb{Z}, i\le 0\\
 0,\, i\in\mathbb{Z},\, i>0
 \end{array}\right.
\]

One can easily check that $R^\sharp(x)=\sum_{i\in\mathbb{Z}}R_i(x),\forall x\in [0,1]^2\times (0,1]$. Moreover we have, for $i\ge 0$, that the support of $R_{-i}$, $\supp R_{-i}$, is a subset of $\Omega_i := [0,1]^2\times [y_{i+2},y_i]$ and $|R_{-i}(x)|\le |R(x)|$ for all $ x \in \supp R_{-i}$.

For $\varphi\in C_c^\infty(\Omega_i)$ and $i\ge 0$ we have
\begin{align}
\Big|\int_{\Omega_i} R_{-i}(x) \varphi(x)\,dx\Big|
	&= \Big|\int_{\Omega_i}\Big(\int_{x_n}^{y_i}\big(\int_z^{y_i} R_{-i}(x',w_n)\,dw_n\big)\,dz\Big)\partial^2_{x_n}\varphi(x',x_n)\,dx'dx_n\Big |\nonumber\\
	&\le \bigg(\int_{\Omega_i}\Big(\int_{x_n}^{y_i}\big(\int_z^{y_i} |R_{-i}(x',w_n)|\,dw_n\big)\,dz\Big)^2\,dx\bigg)^{\frac{1}{2}}\Big(\int_{\Omega_i}|\Delta\varphi|^2\,dx\Big)^{\frac{1}{2}}\nonumber\\
	&\leq C\,\bigg(\int_{\Omega_i}\Big(\int_{x_n}^{y_i}\big(\int_z^{y_i} e^{-\frac{C\,w_n}{\sqrt{L}}}\,dw_n\big)\,dz\Big)^2\,dx\bigg)^{\frac{1}{2}}\,\|\varphi\|_{H^2(\Omega_i)}\nonumber\\
	&\leq C\,\sqrt{L}\bigg(\int_{\Omega_i}\big(\int_{x_n}^{y_i} e^{-\frac{C\,z}{\sqrt{L}}}\,dz\big)^2\,dx\bigg)^{\frac{1}{2}}\,\|\varphi\|_{H^2(\Omega_i)}\nonumber\\
	&\leq C\,L\bigg( \int_{y_{i+2}}^{y_i} (e^{-\frac{C\,x_n}{\sqrt{L}}}-e^{-\frac{C\,y_i}{\sqrt{L}}})^2\,dx_n\bigg)^{\frac{1}{2}}\|\varphi\|_{H^2(\Omega_i)}\nonumber\\
	&\leq C\,L\Big(e^{-\frac{C\,y_{i+2}}{\sqrt{L}}}-e^{-\frac{C\,y_i}{\sqrt{L}}}\Big)\sqrt{y_i-y_{i+1}}\|\varphi\|_{H^2(\Omega_i)}.
  \label{-22longest}
 \end{align}

We have shown that $\|R_{-i}\|_{H^{-2}(\Omega_i)} \le CL\Big(e^{-\frac{Cy_{i+2}}{\sqrt{L}}}-e^{-\frac{Cy_i}{\sqrt{L}}}\Big)\sqrt{y_i-y_{i+1}}$. On the other hand, as
  $\supp R_{-i}\subset\Omega_i$, one can check that $\|R_{-i}\|_{H^{-2}(\Omega_i)}=\|R_{-i}\|_{H^{-2}(\Omega)}$. We thus have
 \begin{equation}
\|R_{-i}\|_{H^{-2}(\Omega)} \le CL\Big(e^{-\frac{Cy_{i+2}}{\sqrt{L}}}-e^{-\frac{Cy_i}{\sqrt{L}}}\Big)\sqrt{y_i-y_{i+1}},
\label{-22est}
\end{equation} 
with the constant $C$ depending only on the dimension.

We next estimate $\|R_{-i}\|_{H^{-1}(\Omega)}$. For $\varphi\in C_c^\infty(\Omega_i)$ and $i\ge 0$ we have:
\begin{align*}
\Big|\int_{\Omega_i} R_{-i}(x)\varphi(x)\,dx\Big|
	&= \Big|\int_{\Omega_i} \Big(\int_{x_n}^{y_i} R_{-i}(x', w_n)\,dw_n\Big)\partial_{x_n}\varphi(x)\,dx\Big|\nonumber\\
	&\le \Big(\int_{\Omega_i}\big(\int_{x_n}^{y_i} |R_{-i}(x', w_n)|\,dw_n\big)^2\,dx\Big)^{\frac{1}{2}}\Big(\int_{\Omega_i} |\nabla\varphi(x)|^2\,dx\Big)^{\frac{1}{2}}\nonumber\\
	&\leq C\,\Big(\int_{\Omega_i}\big(\int_{x_n}^{y_i} e^{-\frac{C\,w_n}{\sqrt{L}}}\,dw_n\big)^2\,dx\Big)^{\frac{1}{2}}\,\|\varphi\|_{H^1_0(\Omega_i)}\nonumber\\
	&\leq C\,\sqrt{L}\Big(\int_{y_{i+2}}^{y_i}\big(e^{-\frac{C\,x_n}{\sqrt{L}}}-e^{-\frac{C\,y_i}{\sqrt{L}}}\big)^2\,dx_n\Big)^{\frac{1}{2}}\|\varphi\|_{H^1_0(\Omega_i)}\nonumber\\
	&\leq C\,{\sqrt{L}}\big(e^{-\frac{C\,y_{i+2}}{\sqrt{L}}}-e^{-\frac{C\,y_i}{\sqrt{L}}}\big)\sqrt{y_i-y_{i+1}}\|\varphi\|_{H^1_0(\Omega_i)}.
\end{align*}
On the other hand, noting that $y_i=2^{-i}$, one finds
\[
\sqrt{L}\Big(e^{-\frac{C\,y_{i+2}}{\sqrt{L}}}-e^{-\frac{C\,y_i}{\sqrt{L}}}\Big)\sqrt{y_i-y_{i+1}} \leq C\,L\Big(e^{-\frac{C'\,y_{i+2}}{\sqrt{L}}}-e^{-\frac{C'\,y_i}{\sqrt{L}}}\Big)\frac{\sqrt{2}}{\sqrt{y_i-y_{i+1}}}.
\]
Combining the last two inequalitites  we get
\begin{equation}
\|R_{-i}\|_{H^{-1}(\Omega)} \le \frac{CL\Big(e^{-\frac{C'\,y_{i+2}}{\sqrt{L}}}-e^{-\frac{C'\,y_i}{\sqrt{L}}}\Big)}{\sqrt{y_i-y_{i+1}}}.
\label{-12est}
\end{equation}

Taking (\ref{-22est}) and (\ref{-12est}) into account and noting that $s$ $\in$ $[-2,-3/2)$ and $y_i=2^{-i}$ for $i\ge 0$, we obtain:
\begin{align*}
\|R^\sharp\|_{H^{s}(\Omega)}
	&\leq \Big\{ \sum_{i\ge 0} 2^{2(s+2)i}\|R_{-i}\|_{H^{-2}}^2+\sum_{i\ge 0}2^{2(s+1)i}\|R_{-i}\|_{H^{-1}}^2 \Big\}^{1/2}\\
	&\leq CL\Big\{\sum_{i\ge 0} 2^{(2s+3)i}\big(e^{-\frac{C\,y_{i+2}}{\sqrt{L}}}-e^{-\frac{C\,y_i}{\sqrt{L}}}\big) + \sum_{i \ge 0} 2^{(2s+3)i}\,\big(e^{-\frac{C'\,y_{i+2}}{\sqrt{L}}}-e^{-\frac{C'\,y_i}{\sqrt{L}}}\big)\Big\}^{1/2}\\
	&\leq CL\Big\{\sum_{i\ge 0} \big(e^{-\frac{C\,y_{i+2}}{\sqrt{L}}}-e^{-\frac{C\,y_i}{\sqrt{L}}}\big) + \sum_{i \ge 0} \big(e^{-\frac{C'\,y_{i+2}}{\sqrt{L}}}-e^{-\frac{C'\,y_i}{\sqrt{L}}}\big)\Big\}^{1/2}\\
	&\leq CL,
\end{align*}
which proves \eqref{RinX}.
\eproof

\appendix
\section{Proof of Proposition \ref{ProjEqn} (completed)}

Recall that we need to solve for $X$ from \eqref{ProjEqn03} and \eqref{ProjEqn04}. We note that by \eqref{MysIds::Id3},
\[
K\Big[\frac{1}{s_+}\,Q^\sharp + \frac{1}{3}\,\IdMa\Big] = \frac{1}{3}\tr(K)\Big[\frac{1}{s_+}Q^\sharp + \frac{1}{3}\,\IdMa\Big],
\]
which implies
\begin{equation}
Q = K\,Q^\sharp = -\frac{1}{3}s_+\,K + \frac{1}{3}\tr(K)\,Q^\sharp + \frac{1}{9}s_+\tr(K)\,\IdMa.
\label{ProjEqn05}
\end{equation}
Hence, \eqref{ProjEqn03} is equivalent to
\begin{equation}
-\frac{1}{3}s_+(K\,X - X\,K) + \frac{1}{3}\tr(K)(Q^\sharp\,X - X\,Q^\sharp) = W.
\label{ProjEqn06}
\end{equation}
Using \eqref{ProjEqn04}, we thus get
\begin{equation}
-\frac{1}{3}s_+(K\,X - X\,K) + \frac{1}{3}\tr(K)(2Q^\sharp\,X - \frac{1}{3}s_+\,X) = W.
\label{ProjEqn06bis}
\end{equation}
Multiplying \eqref{ProjEqn06} to the left and to the right by $Q^\sharp$ and then taking difference yields
\begin{multline}
-\frac{1}{3}s_+(Q\,X + X\,Q) + \frac{1}{3}s_+(Q^\sharp\,X\,K + K\,X\,Q^\sharp)\\
	+ \frac{1}{3}\tr(K)((Q^\sharp)^2\,X + X\,(Q^\sharp)^2 - 2Q^\sharp\,X\,Q^\sharp) = Q^\sharp\,W - W\,Q^\sharp.
\label{ProjEqn07}
\end{multline}
On the other hand, by \eqref{ProjEqn04} and $(Q^\sharp)^2 - \frac{1}{3}s_+\,Q^\sharp - \frac{2}{9}s_+^2\,\IdMa$ $=$ $0$ (see Lemma \ref{DiffRepLimSurf}(iv)), we have
\[
(Q^\sharp)^2\,X + X\,(Q^\sharp)^2 = \frac{5}{9}s_+^2\,X,
\]
and
\[
Q^\sharp\,X\,K + K\,X\,Q^\sharp = \frac{s_+}{3}(K\,X + X\,K) - (Q\,X + X\,Q).
\]
Substituting these relations into \eqref{ProjEqn07} and noting that $Q^\sharp\,X\,Q^\sharp$ $=$ $-\frac{2}{9}s_+^2\,X$ (by \eqref{ProjEqn04} and Lemma \ref{DiffRepLimSurf}(iv)), we get
\begin{equation}
- \frac{2}{3}s_+(Q\,X + X\,Q) + \frac{1}{9}s_+^2(K\,X + X\,K) + \frac{1}{3}s_+^2\,\tr(K)\,X = Q^\sharp\,W - W\,Q^\sharp.
\label{ProjEqn08}
\end{equation}
Performing $\eqref{ProjEqn08} - \eqref{ProjEqn03}\times \frac{2}{3}s_+ - \eqref{ProjEqn06bis} \times \frac{1}{3}s_+$, we arrive at
\[
- \frac{4}{3}s_+\,Q\,X + \frac{2}{9}s_+^2\,K\,X + \frac{1}{3}s_+^2\,\tr(K)\,X - \frac{1}{9}s_+\,\tr(K)(2Q^\sharp\,X - \frac{1}{3}s_+\,X) = Q^\sharp\,W - W\,Q^\sharp - s_+\,W,
\]
i.e.
\[
\Big[-\frac{4}{3}s_+\,K\,Q^\sharp + \frac{2}{9}s_+^2\,K - \frac{2}{9}s_+\,\tr(K)\,Q^\sharp + \frac{10}{27}s_+^2\,\tr(K)\,\IdMa\Big]X = Q^\sharp\,W - W\,Q^\sharp - s_+\,W.
\]
It thus follows from \eqref{ProjEqn05} that
\begin{equation}
-2s_+\Big[Q - \frac{2}{9}s_+\,\tr(K)\,\IdMa\Big]X = Q^\sharp\,W - W\,Q^\sharp - s_+\,W.
\label{ProjEqn09}
\end{equation}

Note that $Q - \frac{2}{9}s_+\,\tr(K)\,\IdMa$ is singular because \eqref{ProjEqn05} implies that
\begin{equation}
\Big[Q - \frac{2}{9}s_+\,\tr(K)\,\IdMa\Big]\Big[\frac{1}{s_+}\,Q^\sharp + \frac{1}{3}\,\IdMa\Big] = 0.
\label{ProjEqn10}
\end{equation}
To solve \eqref{ProjEqn09} for $X$, consider $T$ defined by \eqref{ProjEqn::Est3}. Observe that $T$ is invertible and
\begin{equation}
T^{-1}\Big[Q - \frac{2}{9}s_+\,\tr(K)\,\IdMa\Big] = \IdMa - \beta\,T^{-1}\,\Big[\frac{1}{s_+}\,Q^\sharp + \frac{1}{3}\,\IdMa\Big].
\label{ProjEqn11}
\end{equation}
Multiplying \eqref{ProjEqn11} to the right by $X$ and using \eqref{ProjEqn09}, we get
\begin{equation}
X - \beta\,T^{-1}\,\Big[\frac{1}{s_+}\,Q^\sharp + \frac{1}{3}\,\IdMa\Big]X = -\frac{1}{2s_+}T^{-1}(Q^\sharp\,W - W\,Q^\sharp - s_+\,W).
\label{ProjEqn12}
\end{equation}
Multiplying \eqref{ProjEqn11} to the right by $\big[\frac{1}{s_+}\,Q^\sharp + \frac{1}{3}\,\IdMa\big]X$ and noting that $\big[\frac{1}{s_+}\,Q^\sharp + \frac{1}{3}\,\IdMa\big]^2$ $=$ $\big[\frac{1}{s_+}\,Q^\sharp + \frac{1}{3}\,\IdMa\big]$ and using \eqref{ProjEqn10}, we get
\begin{equation}
0 = \Big[\frac{1}{s_+}\,Q^\sharp + \frac{1}{3}\,\IdMa\Big]X - \beta\,T^{-1}\,\Big[\frac{1}{s_+}\,Q^\sharp + \frac{1}{3}\,\IdMa\Big]X.
\label{ProjEqn13}
\end{equation}
Combining \eqref{ProjEqn12} and \eqref{ProjEqn13}, we obtain
\begin{equation}
\Big[\frac{1}{s_+}\,Q^\sharp - \frac{2}{3}\,\IdMa\Big]X = \frac{1}{2s_+}T^{-1}(Q^\sharp\,W - W\,Q^\sharp - s_+\,W).
\label{ProjEqn14}
\end{equation}

On the other hand, \eqref{ProjEqn03} and \eqref{ProjEqn09} imply
\[
-2s_+\,X\Big[Q - \frac{2}{9}s_+\,\tr(K)\,\IdMa\Big] = Q^\sharp\,W - W\,Q^\sharp + s_+\,W.
\]
Arguing as before, we get
\[
X\Big[\frac{1}{s_+}\,Q^\sharp - \frac{2}{3}\,\IdMa\Big] = \frac{1}{2s_+}(Q^\sharp\,W - W\,Q^\sharp + s_+\,W)T^{-1}.
\]
Together with \eqref{ProjEqn04}, this implies that
\begin{equation}
\Big[\frac{1}{s_+}\,Q^\sharp + \frac{1}{3}\,\IdMa\Big]X = -\frac{1}{2s_+}(Q^\sharp\,W - W\,Q^\sharp + s_+\,W)T^{-1}.
\label{ProjEqn15}
\end{equation}
Summing up \eqref{ProjEqn14} and \eqref{ProjEqn15} we arrive at
\[
X = -\frac{1}{2s_+}\Big[T^{-1}(Q^\sharp\,W - W\,Q^\sharp - s_+\,W) + (Q^\sharp\,W - W\,Q^\sharp + s_+\,W)T^{-1}\Big].
\]
To rewrite this in a better form, note that \eqref{ProjEqn03} and \eqref{ProjEqn04} and implies that $Q^\sharp\,W + W\,Q^\sharp$ $=$ $\frac{1}{3}s_+\,W$. Therefore
\begin{equation}
X = -\Big[T^{-1}\big(\frac{1}{s_+}\,Q^\sharp - \frac{2}{3}\,\IdMa\big)W - W\big(\frac{1}{s_+}\,Q^\sharp - \frac{2}{3}\,\IdMa\big)T^{-1}\Big].
\label{ProjEqn16}
\end{equation}

By \eqref{ProjEqn01bis}, \eqref{ProjEqn02} and \eqref{ProjEqn16}, we conclude the proof.
\eproof


\def\cprime{$'$}

\end{document}